\documentclass[reqno]{amsart}
\usepackage{latexsym,amssymb,amsthm,amsmath}

\theoremstyle{plain}

\theoremstyle{remark}

\numberwithin{equation}{section}

\def\<{\left < }
\def\>{\right >}
\def\({\left ( }
\def\){\right )}

\def\E{\mathbb E}
\def\Cal{\mathcal }

\def\cf{ $C^{\infty}(M)$ }
\def\x{ $x: M\rightarrow {\mathbb E}^m$ }
\def\s{ $S^{n+1}$ }

\def\prop{ Proposition }

\def\onf{ $e_1 ,\ldots, e_n$ }
\def\cal{\Cal}
\def\lap{\Delta}

\def\CC{\Cal{C}}
\def\l{\lambda}
\def\demo{{\bf Proof. }}
\def\sq{$\square$ }

\begin{document}

\parskip=.1in
\scriptspace=.02in

\noindent{\huge {\bf Laplace Transformations of Submanifolds}}

\vskip.8in

\noindent {\Large{By Bang-Yen Chen and Leopold Verstraelen}}

\vskip.8in

\noindent{\Large{\bf Center for Pure and Applied Differential Geometry (PADGE)}}
\vskip.4in

\noindent{Katholieke Universiteit Leuven and Katholieke Universiteit Brussel}
\vskip3.8in

\noindent{Volume 1, 1995}

\vfill\eject

\noindent CONTENTS

\vskip.1in 
\noindent Preface \dotfill 5

\noindent
\noindent {Chapter I: Introduction}\dotfill 2
\vskip.1in

\noindent Chapter II: Submanifolds of Finite Type \dotfill 8

{\S1.} Basic formulas \dotfill8

{\S2.} Order of submanifolds \dotfill11

{\S3.} Minimal polynomial and  examples \dotfill17

{\S4.} Order, mean curvature and isoparametric hypersurfaces
\dotfill 21

{\S5.} Classification of submanifolds of finite type
 \dotfill25

{\S6.} Linearly independent and orthogonal immersions
 \dotfill29

{\S7.} Variational minimal principle\dotfill 34

\vskip.1in
\noindent Chapter III: Laplace Maps of Small Rank \dotfill 40

{\S1.} Curves in $\Bbb E^m$ \dotfill 40

{\S2.} Laplace maps of small rank \dotfill 42

{\S3.} Ruled surfaces in $\Bbb E^m$ \dotfill46
\vskip.1in

\noindent Chapter IV: Immersions with Homothetic Laplace
Transformations \dotfill50

{\S1.} Some general results \dotfill50

{\S2.} Some classification theorems \dotfill57

{\S3.} Surfaces with homothetic Laplace
transformation \dotfill67

 \vskip.1in
\noindent Chapter V: Immersions with Conformal Laplace
Transformations \dotfill74

{\S1.} Some general results \dotfill74

{{\S2.} Surfaces in $\Bbb E^m$ with conformal Laplace
transformation \dotfill80

{\S3.} Surfaces in $\Bbb E^3$ with conformal Laplace
transformation \dotfill83 \vskip.1in 

\noindent Chapter VI: Geometry of Laplace Images\dotfill89

{\S1.} Immersions whose Laplace images lie in a
cone \dotfill89

{\S2.} Laplace images of  surfaces of revolution
\dotfill96

{\S3.} Laplace images of curves \dotfill99

{\S4.} Laplace images of totally real submanifolds\dotfill
 101
\vskip.1in 

\noindent Chapter VII: Submanifolds with Harmonic Laplace
Maps and Transformations \dotfill 105

{\S1.} Submanifolds with harmonic Laplace map \dotfill
105

{\S2.} Submanifolds with harmonic mean curvature
function\dotfill 110
\vskip.1in

\noindent Chapter VIII: Relations between Laplace and Gauss Images
\dotfill116

{\S1.}  Laplace and Gauss images of submanifolds in
$E^m$\dotfill116

{\S2.} $LG$--transformation of spherical submanifolds
\dotfill 120

 \vskip.1in 
\noindent Chapter IX: Laplace Map and 2--type
Immersions\dotfill124

{\S1.} Spherical Laplace map \dotfill 124

{\S2.}  2--type immersions \dotfill 126

{\S3.} Laplace map of 2--type immersions\dotfill 132

\vskip.1in

\noindent Bibliography\dotfill135

\vfill\eject 
{.}
\vskip1.5in

\noindent {\LARGE {\bf Dedicated to}}
\vskip.2in
\noindent {\LARGE{\bf Pi-mei and Godelieve,}}
\vskip.3in

\noindent {\LARGE{\sl with our sincere thanks for their devotion to our families,
their
permanent
support to our work and their hospitality to our
visiting colleagues and their families.}}

\vfill\eject

\noindent {\bf PREFACE}

\bigskip
\noindent{After corresponding for some time on various problems on the differential geometry of submanifolds, the authors started their joint research on submanifold theory at the Michigan State University (MSU) at East Lansing during the academic year 1975-1976, when the second author, as a Fulbright-Hays Grantee, Nato Research Fellow and Researcher of the Belgian National Science Foundation, did post-doctoral work there, under the guidance of the first author, who then was a young Professor at MSU.}

We are happy to say that, ever since, we stayed in close scientific contact, and
that we continuously had opportunities over the past 20 years to effectively do
joint research at East Lansing and at Leuven-Brussel, thanks to the support of
mainly MSU, the Katholieke Universiteiten of Leuven and Brussel (KUL and KUB) -
most in particular through their Research Councils -, the NSF of Belgium and the
Fulbright Programm of the CEE between the USA and Belgium.  Moreover, we are
very glad that over the years our scientific contacts naturally extended to many
of our coworkers, of whom, in particular, we would like to mention here the
first author's colleagues at MSU and both our friends D.E. Blair and G. Ludden,
and the second author's collegues and former students at KUL and KUB, of whom we
specially would like to mention at this stage both our friends P. Verheyen, F.
Dillen and L. Vrancken.
Through the frequent visits of the first author to Leuven-Brussel, the second
author would like to emphasize that most of his $\pm$ 20 doctoral students till
now, from Belgium and several other countries, certainly were largely influenced
by the highly innovative work of the first named author, and that also many of
them enjoyed the great hospitality of the MSU-geometers and their families at
East Lansing.

One more side effect of our cooperation, we would like to mention here is the
following.  Around New Year 1986, J.M. Morvan and the authors discussed
at MSU, the organisation of international
meetings on differential geometry.  This later developed into what lately became
the yearly congresses in ``Pure and Applied Differential Geometry", the last
four being held at Leuven-Brussel.  We are very gratefull for all the work done
by many to make these events successfull, scientifically and socially, and in
particular would like to mention here J.M. Morvan, A. West and S. Carter, U.
Simon and M. Magid, and from the Center for Pure and Applied Differential
Geometry at KUL-KUB (PADGE) : F. Dillen and I. Van de Woestijne.

In the spring term of 1991, at MSU, the authors continued their joint work on
the theory of submanifolds of finite type, which was created some decade before
by the first named author, and also worked on biharmonic submanifolds.  In the
course of this work, we got engaged in the study of the Laplace transformations
of submanifolds and some related topics.  A first version of this work, just
stating part of our results on this at that time, was written down in a paper
for the book ``Differential Geometry in honor of Radu Rosca", published by the
Department of Mathematics of the KUL, dedicated to both our friend, and first
teacher in research of the second author, at the occasion of his 80th birthday.
The present monograph is an extension of this work, including also the proofs of
the results.  We are convinced that many more results could be obtained on the
topics initiated in this booklet and we hope that its reading would inspire some
other mathematicians to continue research on them.

The KUB indeed is, at this moment, but a small Flemish University at Brussel,
which employs only a small number of mathematicians and physicists, and one
chemist in its Group of Exact Sciences, mainly for teaching in their Faculty of
Economical and Applied Economical Sciences.  The second author is as such
teaching there, together with I. Van de Woestijne, the first year's course on
mathematics to the economical engineers.  The mathematicians doing research at
the KUB on submanifold theory at the moment are I. Van de Woestijne, J. Walrave,
J. Verstraelen and the second author.  As member of the group of Exact Sciences
at the KUB, the second author is very happy to state here that, despite the
rather limited possibilities in general of this university, the working
atmosphere there is extremely positive, the contacts with the colleagues and
students are truly enjoyable and all Deans of this Faculty so far as well as the
former and present Rector F. Van Hemelrijck and F. Gotzen, have always fully
supported the research activities of our Group.  In particular, the last three
meetings so far of the above mentioned series of congresses in differential
geometry, would not have been possible without their support.

The second author would like to make a general comment concerning the research
activities of PADGE, which basically is an international group of pure and
applied mathematicians, centered essentially at the Department of Mathematics of
the KUL (in its section with the same name) and at the Group of Exact Sciences
of the KUB: namely, we see our Center as sort of a prototype of effective
cooperation on research and teaching between the KUL and KUB.  Besides our
indebtedness to the authorities of the KUB mentioned before, we would also
like to thank very much A. Warrinnier, the Chairman of the Department of
Mathematics of the KUL, and A. Oosterlinck, the Director of the Group of Exact
Sciences of the KUL, for their support to PADGE since its creation.

The Group of Exact Sciences of the KUB intends to publish a series of monographs
on the work in which they are involved, and of which this is the first volume.
We hope, also in this way, to express our deep gratitude towards the authorities
of the KUB for their generous support of our work.

\vfill\eject
\noindent{\bf Chapter I: INTRODUCTION}

Let $x : M^n \rightarrow \Bbb E^m$ be an isometric immersion
of an $n$-dimensional $(n>0)$ connected Riemannian manifold
$M^n$ into a Euclidean $m$-space. Denote by $\Delta$ the
{\it Laplace operator\/} of $M^n$ with respect to the
Riemannian structure. Then $\Delta$ gives rise to a
differentiable map $L :M^{n} \rightarrow \Bbb E^m$,
called the {\it Laplace map\/}, defined by
$L(p)=(\Delta x)(p)$, for any point $p\in M^n$, where
$x$ denotes the immersion of $M^n$ into ${\mathbb E}^m$ as well as
the position vector field of the submanifold $M^n$ in
${\mathbb E}^m$. We call the image $L(M^n )$ of this map the
{\it Laplace image\/},   and  we
call the transformation ${\Cal L} :M^n \rightarrow
L(M^n )$ from $M^n$ onto its Laplace image
$L(M^n )$ via $\Delta$
 the {\it Laplace transformation\/} of the immersion  $x :
M^n \rightarrow \Bbb E^m$ or of the submanifold $M^n$ in ${\mathbb E}^m$.
Similar definitions of course can be given for
submanifolds in pseudo-Euclidean spaces exactly in the
same way. Many of the results mentioned in this article
can be easily extended to submanifolds in pseudo-Euclidean
spaces. \vskip.05in

Let $H$ denote the mean curvature vector field of the
immersion $x : M^n \rightarrow {\mathbb E}^m$ or the submanifold
$M^n$. As is  well-known, $H$ is a natural and
canonically determined normal vector field on $M^n$ in
${\mathbb E}^m$. Its length $\alpha = || H ||$ is called the {\it
mean curvature function\/} of $x$ or of $M^n$, or simply,
their mean curvature. Much of the geometry of the
submanifold $M^n$ is determined by the properties of $H$
and $\alpha$. One has the following fundamental
 formula of Beltrami: 
 $$\Delta x = -nH;\leqno(1.1)$$ 
 (for pseudo-Riemannian submanifolds in pseudo-Euclidean ambient spaces,
see [C5,
C7]). From the Beltrami formula, it follows that the
Laplace image $L(M^{n})$ of a submanifold $M^n$ in ${\mathbb E}^m$ is
obtained by first parallel translating the mean curvature
vector field $H$ of $M^n$ to a vector field with the
origin 0 as its initial point, and by then performing a
homothety in ${\mathbb E}^m$ with center 0 and factor $-n$. Of
course, the geometries of the Laplace image $L(M^n )$
and of the submanifold in ${\mathbb E}^m$ generated by the mean
curvature vector field $H$ of $M^n$ when centered at the
origin 0, which we could call the {\it $H$-image\/} of
$M^n$ in ${\mathbb E}^m$, are basically the same.

From the Beltrami formula, one knows, for instance, that
the immersion $x$ is minimal, or that $M^n$ is a minimal
submanifold in ${\mathbb E}^m$, if and only if their Laplace image
$L(M^n )$ is the point 0 (namely the origin of the
ambient space ${\mathbb E}^m$). Also, it is clear that  $x$ or
$M^n$ have non-vanishing constant mean curvature function
$\alpha$ if and only if their Laplace image $L(M^n )$ is
contained in a hypersphere of $E^m$ centered at 0. In
the same way, for Riemannian or pseudo-Riemannian
submanifolds in, for instance, a Minkowski space-time, it
is clear that they are pseudo-minimal in the sense of R.
Rosca, {\it i.e.,\/} $H\not= 0$ but $|| H ||=0$, if and
only if their Laplace image lies on the light-cone minus
the origin. 

 \vskip.05in

Our main purpose of this article is to initiate
the study of  the following geometric problem:

\vskip.05in {\it ``To what extent do the properties of the
Laplace transformation and/or the Laplace image of the
 immersion $x : M^n \rightarrow E^m$
determine the immersion?''}

\vskip.05in
\noindent and we will report here on  the results of
our joint work on this problem: Annoucements of some of 
these results were made in the papers [CV1,CV2].

Here, we would like to mention that the formula of
Beltrami plays a rather important role in the theory of
submanifolds of finite type, which was originated by the
first author (see, for instance, [C3-C7]). And the
main topic of our present study actually came up when
studying  questions on finite type submanifolds and
biharmonic submanifolds; the latter one is in fact one of
the off-springs of the theory of finite type.

Following up on some studies initiated in this work, further
research has been done in the mean time by various people. For
instance, following up on our study of submanifolds with
harmonic mean curvature function in Chapter VII, further
research has been done by e.g. also O. Garay, M. Barros and G.
Zafindratafa, and also, this theory now fits in as a special
case in the variational theory of $k$-minimal submanifolds, as
studied first by the authors together with F. Dillen and L.
Vrancken, and later also by D. E. Blair, M. Petrovic and J.
Verstraelen.

The contents of this monograph is as follows. In Chapter II,
we recall some basic facts and formulas on submanifolds in
general, as far as useful in the work later on. In particular,
this chapter offers as introduction to the theory of
submanifolds of finite type and several related notions. The
main purpose of Chapter III is to study submanifolds whose
Lapace maps have small rank, or more precisely, have a
constant rank smaller their dimension. In particular, we
consider submanifolds whose Laplace maps have constant rank 1.
Also, the Laplace maps of ruled surfaces is studied. In
Chapter IV, we study submanifolds for which the Laplace
transformations is homothetic. Among others, we do so in
relation with the theory of finite type. Also, classification
results are obtained for surfaces and hypersurfaces with
homothetic Laplace transformations. Chapter V deals with
submanifolds with conformal Laplace transformations.
Characterizations are given for such submanifolds to be
minimal in hyperspheres and for their Laplace images to be
minimal in hyperspheres centered at the origin. Also, in this
context, conformally flat hypersurfaces are studied. Then, in
particular, surfaces with conformal Laplace transformations
are studied, revealing among others relations with the
property to be biharmonic. In Chapter VI, we study the
geometry of Laplace images of submanifolds. As such, we study,
e.g. surfaces whose Laplace images lie in a plane, a sphere, a
cylinder or a cone, and also curves whose Laplace images lie
in a line or a circle. Finally, for totally real surfaces, we
study when their Laplace images are also totally real.
In Chapter VII, we point out the relation between submanifolds
having harmonic Laplace maps and biharmonic submanifolds. In
particular, an old result of the first author, according to
which every biharmonic surface in $\Bbb E^3$ is minimal, is
also included here. Then we study curves and submanifolds with
parallel mean curvature vector field having harmonic Laplace
transformations. Finally, we initiate the study of the
(no-compact) submanifolds with harmonic mean curvature
functions, and, in particular, classify the flat surfaces in
$\Bbb E^3$ with harmonic mean curvature functions. In Chapter
VIII, we study curves, surfaces and hypersurfaces whose
Laplace-Gauss-transformations, for short $LG$-transformations,
are conformal or, in particular, homothetic. Here, by the
$LG$-transformation of a submanifold $M^n$ in $\Bbb E^m$, we
mean the natural map $LG:L(M^n)\rightarrow G(M^n)$ from the
Laplace image $L(M^n)$ of $M^n$ to the Gauss image $G(M^n)$,
i.e. the image of the Gauss map of $M^n$. Lastly, in Chapter
IX, we study the following, first, we study the problem when
the spherical Laplace map of a submanifold $M^n$ in $\Bbb E^m$
with nonzero constant mean curvature $\alpha$ is harmonic;
here, the spherical Laplace map $L_S$ is the natural function
$L_S:M^n\rightarrow S^{m-1}(n\alpha)\subset \Bbb E^m$,
observing that for such submanifolds the Laplace image $L(M^n)$
is always contained in a hypersphere of $\Bbb E^m$ with radius
$n\alpha$ and centered at the origin. Then we study mostly
2-type submanifolds in relation with their Laplace map, in the
course of which we introduce the notions of conjugate and dual
2-type submanifolds.  In this work we follow the notations of
[C3,C4] closely.

This work was started when the second author was visiting the
Michigan State University at East Lansing in 1991. We had
further discussions on it when the first author was visiting
the Katholieke Universiteit of Leuven and of Brussel in 1991,
1993 and 1994. Also the second author worked on it when
visiting our colleagues P. Embrechts and K. Voss at the ETH
Z\"urich in 1994. We would like to thanks the ETH Z\"urich and
our own universities for their support making this work
possible. Moreover, we would like to thank the Reseach Council
of the Katholieke Universiteit Brussel for its support to the
mathematicians of its Group of Exact Sciences in general, and,
in particular, for publishing this monograph. (March 4, 1995)

\eject\vfill

\noindent {\bf Chapter II: SUBMANIFOLDS OF FINITE  TYPE}
\vskip.1in

The main purpose of this chapter is to review some results on submanifolds of
finite type for later use. For the general references, see for instances,
[C5,C6,C18].
\vskip.2in

\noindent  {\bf \S1. Basic formulas.}
\vskip.1in

Throughout this monograph,  submanifolds are assumed to
be connected, smooth, and of positive
dimension unless mentioned otherwise.

Let $x:M\rightarrow  {\mathbb E}^m$  be an isometric immersion
of an n-dimensional Riemannian manifold $M$ into
$E^m$. Let $\nabla$ and ${\tilde \nabla}$ denote the
Levi-Civita connections of $M$ and  $E^m$,
respectively. For any vector fields $X, Y$ tangent to
$M$, the Gauss formula is given by
$${\tilde \nabla}_{X}Y = \nabla_{X}Y +
h(X,Y),\leqno(2.1)$$
where $h$ denotes the second fundamental form of $M$ in 
${\mathbb E}^m$. The mean curvature vector $H$ is given by $H =
{1\over n}{\hbox{trace}\, h}, n=\dim \,M$. The length of the mean curvature
vector is called the {\it mean curvature}.

Let $A$ denote the Weingarten map.
Then $h$ and $A$ are related by $$<h(X,Y),\xi>=<A_{\xi}X,Y>
\leqno(2.2)$$ for $X,Y$ tangent to $M$ and $\xi$ normal to
$M$, where $<\,,\,>$ is the inner product in  ${\mathbb E}^m$.
Let $D$ denote the normal connection of $M$ in  ${\mathbb E}^m$.
Then the Weingarten formula is given by
$${\tilde \nabla}_{X}\xi =
-A_{\xi}X+D_{X}\xi.\leqno(2.3)$$ Let $R$ denote the
Riemannian curvature tensor of $M$. The equation of Gauss
is then given by
$$<R(X,Y)Z,W>=<h(X,W),h(Y,Z)>-<h(X,Z),h(Y,W)>.\leqno(2.4)$$
For the second fundamental form $h$, the
covariant derivative of $h$ is defined by $$({\bar
\nabla}_{X}h)(Y,Z)=D_{X}h(Y,Z)-h(\nabla_{X}Y,Z)-h(Y,\nabla_{X}Z).
\leqno(2.5)$$ The equation of Codazzi is given by
$$({\bar \nabla}_{X}h)(Y,Z)=({\bar
\nabla}_{Y}h)(X,Z).\leqno(2.6)$$
If we denote by $R^D$ the curvature tensor associated
with the normal connection $D$ of $M$ in ${\mathbb E}^m$, then
the equation of Ricci is given by
$$<R^{D}(X,Y)\xi,\eta>=<[A_{\xi},A_{\eta}]X,Y>\leqno(2.7)$$
for $X,Y$ tangent to $M$ and $\xi,\eta$ normal to $M$.
Since $M$ is Riemannian, one may choose
orthonormal local frame fields
$\{e_{1},\ldots,e_{n},e_{n+1},\ldots,e_{m}\}$ on $M$ such
that $e_{1},\ldots,e_{n}$ are tangent to $M$ and $e_{n+1},\ldots,
e_{m}$ are normal to M. Let
$\omega^{1},\ldots,\omega^{n}$ be the fields of dual
1-forms of $e_{1},\ldots,e_n$. The connection form
($\omega_{B}^A$) is then given by
$${\tilde \nabla}e_{A}= \sum_{B=1}^{m}
\omega_{A}^{B}e_{B},\,\,\,\,\,\omega_{A}^{B}=-
\omega_{B}^{A},\,\,\,A,B,C=1,\ldots,m.
\leqno(2.8)$$
From (2.8) we find
$$d\omega_{A}^{B}=-\sum_{C=1}^{m}\omega_{A}^{C}\wedge
\omega_{C}^{B}.\leqno(2.9)$$
Put $$h=\sum_{i,j,r}
h_{ij}^{r}\omega^{i}\omega^{j}e_{r},\,\,\,\,i,j=1,\ldots,n,\,\,\,
r=n+1,\ldots,s.\leqno(2.10)$$
Then we have
$$\omega_{i}^{r}=\sum_{j} h_{ij}^{r}\omega^{j},
\,\,\,<A_{r}e_{i},e_{j}>=\epsilon_{r}h_{ij}^{r},
\leqno(2.11)$$
where $A_{r}=A_{e_{r}}.$

\def\sq{$\square$}

An $n$--dimensional submanifold $M$ of a Riemannian manifold $\tilde M$ is called a
{\it totally umbilical submanifold\/} if the Weingarten map $A_\xi$ is proportional
to the identity map for any normal vector $\xi$. It is known that 
 totally umbilical submanifolds of a Euclidean $E^m$ are open 
portions of $n$--dimensional affine subspaces and open portions of
hyperspheres of $(n+1)$--dimensional affine subspaces of ${\mathbb E}^m$. 
An $n$--dimensional submanifold $M$ of a Riemannian manifold
 $\tilde M$ is called a {\it pseudo--umbilical submanifold\/} if 
 mean curvature vector $H$ is nowhere zero and the Weingarten map
$A_H$ at $H$ is proportional to the identity map.

Here, we mention two results concerning pseudo--umbilical submanifolds for later use.

{\bf Theorem 1.1.} {\rm [CY]} {\it Let $x:M\rightarrow {\mathbb E}^m$ be an
isometric immersion of an $n$--dimensional Riemannian manifold $M$ into
$E^m$. Then $M$ is a pseudo--umbilical submanifold with parallel mean
curvature vector if and only if $M$ is immersed as a minimal submanifold
in a hypersphere of ${\mathbb E}^m$}. \sq 

{\bf Theorem 1.2.} {\rm [C1]} {\it Let $x:M\rightarrow
{\mathbb E}^{n+2}$ be an isometric immersion of an $n$--dimensional
Riemannian manifold $M$ into ${\mathbb E}^{n+2}$. Then $M$ is a
pseudo--umbilical submanifold with constant mean curvature if
and only if $M$ is immersed as a minimal submanifold in a
hypersphere of ${\mathbb E}^m$.} \sq 

We need the following formula for later use, too.

{\bf Theorem 1.3.} {\rm [C5,C6]} {\it Let $M$ be an
$n$--dimensional submanifold of ${\mathbb E}^m$. Then we have
$$\Delta H=\Delta^D H+||A_{n+1}||^2H+a(H)+ \hbox{trace}(\overline{
\nabla} A_H),$$
where $$\overline{ \nabla} A_H=\nabla A_H+A_{DH},\leqno(2.12)$$
$$\hbox{trace}(\overline{\nabla} A_H)=\sum_{i=1}^n \{A_{D_{e_i}H}e_i
+(\nabla_{e_i}A_H)e_i\}.\leqno(2.13)$$
 $$a(H)=\sum_{r=n+2}^m \hbox{trace}(A_HA_r)e_r.$$
 $$\Delta H=\Delta^{D}H + \sum_{i=1}^{n}
h(A_{H}e_{i},e_{i}) +  {{n}\over {2}}grad<H,H> + 2\,
trace\,A_{DH},\leqno(1.24)$$ where
 $\Delta^D$ is the Laplacian operator
associated with the normal connection $D$.}  \sq

If $M$ is a spherical submanifold, we have the following
[C5,C6,C8].

{\bf Proposition 1.4.} {\it Let $M$ be an n--dimensional
submanifold of a hypersphere $S^{m-1}(r)$ of radius r and
centered at the origin of ${\mathbb E}^m$. Then we have
$$\Delta H=\Delta^{\bar D} {\bar H} +(n+ ||A_\xi ||^2)\bar H +{\Cal
A}(\bar H) +\hbox{trace}(\overline{\nabla} A_H)-{{n\alpha^2}\over
{r^2}}x,\leqno(2.15)$$ where $\bar H$ is the mean curvature vector of
$M$ in $S^{m-1}(r)$,  $\xi$  the unit vector in the direction of
$\bar H$, and $\bar D$, ${\Cal A}(\bar H)$
the normal connection and  the allied mean curvature vector of $M$ in
$S^{m-1}(r)$, respectively.} \sq
\vskip.2in

\noindent  {\bf \S2. Order of submanifolds.}
\vskip.1in

An algebraic manifold or an algebraic variety is defined by algebraic
equations. Thus, one may define the notion of degree of an algebraic
manifold by its algebraic structure (which can also be defined by
using homology). The concept of degree is both important and fundamental
in algebraic geometry. On the other hand, one cannot talk about
the degree of an arbitrary submanifold in a Euclidean m--space ${\mathbb E}^m$.
In this section, we
will use the induced Riemannian structure on a submanifold $M$ of
${\mathbb E}^m$ to introduce two well--defined numbers $p, q$ associated with
the submanifold $M$; more precisely, $p$ is a positive integer and
$q$ is either $+\infty$ or an integer $\geq p$. We call the pair
$[p,q]$ {\it the order of the submanifold M\/}; moreover, $p$ is
called the {\it lower order\/} and $q$ the {\it upper order} of the
submanifold $M$.  The submanifold $M$ is said to be {\it of finite
type\/} if the upper order $q$ is finite.  And $M$ is of {\it
infinite type\/} if the upper order $q$ is $\infty$.

Although  it remains possible to define the notion of submanifolds of
finite type and the related notions of order,... for those
submanifolds, for simplicity, we  define them here only for compact 
submanifolds.

Let $(M,g)$ be a compact Riemannian manifold, denote the Levi--Civita
connection of $(M,g)$ by $\nabla$ and let $\Delta=-\hbox{trace}\,
\nabla^2$ be the Laplacian of $(M,g)$ acting as a differential
operator on $C^{\infty}(M)\subset L^2(M,\mu_g)$. It is well--known
that the eigenvalues of $\Delta$ form a discrete infinite sequence:

$$0=\lambda_0<\lambda_1<\lambda_2<\ldots \nearrow\infty.$$
Let $V_k=\{f\in C^{\infty}(M):\Delta f=\lambda_k f\}$ be the 
eigenspace of $\Delta$ with eigenvalue $\lambda_k$. Then $V_k$ is
finite--dimensional. We define an inner product $(\;,\;)$ on 
\cf by
$$(f,h)=\int_M fh dV$$
where $dV$ is the volume element of $(M,g)$. 

Then $\sum_{k=0}^{\infty} V_k$ is dense in \cf (in $L^2$--sense).
 Denoting
by $\hat{\oplus} V_k$ the completion of $\sum V_k$, we have
$$C^{\infty}(M)=\hat\oplus_k V_k.$$

For each function $f\in$ \cf let $f_t$ be the projection of $f$ onto
the subspace $V_t$ $(t=0,1,2,\ldots).$ Then we have the following
spectral decomposition:
$$f=\sum_{t=0}^{\infty} f_f\quad (\hbox{in $L^2$-sense}).$$

Becuase $V_0$ is 1--dimensional, for any non--constant function
$f\in $\cf there is a positive integer $p\geq 1$ such that
$f_p\not=0$ and
$$f-f_0=\sum_{t\geq p}f_t,$$
where $f_0\in V_0$ is a constant. If there are infinite $f_t$'s which
are nonzero, we put $q=\infty$, Otherwise, there is an integer $q$,
$q\geq p$, such that $f_q\not=0$ and
$$f-f_0=\sum_{t=p}^q f_t.\leqno(2.1)$$
If we allow $q$ to be $\infty$, we have the decomposition (2.1) in 
general.

The set 
$$T(f)=\{ t\in N_0:f_t \not=0\}\leqno(2.2)$$
is called {\it the type of f.} A function $f$ in \cf is said to be
of {\it finite type} if $T(f)$ is a finite set, {\it i.e.,} if its 
spectral decomposition contains only finitely many non--zero terms.
Otherwise $f$ is said to be {\it of infinite type.} The
smallest element in T(f) is called 
the {\it lower order of f} (notation: l.o.($f$)) and the supremum
of $T(f)$ is called the {\it upper order of f} (notation: u.o.($f$)).
$f$ is said to be of {\it k--type} if $T(f) $ contains exactly
$k$ elements. 

For an isometric immersion $x: M\rightarrow {\mathbb E}^m$ of a compact
Riemannian manifold $M$ into a Euclidean $m$--space, we put
$$x=(x_1,\ldots,x_m)$$
where $x_A$ is the $A$--th Euclidean coordinate function of $M$ in
${\mathbb E}^m$. For each $x_A$ we have
$$x_A-(x_A)_0=\sum_{t=p_A}^{q_A}(x_A)_t,\quad A=1,\ldots,m.$$
For the isometric immersion $x:M\rightarrow {\mathbb E}^m$, we put
$$p=\inf_{A} \,\{p_A \},\quad q=\sup_A\, \{q_A\}$$
where $A$ ranges among all $A$ such that $x_A-(x_A)_0\not=0.$ It 
is easy to see that $p$ and $q$ are well--defined geometric invariants such
that $p$ is a positive integer and $q$ is either $\infty$ or an 
integer $\geq p$. By using the above notation, we have the following
spectral decomposition of $x$ in vector form:
$$x=x_0+\sum_{t=p}^q x_t.\leqno(2.3)$$
We define $T(x)$ by
$$T(x)=\{t\in N_0: x_t\not=0\}.$$
The immersion $x$ or the submanifold $M$ is said to be of {\it k--type}
if $T(x)$ contains exactly $k$ elements. Similarly we can define
the lower order of $x$ and the upper order of $x$.
The immersion $x$ is said to be {\it of finite type\/} if $q$, the upper
order of $x$, if finite; and $x$ is {\it of infinite type\/} if $q$ is
$\infty$.

{\bf Remark 2.1.}  It is easy to see
that $T(x)=T(y)$ if $x$ and $y$ are congruent isometric immersions of
$(M,g)$ in ${\mathbb E}^m$, so that all these notions do in fact have geometric
significance. \sq 

{\bf Remark 2.2} Let $x:M\rightarrow {\mathbb E}^m$ be an isometric
immersion and $\bar x:M\rightarrow {\mathbb E}^m\subset {\mathbb E}^{\bar m}$. Then we have
$T(x) =T(\bar x)$. \sq

We give the following lemmas for later use.

{\bf Lemma 2.1.} {\it  Let $x:M\rightarrow {\mathbb E}^m$ be an isometric 
immersion of a compact Riemannian manifold M into $E^m$. Then $x_0$ 
is the center of mass of M in ${\mathbb E}^m$.}

 Consider the spectral decomposition:
$$x=x_0+\sum_{t=p}^q x_t,\quad \Delta x_t=\lambda_tx_t.$$
By Hopf' lemma, we have 
$$\int_M x_tdV={1\over {\lambda_t}}\int_M\Delta x_tdV=0,
\quad t\not=0.$$
Since $x_0$ is a constant vector in $E^m$, we obtain that
$$x_0={1\over{vol(M)}}\int_MxdV.$$
This shows that $x_0$ is the center of mass of $M$ in ${\mathbb E}^m$. \sq

For two ${\mathbb E}^m$-valued functions $v,w$ on $M$, we define the inner
product $v, w$ by
$$(v,w)=\int_M\left< v,w\right>dV$$
where $\left<v,w\right>$ denotes the Euclidean inner product of
$v,w$.

We then have the following.

{\bf Lemma 2.2.}  {\it Let $x:M\rightarrow {\mathbb E}^m$ be an isometric
immersion of a compact Riemannian manifold M into ${\mathbb E}^m$. Then we have
$$(x_t,x_s)=0,\quad t\not=s,$$
where $x_t$ is the t--th component of x with respect to the spectral
decomposion of x.}

{\bf Proof. } Since $\Delta$ is self--adjoint, we have
$$\lambda_t(x_t,x_s)=(\Delta x_t,x_s)=(x_t,\Delta x_s)=\lambda_s
(x_t,x_s).$$
Because $\lambda_t\not=\lambda_s$, we obtain the lemma. \sq

For  general Riemannian
manifolds, in general one cannot make a spectral decomposition of
a function. However, it remains possible to define the notion of a 
function of finite type and the related notions of order,.... for
those functions. A function $f$ in general is said to be {\it of
finite type\/} if it is a finite sum of eigenfunctions of the
Laplacian. More precisely, a function of finite type can be written
as
$$f=f_0+\sum_{i=1}^k f_i$$
where $\Delta f_0=0,\; f_i$ is a nonzero eigenfunction of
 $\Delta$ with
nonzero eigenvalue $\lambda_{t_i}$ and all $\lambda_{t_i}$ are mutually
distinct. We define $T(f)$ to be the set $\{\lambda_{t_1},\ldots,
\lambda_{t_k}\}$, u.o.$(f)=\sup\;T(f)$, l.o.$(f)=\inf \;T(f)$.
Similar notions can be defined for an isometric immersion from
a general Riemannian manifold into a Euclidean space.

Here we mention the following two well-known results.

{\bf Proposition 2.3.} {\rm (E. Beltrami)} 
 {\it Let \x be an isometric immersion. Then we have
$$\Delta x=-nH,\leqno (2.6)$$
where H denotes the mean curvature vector of x.}

{\bf Proof. } Let $h, \nabla,\tilde\nabla$ be the second fundamental
form of $x$, the Levi--Civita connection of $M$ and the Levi--Civita
connection of ${\mathbb E}^m$, respectively. Denote by $e_1,\ldots,e_n$ an
orthonormal local frame fields of $M$. Then we have
$$\Delta x=-\sum_{i=1}^n \{e_ie_ix-(\nabla_{e_i}e_i)x\}
=-\sum_{i=1}^n h(e_i,e_i)=-nH. \;\;\square$$

{\bf Theorem 2.4.} {\rm [Ta1]}  {\it Let \x be an isometric immersion. If $\Delta x=
\lambda x,\;\lambda \not=0$, then
\begin{itemize}
\item[(i)] $\lambda>0$;
\item[(ii)] $x(M)\subset S^{m-1}_0(r)$, where $S^{m-1}_0(r)$
is a hypersphere of ${\mathbb E}^m$ centered at the origin $0$ and with radius
$r=\sqrt{n/\lambda}$;
\item[(iii)]  $x(M)\subset S^{m-1}_0(r)$ is a minimal immersion;
moreover, if $x(M)\subset S^{m-1}_0(r)$ is a minimal immersion,
then $\Delta x=({n\over {r^2}})x.$
\end{itemize}}

{\bf Proof. } If $\Delta x=\lambda x,\,\lambda\not=0$, then we have
$H=-({{\lambda}\over n})x.$ Let $X$ be a vector field tangent to $M$. We 
have
$\left<x,X\right>=0.$ Thus, $X\left<x,x\right>=2\left<x,X
\right>=0.$ Therefore, $\left<x,x\right>$ is constant on $M$.
This proves that $M$ is immersed into a hypersphere $S^{m-1}(r)$
of $E^m$ centered at the origin with radius $r$. Let $h $, $h'$
and  $\tilde h$ be the second fundamental forms of $M$ in
${\mathbb E}^m,$ $M$ in  $S^{m-1}(r)$, and $S^{m-1}(r)$ in ${\mathbb E}^m$,
respectively. Then we have $$h(X,Y)=h '(X,Y)+\tilde h(X,Y)$$ for
$X,Y$ tangent to $M$. Thus, the mean curvature vector
$H=H'-{1\over  {r^2}}x$. Since $x$ is perpendicular to
$S^{m-1}(r)$, and $H$ is parallel to $x$, we have $H'=0$. Thus
$M$ is minimal in $S^{m-1}(r)$. Furthermore, we have $$\Delta
x=-nH={n\over {r^2}}x.$$ Therefore, $\lambda={n\over {r^2}}.$
This proves statements (i)--(iii). The remaining part is easy to
verify. \sq

{\bf Remark 2.3.} The pseudo-Riemannian version of Theorem 2.4
was obtained in [C9]. \sq

In terms of finite type submanifolds, Proposition 2.3 and Theorem 2.4
give the following  result.

{\bf Proposition 2.5} {\it Let \x be an isometric immersion. Then $x$
is of 1--type if and only if either $M$ is a minimal submanifold of 
${\mathbb E}^m$ or $M$ is a minimal submanifold of a hypersphere of ${\mathbb E}^m$.}
\sq

In particular, Proposition 2.5 shows that a linear subspace and
 a hypersphere of a Euclidean space are of 1--type. 

We give some  more examples of finite type submanifolds.

{\bf Example 2.1.} (Product submanifolds) If $M$ is
 a finite type submanifold of ${\mathbb E}^m$ and
$N$ a finite type submanifold of ${\mathbb E}^N$, then the product submanifold
$M\times N$ in ${\mathbb E}^{m+N}$ is of finite type.

In particular, the product of two plane circles $S^1(a)\times S^1(b)$
is of finite type in ${\mathbb E}^4$. In fact, it is of 1--type if $a=b$ and
is of 2--type if $a\not=b$. Also a circular cylinder $R\times S^1$
is a finite type surface in ${\mathbb E}^3$, in fact, a 2--type surface in
${\mathbb E}^3$. \sq

{\bf Example 2.2.} (Diagonal immersions) Let  $x_i:M\rightarrow
{\mathbb E}^{n_i}$, $i=1,\ldots,k$ be k isometric immersions from a Riemannian
manifold $M$ into ${\mathbb E}^{n_i}$ respectively. For any $k$ real numbers
$c_1,\ldots,c_k$
 with $c_1^2+\ldots+c_k^2=1$, the immersion:
$$\tilde x=(c_1x_1,\ldots,c_kx_k):M\rightarrow {\mathbb E}^{n_1+\cdots+n_k}$$
is an isometric immersion from $M$ into ${\mathbb E}^{n_1+\cdots+n_k}$. Such an
immersion is called a {\it diagonal immersion.}  If $x_i$ are all of finite
type, then also $\tilde x$ is of finite type. \sq

There are ample examples of finite type submanifolds which are not
of the above kinds. Here we give some such  examples.

{\bf Example 2.3.} (2--type curves in ${\mathbb E}^3$) For each positive number
$\epsilon$ we put
$$\gamma_\epsilon(s)={{12}\over{\epsilon^2+36}}(\epsilon\sin s,
-{{\epsilon^2}\over{12}}\cos s+\cos 3s,-{{\epsilon^2}\over{12}}
\sin s+\sin 3s).$$
Then $\gamma_\epsilon$ is of 2--type. According to a result of [CDV],
up to homothetic transformations of ${\mathbb E}^3$,
a 2--type curve in $E^3$ is either a right circular helix or it is
congruent to the curve $\gamma_\epsilon$ for some positive number
$\epsilon$. \sq

{\bf Example 2.4.} (A flat torus in ${\mathbb E}^6$) Consider the flat torus
$$T_{ab}=S^1(a)\times S^1(b),\quad a^2+b^2=1.$$
Let $x:T_{ab}\rightarrow {\mathbb E}^6$ be defined by
$$x=x(s,t)=(a\sin s,b\sin s\sin {t\over b},b\sin s\cos{t\over b},
a\cos s,b\cos s\sin{t\over b},b\cos s\cos{t\over b}).$$

By a direct compution, we can see that
 $T_{ab}$ is of 2--type in ${\mathbb E}^6$. \sq

This example was first given by N. Ejiri [Ej1] in answering an open
question of Weiner.

{\bf Example 2.5.} (Spheres in ${\mathbb E}^m$) Let $S^n$ be the unit hypersphere
of ${\mathbb E}^{n+1}$ defined by $y_1^2+y_2^2+\ldots+y_{n+1}^2=1$, where
$(y_1,\ldots,y_{n+1})$ is a Eucldean coordiante system of ${\mathbb E}^{n+1}$.
Let $$x=(x_1,\ldots,x_m):S^n\rightarrow {\mathbb E}^m$$
be an isometric immersion of $S^n$ into ${\mathbb E}^m$. Then the coordinate
functions
$$x_A=x_A(y_1,\ldots,y_{n+1})$$ are funtions of $y_1,\ldots,y_{n+1}$.
Since the eigenspace of $\Delta$ of $S^n$ associated with $\lambda_k$
is spanned  by harmonic homogeneous polynomials of degree $k$ on
${\mathbb E}^{n+1}$, restricted to $M$,
 the isometric immersion \x is thus of finite type if and only
if each of $x_A(y_1,\ldots,y_{n+1})$ is a polynomial in $y_1,\ldots,y_{n+1}$.
Moreover, when \x is of finite type, then the lower 
order $p$ of the immersion
$x$ is nothing but the lowest (nonzero) degree of $\{x_A(y_1,\ldots,
y_{n+1}):A=1,\ldots,m\}$ and the upper order $q$ of the immersion $x$
is the highest degree of  $\{x_A(y_1,\ldots,
y_{n+1}):A=1,\ldots,m\}$.

Therefore, in this case, the order of the immersion is nothing but the
degree. \sq

{\bf Example 2.6.} (Symmetric Spaces in $\Bbb E^m$) 

For a compact rank one symmetric space, we have the following two
theorems of [CDV1].

{\bf Theorem 2.6.}  {\it If \x is an isometric immersion
of a compact rank one symmetric space  $M$ into a Euclidean space,
then $x$ is of finite type if and only if  all geodesics of $M$ are
curves of finite type in ${\mathbb E}^m$ via $x$.} \sq

{\bf Theorem 2.7.} {\it Let \x be an isometric immersion of finite type
of a compact rank one symmetric sapce  $M$ with upper order q. Then
\begin{itemize}
\item[(i)] every geodesic is mapped to a curve of finite type
of upper order at most q; and
\item[(ii)] through each point of M, there is a geodesic
of upper order q. \sq
\end{itemize}}

Theorem 2.7 was generalized by J. Deprez in  [De1] to the following.

{\bf Theorem 2.8.} {\it Let M be the Riemannian product of a symmetric space
of compact type and a number of circles. Then an isometric immersion
\x is of finite type if and only if \x maps all geodesics to curves of
finite type.}

\vskip.2in
\noindent  {\bf \S 3. Minimal polynomial and examples.}
\vskip.1in

The following criterion of [C5,C6] for finite type submanifolds is quite
useful.

{\bf Theorem 3.1.} (Minimal Polynomial Criterion) {\it
 Let \x be an isometric immersion of a compact
Riemannian manifold M into ${\mathbb E}^m$ and H the mean curvature vector of
M in ${\mathbb E}^m$. Then 
\begin{itemize}
\item[(i)] M is of finite type if and only if there is
a nontrivial polynomial $Q(t)$ such that $Q(\Delta)H=0;$
\item [(ii)] if M is of finite type, there is a unique
monic polynomial P(t) of least degree with $P(\Delta)H=0$;
\item[(iii)] if M is of finite type, then M is of k--type
if and only if $\deg(P)=k$.
\end{itemize}
The same results holds if $H$ is replaced by $x-x_0$.}

{\bf Proof. } Let \x be an isometric immersion of a compact Riemannian
manifold $M$ into $E^m$. Consider the  spectral 
decomposition:
$$x=x_0+\sum_{t=p}^qx_t,\quad \Delta x_t=\lambda_tx_t.\leqno(3.1)$$
If $M$ is of finite type, then $q<\infty.$ Because $\Delta x=-nH$,
 we have
$$-n\Delta^iH=\sum_{t=p}^q\lambda^{i+1}_tx_t,\quad i=0,1,
\ldots.\leqno(3.2)$$
Let
$$c_1=-\sum_{t=p}^q\lambda_t,\quad c_2=\sum_{t<s}\lambda_t\lambda_s,
\ldots,c_{q-p+1}=(-1)^{q-p+1}\lambda_p\cdots\lambda_q.$$
Then by direct computation we find
$$\Delta^kH+c_1\Delta^{k-1}H+\ldots+c_kH=0,$$
where $k=q-p+1.$ Conversely, if $H$ satisfies $Q(\Delta)H=0$ for 
some monic polynomial 
$$Q(t)=t^k+c_1t^{k-1}+\ldots+c_{k-1}t+c_k,\leqno(3.3)$$
then from (2.2) we have
$$\sum_{t=1}^\infty \lambda_t(\lambda^k_t+c_1\lambda^{k-1}_t
+\ldots+c_{k-1}\lambda_t+c_k)x_t=0.\leqno(3.4)$$
For each positive integer $s$, (3.4) yields
$$\sum_{t=1}^\infty \lambda_t(\lambda^k_t+c_1\lambda^{k-1}_t
+\ldots+c_k)\int_M\<x_s,x_t\>dV=0.$$
Therefore, by applying Lemma 1.2, we obtain
$$(\lambda^k_s+c_1\lambda^{k-1}_s+\ldots+c_k)||x_s||^2=0,\leqno(3.5)$$
where $||x_s||^2=(x_s,x_s).$ If $x_s\not=0$, then (3.5) implies
$$\lambda^k_s+c_1\lambda^{k-1}_s+\ldots+c_k=0.$$
Since this equation has at most $k$ real solutions and equation 
(3.5) holds for any positive integer $s$, at most $k$ of the $x_t$'s
are nonzero. Thus the decomposition (3.1) is in fact a finite
decomposition. Consequently, the immersion $x$ is of finite type. This 
proves (i).

Suppose that $M$ is of finite type. Let $P=P(t)$ be a monic polynomial
of least degree with $P(\Delta)H=0$. If $K$ is another such polynomial,
then $\deg(P)=\deg(K).$ Since $R=P-K$ is a polynomial of smaller 
degree satisfying $R(\Delta)H=0$, we have $R=0$. This implies that
$P=K$. So (ii) is proved.

(iii) follows from the proof of (i).

The same proof applies if we replace $H$ by $x-x_0$. \sq

{\bf Definition 3.1.} The unique monic polynomial $P$ given  in (ii) of
Theorem 2.1 is called the {\it minimal polynomial\/} of the finite
type submanifold $M$. \sq

 If the submanifold $M$ is not compact, then the existence 
of a nontrivial polynomial $Q$ such that $Q(\Delta)H=0$ does not imply
that $M$ is of finite type in general. However, we have the following
three results of [CP].

{\bf Theorem 3.2.} {\it Let $\gamma$ be a curve in ${\mathbb E}^m$. If there is 
a nontrivial polynomial $Q$ of one variable such that $Q(\Delta)H=0$, 
then $\gamma$ is a curve of finite type.} \sq

{\bf Theorem 3.3.} {\it Let $x: M\rightarrow {\mathbb E}^m$ be an isometric
immersion. If there exists a polynomial $Q$ such that $Q(\Delta)H=0$,
then either M is of infinitely type or is of k--type with 
$k\leq\deg P +1$.}  \sq

{\bf Theorem 3.4.} {\it Let $x:M\rightarrow {\mathbb E}^m$ be an isometric
immersion. If there exist a vector $c\in {\mathbb E}^m$ and a polynomial
$P(t)=\prod_{i=1}^k (t-\ell_i)$ with mutually distinct 
$\ell_1,\ldots,
\ell_k$ such that $P(\Delta)(x-c)=0$, then M is of finite type.}
\sq

{\bf Remark 3.1.} It is clear that congruent submanifolds of finite type
have the same minimal polynomial. Furthermore, if two finite type
submanifolds $M$ and $N$ in $E^m$ have the same minimal polynomial,
then (1) they are of the same type, (2) the Laplacians of $M$ and $N$
have at least $k$ common eigenvalues given by the roots of the minimal
polynomial and (3) the order of $M$ and $N$ are the same whenever $M$
and $N$ are isospectral. As a consequence, the minimal polynomial provides
us some important imformation on the spectrum of a finite type submanifold.
\sq

{\bf Definition 3.2.} Let $M$ be a manifold and $G$ a closed subgroup
of  the group $I(M)$ of the isometries acting transitively on $M$ and
$\tilde M$ a Riemannian manifold with  group $I(\tilde M)$ as the group of
isometries. An immersion $f:M\rightarrow \tilde M$ of $M$ into $\tilde M$
is called {\it G--equivariant\/} if there is a homomorphism
$\zeta:G\rightarrow I(\tilde M)$ such that
$$f(a(p))=\zeta(a)f(p)$$
for $a\in G$ and $p\in M$. \sq

{\bf Theorem 3.5.} {\rm [C5]} {\it Let $M$ be a compact homogeneous Riemannian
manifold. If M is equivariantly isometrically immersed in ${\mathbb E}^m$, then M
is of finite type. Moreover, M is of k--type with $k\leq m$.}

{\bf Proof. } Let $p$ be an arbitrary point of $M$. Then the $m+1$
vectors $H,\Delta H,\ldots,$ $ \Delta^mH$ at $p$ are linearly independent.
Thus, there is a polynomial $Q(t)$ of degree $\leq m$ such that
$Q(\Delta)H=0$ at $p$. Because $M$ is equivariantly isometrically
immersed in $E^m$, $Q(\Delta)H=0$ at every point of $M$. Thus, by
Theorem 2.1, $M$ is of finite type. Moreover, because the minimal
polynomial of $M$ is of degree $\leq m$, $M$ is of $k$--type with
$k\leq m$. \sq

In [De1,Ta2], J. Deprez and T. Takahashi studied equivariant isometric immersions
of compact {\it irreducible\/} homogeneous Riemannian manifolds
and they proved the following.

{\bf Theorem 3.6.} {\it M be a compact homogeneous Riemannian 
manifold with irreducible isotropy action. Then an equivariant 
isometric immersion $x:M\rightarrow {\mathbb E}^m$ of M into ${\mathbb E}^m$ is
the diagonal immersion of some 1--type isometric immersions of
$M$.} \sq

A submanifold $M$ of a hypersphere $S^{m-1}$ of ${\mathbb E}^m$ is said to be
{\it mass-symmetric\/} if the center of gravity of $M$ coincides with
the center of the hypersphere $S^{m-1}$ in ${\mathbb E}^m$.

For mass--symmetric 2--type immersions of a topological sphere we have the
following result of M. Kotani.

{\bf Theorem 3.7.} {\rm [Ko]}  {\it Any mass--symmetric 2--type immersion of a
topological 2--sphere into a hypersphere of ${\mathbb E}^m$ is the diagonal
immersion of two 1--type immersions.} \sq

Following from this theorem of Kotani, we have the following two
consequences.

{\bf Corollary 3.8.} {\rm [Ko]} {\it If the immersion in Theorem 3.7
 is full,
then m is odd and greater than 5.} \sq

{\bf Corollary 3.9.} {\rm [Ko]} {\it If a 2--sphere admits a mass--symmetric 
2--type immersion into $S^9$, then the 2--sphere is of constant
curvature.} \sq

Although there do exist mass--symmetric 2--type surfaces in $S^3$ and in 
$S^5$, it is interesting to point out that for surfaces in $S^4$, we
 have the following non--existence theorem.

{\bf Theorem 3.10.} {\rm [BC1]} {\it There exists no compact mass--symmetric
2--type surfaces which lie fully in $S^4\subset {\mathbb E}^5$.} \sq

In view of Theorem 3.7 of Kotani, it is interesting to point out that
{\it there exist 2--type isometric immersions of an ordinary 2--sphere 
into ${\mathbb E}^{10}$
with order $\, [1,3]\,$ that are not diagonal immersions.\/} (cf. [CDV3]).

We mention the following result of [C11,C12] for later use.

{\bf Theorem 3.11} {\it Let $x:M\rightarrow {\mathbb E}^m$ be an isometric
immersion from an $n$--dimensional Riemannian manifold into ${\mathbb E}^m$. Then
the mean curvature vector field of $x$ is an eigenvector of the Laplace
opertor, i.e., $\Delta H=\lambda H$ for some constant $\lambda$,
if and only if $M$ is one of the following submanifolds:
\begin{itemize}
\item[(1)] a 1-type submanifold of ${\mathbb E}^m$;
\item[(2)] a biharmonic submanifold; or
\item[(3)] a null 2--type submanifold. 
\end{itemize}
In particular,  a surface $M$ in ${\mathbb E}^3$ satisfies the condition
$\Delta H=\lambda H$  for some constant $\lambda$ if and only if $M$
is either a minimal surface or an open portion of circular cylinder.}
\sq  

{\bf Remark 3.2.} In fact, the exact proof given in [C12] works for
pseudo-Rieman\-nian submanifolds in pseudo-Euclidean space as
well. Hence,
 the pseudo-Rieman\-nian version of Theorem 3.11 also holds
(cf. [C12,C23]). \sq

{\bf Remark 3.3.}  The classifications of submanifolds in
hyperbolic spaces and in de Sitter space-times which satisfy
condition $\Delta H=\lambda H$ are obtained in [C21] and in
[C23], respectively.

{\bf Remark 3.4.} For the classification of spherical 2-type
surfaces with constant sectional curvature, see [Mi]; and
for the classification of flat 2-type spherical Chen surfaces,
see [G1].
\vskip.2in

\noindent {\bf \S 4. Order,  mean Curvature, and isoparametric
hypersurfaces.}

\vskip.1in

In this section, we will give some relations between order,
mean curvature and type for spherical submanifolds.
For simplicity, we will do this here only for spherical
hypersurfaces.

In the following, we denote by $S^{n+1}$ the unit hypersphere of
${\mathbb E}^{n+2}$ centered at the origin.

{\bf Theorem 4.1.} {\it Let $M$ be a compact 
hypersurface of 
 $S^{n+1}$. Then $M$ has nonzero constant  mean curvature in 
\s and
constant scalar curvature if and only if $M$ is mass--symmetric and
of 2--type, unless M is a small hypersphere of $S^{n+1}.$}

{\bf Proof. } First we give the following [C8].

{\bf Lemma 4.2.}  {\it Let M be an n--dimensional
hypersurface of $S^{n+1}$. Then we have
$$\hbox{trace}(\overline{\nabla} A_H)={n\over 2}\hbox{grad}\;\alpha^2
+2\,\hbox{trace}\;A_{DH}.\leqno(4.1)$$}

{\it Proof of Lemma 4.2.} Let $\xi$ be a unit normal vector field of
$M$ in $S^{n+1}$ and $e_1,\ldots,e_n$  $n$ orthonormal
 eigenvectors
of $A_\xi$ with eigenvalues $\rho_1,\ldots,\rho_n$, respectively.
Denote by $H'$ the mean curvature vector field of $M$ in $S^{n+1}$ and
$H'=\alpha'\xi$.
Then we have
$$A_{H'}e_i=\alpha'\rho_ie_i.\leqno(4.2)$$
We put $\nabla e_j=\sum \omega^k_j e_k.$ Then by the Codazzi equation
we may obtain
$$e_j\rho_i=(\rho_i-\rho_j)\omega^j_i(e_i),\quad j\not= i.\leqno(4.3)$$
On the other hand, by applying
$$H=H'-x,\quad A_x=-I,\quad Dx=0,$$
and 
$$(\nabla_{e_i}A_{H'})e_j=(e_i\alpha')\rho_je_j+\alpha' (e_i\rho_j)e_j
+\sum \alpha'(\rho_j-\rho_k)\omega^k_j(e_i)e_k,$$
we may obtain
$$\hbox{trace}(\overline{\nabla}A_H)=\sum_i\{2(e_i\alpha')\rho_i
+\alpha'(e_i\rho_i)+\sum_k \alpha'(\rho_k-\rho_i)\omega^i_k(e_k)\}
e_i.\leqno(4.4)$$
Substituting (4.3) into (4.4) and making a direct computation, we
may obtain (4.1). \sq

Now, we really start the proof of Theorem 4.1.

Combining Lemma 3.4 with Lemma 4.2 we find
$$\Delta H=\Delta^D H+a(H)+||A_{n+1}||^2H+{n\over 2}\,\hbox{grad}\;\alpha^2
+2\,\hbox{trace}\;A_{DH}.\leqno (4.5)$$
Since $H=H'-x=\alpha'\xi-x$, we may also obtain from (4.5) the following.
$$\Delta H=(\Delta\alpha') \xi+(||A_\xi ||^2+n)H'-n|H|^2x
+{n\over 2}\,\hbox{grad}\;\alpha^2
+2\,\hbox{trace}\;A_{DH}.\leqno(4.6)$$
In particular, if $M$ is mass--symmetric and of 2--type, then we
have the following spectral decomposition:
$$x=x_p+x_q,\quad \Delta x_p=\lambda_p x_p,\quad \Delta x_q=\lambda_q
x_q.$$
This implies
$$n\Delta H=n(\lambda_p+\lambda_q)H+\lambda_p\lambda_q x.\leqno(4.7)$$
Because $H=H'-x$ and since $\xi, x$ are orthonormal, (4.6) and (4.7)
imply
$$|H|^2=1+(\alpha')^2={{\lambda_p+\lambda_q}\over n} -{{\lambda_p
\lambda_q}\over {n^2}},\leqno (4.8)$$
which implies that $M$ has constant mean curvature.
And  hence (4.6) and (4.7) also yield
$$||h||^2=||A_\xi ||^2 +n=\lambda_p+\lambda_q.\leqno(4.9)$$
Consequently, by the Gauss equation, we conclude that the scalar curvature
of $M$ is given by
$$\tau= {1\over n}(\lambda_p+\lambda_q)-{1\over{n(n-1)}}\lambda_p
\lambda_q.\leqno(4.10)$$

Conversely, if $M$ has constant mean curvature and constant scalar
curvature, then by
(4.6) we have
$$\Delta^2 x=-n\Delta H=-n||h||^2 H'+n^2|H|^2 x
=||\sigma ||^2\Delta x+(n^2|H|^2-n||\sigma ||^2) x.\leqno(4.11)$$
Since $||h || $ and $|H|$ are constant, this implies that
$M$ has a minimal polynomial of degree $\leq 2$. Hence, by Theorem 2.1,
$M$ is either of 1--type or of 2--type. 

If $M$ is of 1--type and $M$ is not minimal in $S^{n+1}$, then $M$
is contained in the intersection of two hyperspheres of $E^{n+2}$. Thus,
$M$ is a small hypersphere of $S^{n+1}$. If $M$ is of 2--type, then
(4.11) and Hopf's lemma imply that $M$ is mass-symmetric in $S^{n+1}$. \sq

{\bf Remark 4.1.} Theorem 4.1 was first proved in [C5,C8] (see
also [C14]).  Recently, it was proved in [HV2] that Theorem
4.1 is fact of local natural, so that Theorem 4.1 still holds
if the compactness and the mass-symmetric condition in Theorem
4.1 were omitted. 

For spherical 2--type hypersurfaces, we  have the following [C5,C8].

{\bf Theorem 4.3.} {\it Let M be a compact 2--type 
hypersurface of $S^{n+1}$. Then the  geometric invariants:
the mean curvature, the scalar curvature and the length of se\-cond
fundamental form of M in ${\mathbb E}^{n+2}$ are completely determined by 
the order of $M$; more precisely, we have
\begin{itemize}
\item[(i)] the mean curvature $\alpha$ of M in ${\mathbb E}^{n+2}$ is
constant and given by
$$\alpha^2={1\over n}(\lambda_p+\lambda_q)-{1\over{n^2}}\lambda_p
\lambda_q;$$
\item[(ii)] the scalar curvature $\tau$ of M is constant and given by
$$\tau= {1\over n}(\lambda_p+\lambda_q)-{1\over{n(n-1)}}\lambda_p
\lambda_q;$$
\item[(iii)] the length of the second fundamental form $\sigma$
of M in ${\mathbb E}^{n+2}$ is constant and given by
$$||h ||^2=\lambda_p+\lambda_q.$$
\end{itemize}}

Theorem 4.3 shows that the notion of order is fundamental for submanifolds
in Euclidean spaces.

A hypersurface $M$ in \s is called an {\it isoparametric hypersurface\/}
if $M$ has constant principal curvatures. H. F. M\"unzner [M\"u] proved
that the number $g$ of distinct principal curvatures of an isoparametric
hypersurface is 1, 2, 3, 4 or 6. \'E. Cartan has shown that an
isoparametric hypersurface with at most 3  distinct principal
curvature is homogeneous. R. Takagi and T. Takahashi [TT] classified
all homogeneous isoparametric hypersurfaces. For $g=4$ or 6, there
exist non--homogeneous examples. In fact, H. Ozeki and M. Takeuchi 
[OT]  constructed two infinite series of non--homogeneous
isoparametric hypersurfaces. D. Ferus, H.Karcher and H. F. M\"unzner 
[M\"u]
found, for $g=4$, a new type of examples constructed from
representations of a Clifford algebra. J. Dorfmeister and E.
Neher [DN] gave another algebraic approach to isoparametric
hypersurfaces in spheres. 

Since isoparametric hypersurfaces have constant mean curvature and
constant scalar curvature,
Theorem 4.1 implies immediately the following [C5,C8].

{\bf Theorem 4.4.} {\t Every isoparametric 
hypersurface of \s is either of 1--type or of 2--type.} \sq

{\bf Remark 4.2.} Theorem 4.2 shows that there exist abundant examples
of mass--symmetric 2--type hypersurfaces in $S^{n+1}$ \sq.

Although H. B. Lawson proved that there exists many examples of 
minimal surfaces in $S^3$ ({\it i.e.,} 1--type surfaces in $S^3$). 
The following result of [C5,BG,HV1] 
shows that stadard tori
are the only 2--type surfaces in $S^3$.

{\bf Theorem 4.5.} {\it The only 2--type surfaces in $S^3$ are open pieces
of the product of circles: $S^1(a)\times S^1(b)$, $a\not= b$.}
\sq

The following complete classification of compact 2-type hypersurfaces
in the unit 4-sphere is obtained in [C20] (see, also [C22]).

{\bf Theorem 4.6.}  {\it A compact hypersurface
in $S^4(1)\subset {\mathbb E}^5$ is of 2-type if and only if it is one
of the following hypersurfaces:
\begin{itemize}
\item[{(1)}] $S^1(a)\times S^2(b)\subset S^4(1)\subset {\mathbb E}^5$ with
$a^2+b^2=1$ and $(a,b)\not=(\sqrt{1\over 3}
,\sqrt{2\over 3})$ imbedded in the standard way,
\item[{(2)}] a tubular hypersurface with constant radius
$r\not={\pi\over 2}$ about the Veronese minimal imbedding of the real
projective plane $\Bbb RP^2({1\over 3})$ of constant curvature
$1\over 3$ in $S^4(1)$. \sq
\end{itemize}}

For 3-type hypersurfaces we have the following results.

{\bf Theorem 4.7.} {\rm [C13,CL]} {\it Every 3-type
hypersurface in a hypersphere $S^{n+1}\subset {\mathbb E}^{n+2}$ have
non-constant mean curvature. \sq}

{\bf Theorem 4.8.} {\rm [HV4]} {\it A 3-type surface in
the Euclidean 3-space ${\mathbb E}^3$ has non-constant mean curvature.}
\sq 

 For finite type surfaces in $S^3$ we have the following [CD1].

{\bf Theorem 4.9.} {\it Let $M(c)$ be a 
compact surface with constant curvature $c$ in $S^3$. Then $M(c)$
is of finite type if and only if
either (i) $c\geq 1$ and $M(c)$ is totally umbilical in $S^3$ or
(ii) $c=0$ and $M=S^1(a)\times S^1(b)$ with $a^2+b^2=1$.} \sq

Note that the surfaces in (i) are of 1--type and the surfaces in (ii)
are of 2--type unless $a=b$, in which case they are of 1--type.

\vskip.2in
\noindent {\bf \S 5. Classification of submanifolds of finite type.}

\vskip.1in

Let $\gamma: S^1\rightarrow {\mathbb E}^m$ be an isometric immersion from a
circle with radius 1 into a Euclidean $m$--plane ${\mathbb E}^m$. Denote by $s$ 
the arclength of $S^1$. 

For a periodic continuous function $f=f(s)$ with period $2\pi$, $
f(s)$ has a {\it Fourier series expansion\/} given by
$$f(s)={{a_0}\over 2}+a_1\cos s+b_1\sin s +a_2\cos 2s +
b_2\sin 2s +\cdots\,,$$
where $a_k, b_k$ are the Fourier coefficients given by
$$a_k={1\over{\pi}}\int_{-\pi}^{\pi} f(s)\cos ks\,ds,\quad
b_k={1\over{\pi}}\int_{-\pi}^\pi f(s)\sin ks\,ds.$$

In terms of Fourier series expansion, we have the following

{\bf Proposition 5.1.} {\it Let $\gamma: S^1\rightarrow {\mathbb E}^m$ be
a closed smooth curve in ${\mathbb E}^m$. Then $\gamma$ is of finite type
if and only if the Fourier series expansion of each coordinate
function $x_A$ of $\gamma$ has only finite nonzero terms.}

{\bf Proof. } Because $\Delta=-{{d^2}\over{ds^2}}$ in this case,
we have
$$\Delta^j H=(-1)^jx^{(2 j+2)},\quad j=0,1,2,\cdots,$$
where $x^{(j)}={{d^jx}\over {ds^j}}$.
If $\gamma$ is of finite type in ${\mathbb E}^m$, then each Euclidean
coordinate function $x_A$ of $\gamma$ in $E^m$ satisfies the
following homogeneous ordinary differential equation with
constant coefficients:
$$x_A^{(2k+2)}+c_1x_A^{(2k)}+\cdots+c_kx_A^{(2)}=0,\leqno(5.1)$$
for some integer $k\geq 1$ and some constants $c_1,\ldots,c_k$.
Because the solutions of (5.1) are periodic with period $2\pi$,
each solution $x_A$ is a finite linear combination of the following
particular solutions:
$$1,\;\cos n_is,\;\sin m_i s$$
where $n_i,m_i$ are positive integers. Therefore, each $x_A$ is
of the following forms:
$$x_A=c_A+\sum_{p_A}^{q_A}\{a_A(t)\cos ts+b_A(t)\sin ts\}$$
for some suitable constants $a_A(t),b_A(t),c_A$ and some positive
integers $p_A,q_A; A=1,\ldots,m$. Therefore, each coordinate
function $x_A$ has a Fourier series expansion which has only
finite nonzero terms.

Conversely, if each $x_A$ has a Fourier series expansion which 
has only finite nonzero terms, then the position vector $x$ 
of $\gamma$ in ${\mathbb E}^m$ takes the following form:
$$x=c+\sum_{t=p}^q\{a_t\cos ts+b_t\sin ts\} \leqno(5.2)$$
for some constant vectors $a_t,b_t,c$ in $E^m$.
Let $x_t=a_t\cos ts+b_t\sin ts$.
Since $\Delta={{d^2}\over{ds^2}}$, we have $\Delta x_t=t^2 x_t$.
This shows that $\gamma$ is of finite type. \sq

The basic results concerning finite type curves are given in the following.

{\bf Theorem 5.2.}  {\it Let $\gamma$ be a plane curve of finite
type. Then $\gamma$ is a part of a circle or a part of a straight
line. \sq}

In views of Proposition 5.1, Theorem 5.2 can be stated as follows: the
Fourier series expansion of all unit speed curves in ${\mathbb E}^2$, except the
circles and lines, are always infinite. For space curves, the situation
is completely different. In fact, for any $k\in \{2,3,4,\cdots\}$, 
infinitely many distinct $k$--type curves are constructed in [CDDVV].

{\bf Theorem 5.3.}  {\it For every $k\in \{2,3,4,\cdots\}$ there exist
infinitely many non--equivalent curves of k--type in ${\mathbb E}^3$.}
\sq

{\bf Theorem 5.4.}  {\it Every closed space curve of finite type which is
contained in a 2--dimensional sphere is of 1--type, and hence a circle.}
\sq

{\bf Theorem 5.5.}  {\it Every closed space curve of finite type which
is contained in a 3--dimensional sphere is of 1--type, and hence a circle,
or a 2--type W--curve.} \sq

For higher dimensional spheres, a similar result is no longer true. One can
have finite type curves in 4--dimensional spheres that are not 
W--curves. 

{\bf Theorem 5.6.}  {\it Every closed curve of finite type in ${\mathbb E}^3$ having
constant curvature is of 1--type, and hence a circle.} \sq

{\bf Remark 5.1.} Theorem 5.2 was first proved  in [C5] for
closed plane curves. In  [C12], it was pointed out 
that a finite type non--closed curve in a plane is a part of a
straight line without giving detailed proof. 
The detailed proof was then given in  [CDDVV].
Theorems 5.3--5.6 are also given in [CDDVV]. \sq

2--type curves were classified in [CDV] and [DPVV]. It is known that
every 2--type curve in ${\mathbb E}^3$ lies on a hyperboloid of revolution of one 
sheet or on a cone of revolution [DDV]. Curves of finite type in ${\mathbb E}^3$
whose image is contained in a quadratic surface, other than spheres, were
treated in [DDV]. In particular, we have the following.

{\bf Theorem 5.7.} {\rm [DDV]}  {\it  Let $\gamma$ be a closed curve of finite
type in ${\mathbb E}^3$ whose image is contained in a quadratic surface $Q$. Then
\begin{itemize}
\item[(1)] $\gamma$ is a circle, or
\item[(2)] Q is one of the following surfaces:

{\rm (2.a)} an ellipsoid of revolution which is not a sphere;

{\rm (2.b)} a (one-- or two--sheeted) hyperboloid of 
revolution, or

{\rm (2.c)}  a circular cone. \sq
\end{itemize}}

Examples of curves of finite type lying on certain ellipsoids of revolution,
on circular cones or on one--sheeted hyperboloids of revolution are given in
[DDV]. It is not known whether there are curves of finite type on 
two--sheeted hyperboloids of revolution.

Curves of finite type in the context of affine differential geometry were
studied by J. Copaert and L. Verstraelen (cf. [V1]). Quite remarkable
is that all affine curve of finite type are affinely equivalent to
skew lines. In Euclidean spaces, a similar statement is false.

Now, we present some classification theorems for
surfaces of finite type in Euclidean spaces.

{\bf Theorem 5.8.}  {\it The only tubes of finite
type in ${\mathbb E}^3$ are circular cylinders.} \sq

{\bf Theorem 5.9.}   {\it A ruled surface M in ${\mathbb E}^m$ is of finite
type if and only if M is a cylinder on a curve of finite type or
M is a part of a helicoid in an affine subspace ${\mathbb E}^3$.} \sq

{\bf Theorem 5.10.}  {\it A ruled surface $M$ in ${\mathbb E}^3$ is of finite
type if and only if M is a part of a plane, a circular cylinder
or a helicoid.} \sq

{\bf Corollary 5.11.}   {\it A flat surface in ${\mathbb E}^3$ is of finite type
if and only if it is a part of plane or a circular cylinder.}
\sq

{\bf Remark 5.2.} Theorem 5.8 was proved in [C10] and 
Theorem 5.9, Theorem 5.10 and Corollary 5.11 are due to 
Chen--Dillen--Vers\-trae\-len--Vrancken [CDVV1]. \sq

For algebraic surfaces of degree 2 we have the following result of
[CD2].

{\bf Theorem 5.12.}  {\it The only quadrics in ${\mathbb E}^3$ which is of
finite type  are the spheres and the circular cylinders.} \sq

The complete classification of finite type algebraic
hypersurfaces of degree 2 is obtained in [CDS].

In 1822, C. Dupin defined a {\it cyclide\/} to be a surface $M$
in ${\mathbb E}^3$ which is the envelope of the family of spheres tangent to
three fixed spheres. This was shown to be equivalent to requiring
that both sheets of the focal set degenerate into curves. The 
cyclides are equivalently characterized by requiring that the
lines of curvatures in both families be arcs of circles or straight
lines. Thus, one can obtain three obvious examples: a torus of
revolution, a circular cylinder and a circular cone. It turns
out that all cyclides can be obtained from these three by inversions
in a sphere in $E^3$. 

Recently, F. Defever, R. Deszcz and L. Verstraelen
[DDV1,DDV2] proved the following.

{\bf Theorem 5.13.}  {\it Cylides of Dupin are of infinite type.}
\sq

All of the above results support the following conjecture.

{\bf Conjecture 5.1.}  {\it The only compact surfaces of finite type
in ${\mathbb E}^3$ are spheres.} \sq

From Theorem 5.2 we know that, in one dimension case, this is true:
 namely, the
{\it circles are the only closed planar curves of finite type.}

{\bf Theorem 5.2'.} If $\gamma: S^1\rightarrow {\mathbb E}^2$ is
an isometric immersion from $S^1$ into ${\mathbb E}^2$. Then $\gamma$ is
of finite type if and only if $\gamma$ is a circle. \sq

\vskip.2in
\noindent {\bf \S 6. Linearly independent and orthogonal immersions. }
\vskip.1in

   Let $x: M\rightarrow {\mathbb E}^m$ be an immersion of $k$-type. Suppose $$x =
c+x_{1}+\ldots +x_{k},\hskip.3in \Delta
x_{i}=\lambda_{i}x_{i},\hskip.2in \lambda_{1}<\ldots
<\lambda_{k},\leqno(6.1)$$ is the spectral decomposition of
the immersion $x$, where $c$ is  a constant vector and
$x_{1},\ldots,x_{k}$ are non-constant maps. For each $i\in
\{1,\ldots,k\}$ we put
 $$E_{i}=Span\{x_{i}(p):p\in M\}.$$
Then each $E_i$ is a linear subspace of ${\mathbb E}^m$.  

 Here, we recall  the  notions of linearly independent
 and orthogonal immersions first introduced in [C15]. 

{\bf Definition 6.1.}
 Let $x:M \rightarrow
{\mathbb E}^m$ be an immersion of $k$-type whose spectral
decomposition is given by (6.1). Then the immersion $x$ is
said to be {\it linearly independent\/} if the subspaces
$E_{1},\ldots,E_{k}$ are linearly independent, that is,
the dimension of the subspace spanned by all vectors in $E_{1}\cup
\ldots \cup E_{k}$ is equal to dim $E_{1}+\ldots+$ dim
$E_k$. And the immersion $x$ is said to be {\it
orthogonal\/} if the subspaces $E_{1},\ldots,E_{k}$ are
mutually orthogonal in ${\mathbb E}^m$. \sq

 Obiviously, every 1-type immersion is an orthogonal
immersion and hence a linearly independent immersion. 
There exist ample examples of orthogonal immersions and
ample examples of linearly independent immersions which
are not orthogonal. For instances, every $k$-type curve
lying fully in ${\mathbb E}^{2k}$ is an example of linearly
independent immersion and the immersion of the curve is
orthogonal if and only if it is a $W$-curve. Moreover, the product of any two
linearly independent immersions is again a linearly
independent immersion.

 Let  $x: M \rightarrow {\mathbb E}^{m}$ be an immersion of $k$-type
whose spectral decomposition is given by (6.1). 
We  choose a Euclidean coordinate system
$(u_{1},\ldots, u_{m})$ on ${\mathbb E}^m$ with $c$ as its 
origin. Then we have the
following spectral decomposition of $x$:
 $$x = x_{1}+\ldots
+x_{k},\hskip.3in \Delta x_{i}=\lambda_{i}x_{i},\quad
\lambda_{1}<\ldots <\lambda_{k}.\leqno(6.2)$$ For each
$i\in \{1,\ldots,k\}\,$ we choose a basis $\{c_{ij}:
j=1,\ldots,m_{i}\}$ of $E_i$, where $m_{i}$ is the
dimension of $E_i$. Put $\ell=m_{1}+\ldots+m_{k}$ and let
$E^{\ell}$ denote the  subspace of ${\mathbb E}^m$ spanned by
$E_{1},\ldots, E_k$. If the immersion $x$ is linearly
independent, then the vectors $\{c_{ij}:\, i=1,\ldots,k;
j=1,\ldots,m_{i}\}$ are $\ell$ linearly independent vectors
in ${\mathbb E}^{\ell}$. Furthermore, we choose the Euclidean
coordinate system $(u_{1},\ldots,u_{m})$ on ${\mathbb E}^m$ such that
${\mathbb E}^{\ell}$ is  defined by $u_{\ell +1}=\ldots =u_{m}=0.$
Regard each $c_{ij}$ as a column $\ell$-vector. We put 
$$S =
(c_{11},\ldots,c_{1m_{1}},\ldots,c_{k1},\ldots,c_{km_{k}}).
\leqno(6.3)$$ Then the matrix $S$ is a nonsingular
$\ell\times \ell$ matrix. Let $D$ denote the diagonal
$\ell\times \ell$ matrix given by
$$D=Diag(\overbrace{\lambda_{1},\ldots,\lambda_{1}}^{m_1\; \hbox{times}}
,\ldots,\overbrace{\lambda_{k},\ldots,\lambda_{k}}^{m_k\; 
\hbox{times}}),\leqno(6.4)$$ 
 where $\lambda_{i}$
repeats $m_i$-times. If we put $A=SDS^{-1}$, then
$Ac_{ij}=\lambda_{i}c_{ij}$ for any $i\in \{1,\ldots,k\}$
and $j\in \{1,\ldots,m_{i}\}$. Therefore, we have $$\Delta
x=Ax\leqno(6.5)'$$  for the immersion $x:M\rightarrow 
{\mathbb E}^{\ell}$ induced from the original immersion
$x:M\rightarrow {\mathbb E}^m$.  By regarding the $\ell\times \ell$
matrix $A$ as an $m\times m$ matrix in a natural way (with
zeros for each of the additional entries), we obtain
$$\Delta x=Ax,\hskip.2in A=(a_{ij})\leqno(6.5)$$ for the
immersion $x:M\rightarrow {\mathbb E}^m$. 

  If $x: M\rightarrow {\mathbb E}^m$ is a minimal immersion, then
we have $A=0$. If $x:M\rightarrow {\mathbb E}^{\ell}\; (\subset
{\mathbb E}^{m})$ is a non-minimal full immersion, then the Euclidean
coordinate functions $u_{1},\ldots,u_{\ell}$ of
${\mathbb E}^{\ell}$, restricted to $M$, do not satisfy any linear
equation. Thus, the coordinate functions 
$u_{1}|_{M},\ldots,u_{\ell}|_{M}$ of $M$ in ${\mathbb E}^{\ell}$ are
linearly independent functions. Therefore, if $B$ is any
$\ell\times {\ell}$ matrix such that $\Delta x=Bx$,
then $A=B$. Hence, the ${\ell}\times \ell$ matrix $A$ 
in $(6.5)'$ defined above is unique. Consequently, if the
(original) immersion $x:M\rightarrow {\mathbb E}^{m}$ is a
non-minimal, linearly independent immersion, then the
$m\times m$ matrix $A$ given in (6.5) is also uniquely
defined (with respect to the Euclidean coordinate system
so chosen).

 By using this   $m\times m$ matrix $A$ in (6.5)
defined above, the first author introduced in [C15]
 the notion of {\it adjoint hyperquadric\/}
as follows .

{\bf Definition 6.2} Let $x:M \rightarrow
{\mathbb E}^m$ be a non-minimal, linearly independent immersion 
whose spectral decomposition is given by (6.1).
 Let $u=(u_{1},\ldots,u_{m})$
 be a Euclidean coordinate system  
on ${\mathbb E}^m$ with $c$ as its origin
 and let $A$ be the $m\times m$ matrix in (6.5) associated
with the immersion $x$ defined above. Then, for any point
$p\in M$, the equation
 $$\<Au,u\>:=\sum_{i,j}^{m} a_{ij}u_{i}u_{j}=c_{p},\,\,\,\,\,
(c_{p}=\< Ax,x\> (p))\leqno(6.6)$$ defines a
hyperquadric $Q_p$
 in ${\mathbb E}^m$. We call the hyperquadric
$Q_p$ the {\it adjoint hyperquadric of the immersion
$x$ at $p$.}
In particular, if $x(M)$ is contained in
an adjoint hyperquadric $Q_p$ of $x$ for some point $p \in
M$, then all of the adjoint hyperquadrics $\{ Q_{p}: p\in
M\}$ give a common adjoint hyperquadric, denoted by $Q$. We
call the hyperquadric $Q$  the {\it adjoint hyperquadric
of the linearly independent immersion $x$\/.} \sq

  The following result from [C15] provides us a necessary and
sufficient condition for a compact, linearly independent 
 submanifold to lie in its adjoint hyperquadric.

{\bf Theorem 6.1.}  {\it  Let $x:M \rightarrow {\mathbb E}^m$ be a
 linearly independent immersion  from a compact
manifold $M$ into ${\mathbb E}^m$ whose spectral decomposition is
given by (6.2). Then $M$ is immersed into the adjoint
hyperquadric of $x$ if and only if  $x(M)$ is
contained in a hypersphere of ${\mathbb E}^m$ centered at the
origin.}

{\bf Proof. } Assume that
$M$ is immersed into a hypersphere $S^{m-1}(r)$ of ${\mathbb E}^m$
with radius $r$ centered at the origin. Denote by $H$ and
${\bar H}$ the mean curvature vectors of $M$ in ${\mathbb E}^m$ and
of $M$ in $S^{m-1}(r)$, respectively. Then we have
$$H={\bar H} - {1\over r}\, x.\leqno(6.7)$$
This implies $\<H,x\>=-r$. Since $\Delta x = -nH, n=dim\, M$,
we obtain $\<\Delta x,x\>=nr$. Therefore, by (6.5), we
conclude that $M$ is immersed into the adjoint
hyperquadric defined by $\<Au,u\>=nr$, where $nr$ is a
constant.

 Conversely, suppose that $x:M\rightarrow {\mathbb E}^m$ is a
linearly independent immersion of a compact manifold $M$  
such that $x(M)$ is contained in an adjoint hyperquadric
$Q_p$ for some point $p$. Then we have $\<Ax,x\>=c_p$ where
$c_p$ is the constant given by $c_{p}=\<Ax,x\>(p)$. Since
$Ax=\Delta x = -nH$, we have
$$\<nH,x\>=-c_{p}.\leqno(6.8)$$
Because $M$ is compact, we also have 
$$\int_{M} \{1+\<H,x\>\}*1=0.\leqno(6.9)$$
Formulas (6.8) and (6.9) imply $c_{p}=-n$. Therefore,
$$\Delta \<x,x\>=2\<\Delta x,x\>-2n = -2n(\<H,x\>+1)=0.$$
Thus, by Hopf's lemma, $\<x,x\>$ is a constant and $M$
is immersed into a hypersphere of ${\mathbb E}^m$ centered at the
origin. \sq

{\bf Remark 6.1.} Although the implication
$(\Longleftarrow )$ in Theorem 6.1 holds in general
without the asssumption of compactness, the implication
$(\Longrightarrow )$ does not hold in general if $M$ is
not compact. For example, let $M$ be the product surface
of 
 a unit plane circle and a line. Then the inclusion map
$x$ of $M$ in $E^3$  defined by 
$$x(u,v)=(\cos u, \sin u, v)$$
is a non-spherical, linearly independent immersion 
whose adjoint hyperquadric  is given by
$u_{1}^{2}+u_{2}^{2}=1$. It is clear that $M$ coincides
with the  adjoint hyperquadric $Q$ of $M$ in ${\mathbb E}^3$. \sq

 The following result from [C15] provides us a necessary and
sufficient condition for a linearly independent immersion
to be an orthogonal immersion in terms of the adjoint
hyperquadric.

{\bf Theorem 6.2.}   {\it Let $x:M \rightarrow
{\mathbb E}^m$ be a non-minimal, linearly independent immersion.
 Then $M$
is immersed by $x$ as a minimal submanifold of the adjoint
hyperquadric if and only if the immersion $x$ is an
orthogonal immersion.} \sq

In general a submanifold obtained from an equivariant
immersion of a compact homogeneous space is not a
minimal submanifold of any hypersphere. However, the
following result from [C15] shows that such a submanifold is
always a minimal submanifold in its adjoint hyperquadric.

{\bf Theorem 6.3.}  {\it Let $x: M\rightarrow
{\mathbb E}^{m}$ be an equivariant isometric immersion of a compact
$n$-dimensional Riemannian homogeneous space $M$ into
${\mathbb E}^m$. Then $M$ is
immersed as a minimal submanifold of the adjoint
hyperquadric.} \sq

The class of 1-type immersions has been classified by T.
Takahashi. In fact, as stated before, he showed that the
submanifolds of ${\mathbb E}^m$ for which
$$\Delta x=\lambda x\leqno(6.10)$$
are precisely either the minimal submanifolds of ${\mathbb E}^m\,
(\lambda=0)$ or the minimal submanifolds of the hyperspheres
$S^{m-1}$ in ${\mathbb E}^m$ (the case when $\lambda\not= 0,$
actually 
$\lambda >0$).
As a geralization of Takahashi's result, Garay 
studied the hypersurfaces $M^n$ in ${\mathbb E}^{n+1}$ for which
$$\Delta x=Ax,\leqno(6.11)$$
where $A$ is a diagonal matrix.

  Dillen, Pas and Verstraelen observed that
Garay's condition is not coordi\-nate--invariant and they
considered the submanifolds in ${\mathbb E}^m$ for which 
$$\Delta x = Ax+b \leqno(6.12)$$
where $A\in  {\mathbb E}^{m\times m}$ and $b\in {\mathbb E}^m$.
This setting generalizes Takahashi's condition in a way
which is independent of the choice of coordinates. In
1988 Garay proved that if a hypersurface $M$ in
${\mathbb E}^{n+1}$ satisfies his condition, it is either minimal in
${\mathbb E}^{n+1}$, or it is a hypersphere, or it is a spherical
hypercylinder.  In 1990 Dillen, Pas and Verstraelen
proved that a surface in ${\mathbb E}^3$ satisfies their condition
if and only if it is an open part of a minimal surface, a
sphere, or a circular cylinder. 

In [CP], Chen and Petrovic  proved the following. 

{\bf Theorem 6.4.} {\rm (Characterization)}   {\it  Let
$x:M\rightarrow {\mathbb E}^m$ be an immersion of finite type. Then
the immersion  $x$ is linearly independent if and only if
$x$ satisfies $\Delta x=Ax+b$ for some $A\in  {\mathbb E}^{m\times
m}$  and $b\in {\mathbb E}^m$.} \sq

{\bf Theorem 6.5.} {\rm (Characterization)}   {\it 
 Let $x:M\rightarrow {\mathbb E}^m$
be an immersion of finite type. Then the immersion  $x$
is orthogonal if and only if $x$ satisfies
$\Delta x=Ax+b$ for some symmetric matrix $A\in 
	{\mathbb E}^{m\times m}$ and $b\in  {\mathbb E}^m$.} \sq

{\bf Theorem 6.6.} {\rm (Classification)}  {\it Let $x: M\rightarrow
{\mathbb E}^{n+1}$ be an isometric immersion from an n--dimensional manifold
 $M$ in ${\mathbb E}^{n+1}$. Then the immersion x is linearly independent if
and only if M is an open portion of minimal
hypersurface, a hypersphere $S^n$, or a spherical
hypercylinder $S^\ell \times {\mathbb E}^{n-\ell},\,\ell \in
\{1,2,\ldots,n-1\}$.} \sq

Theorem 6.6 was obtained in [CP], also independently in [HV3].

The class of linearly independent immersions is contained in a much
larger class of immersions, namely the class of  {\it
immersions of restricted type.} A submanifold of a Euclidean (or
  pseudo-Euclidean) space is said to be of restricted type
if its shape operator with respect to the mean curvature
vector is the restriction of a fixed linear transformation of
the ambient space to the tangent space of the submanifold at
every point of the submanifold. The notion of immersions of
restricted type  is first introduced in [CDVV3]. Further
results can be found, e.g. in [DVVW] and [BBCD].

\vskip.2in
\noindent {\bf \S 7. Variational minimal principle.}
\vskip.1in

In this section we mention the fundamental relationship between the
theory of finite type and calculus of variations obtained in
[CDVV2,CDVV4].

Let \x  be an isometric immersion of a compact Riemannian
manifold $M$ in a Euclidean space ${\mathbb E}^{m}$. Associated to each
${\mathbb E}^{m}$-valued
vector field $\xi$ defined on $M$, there is a deformation $\phi_t$, defined
by
$$
\phi_t(p):=x(p)+t \xi(p),\quad p\in M,\quad t\in (-\epsilon,
\epsilon), $$ where $\epsilon$ is a sufficiently small positive
number. For each $t$, $\phi_t$ gives rise to a submanifold
$M_t=\phi_t(M)$. Let $A(t)$ denote the area of $M_t$.

 Let $\Cal D$ denote the class of all deformations acting on
the submanifold $M$ and let  ${\mathbb E}E$ denote a nonempty subclass of
$\Cal D$. A compact submanifold $M$ in ${\mathbb E}^{m}$ is said to
{\it satisfy the variational minimal principle\/} in the class ${\mathbb E}E$
if $M$ is a critical point of  the volume functional for all
deformations in ${\mathbb E}E$, {\it i.e.,} for each deformation in ${\mathbb E}E$,
one has $A'(0)=0$. A compact submanifold $M$ in ${\mathbb E}^{m}$ is
called stable in the class ${\mathbb E}E$ if it satisfies the variational
minimal principal in the class ${\mathbb E}E$ and if  $A''(0)\geq 0$
 for each deformation in ${\mathbb E}E$.

We consider one type of deformations which occurs frequently, namely the
directional deformations.  Such deformations occur when an object moves in a
fixed direction.

Directional deformations are defined as follows : let $c$ be a fixed
vector in ${\mathbb E}^{m}$ and let $f$ be a smooth function defined on the
submanifold $M$.  Then we have a deformation given by
$$
\phi_t^{fc}(p):=x(p)+t f(p) c,\quad p\in M\quad t\in (-\epsilon,
\epsilon).\leqno(7.1)
$$
Such a deformation is called a {\it directional deformation\/} in the
direction $c$.

For each $q\in \Bbb N$, we define
$\CC_q$ to be the class of all
directional deformations given by smooth functions
$f\in \sum_{i\geq q}V_i$. It is clear that
$\CC_0 \supset \CC_1 \supset \CC_2 \supset\dots\supset\CC_k\supset\cdots$.
Therefore, if a compact submanifold $M$ of ${\mathbb E}^{m}$ satisfies the
variational minimal principle for one class $\CC_k$,
then it automatically satisfies the variational minimal principle
in $\CC_{\ell}$  for $\ell\geq k$.

{\bf Theorem 7.1.}  {\rm [CDVV2,CDVV4]}  {\it We have:
\begin{itemize}
\item[(a)] There are no compact submanifolds in
${\mathbb E}^{m}$ which satisfy the variational minimal principle in the
classes $\CC_0$ and $\CC_1$.

\item[(b)] A compact submanifold  $M$ of ${\mathbb E}^{m}$ is of finite type
if and only if it satisfies the variational minimal principle in
the class $\CC_q$ for some $q\geq2$.
\end{itemize}}

{\bf Proof. } Let $c\in {\mathbb E}^{m}$ and $f\in C^{\infty}(M)$. Consider
 the directional deformation
$\phi_t^{fc}$ defined by \thetag{7.1}. Let $e_1,\ldots,e_n$ be
an orthonormal local frame field on $M$. Extended
$e_1,\ldots,e_n$ by $(\phi_t)_*e_1,\ldots,(\phi_t)_*e_n$. Put
$$g_{ij}(t)=\<(\phi_t)_*e_i,(\phi_t)_*e_j\>.\leqno(7.2)$$
Then $$A(t)=\int_M \sqrt{det(g_{ij}(t))}dA.\leqno(7.3)$$

On the other hand, from (7.1), we have
$$(\phi_t)_*e_i=e_i+t(f_i)c,\quad {\partial\over{\partial
t}}=fc,\leqno (7.4)$$
where $f_i=e_if$. From (7.2) and (7.3) we get
$$g_{ij}(t)=\delta_{ij}+t(f_i\<c,e_j\>+f_j\<c,e_i\>)+t^2f_if_j\<
c,c\>.\leqno(7.5)$$
By  using (7.3) and (7.5), we may obtain
$$
A(t)=\int_M\{ 1+2t\<c,\nabla f\>+ t^2(\< c^{\perp},c^{\perp}\>|\nabla f|^2
+
\<c,\nabla f\>^2) \}^{\tfrac12} dA,\leqno(7.6)
$$
where $\< \phantom{v},\phantom{v}\>$ is the inner product of
${\mathbb E}^{m}$, $c^{\perp}$
the normal component of $c$, $\nabla f$ the gradient of $f$ and $dA$ the
volume
element of $M$.

From \thetag{7.6} we obtain
$$
A'(0)=\int_M\<c,\nabla f\> dA.\leqno(7.7)$$
Let $c^T$ denote the tangential component of $c$ and $(c^T)^{\sharp}$
the one-form on $M$ dual to $c^T$. Then we have
$$
(c^T)^{\sharp}=dh_c,\leqno(7.8)
$$
where $h_c$ denotes the height function of $M$ in ${\mathbb E}^{m}$ with
respect to $c$.

Denote by $d$ and $\delta$ the differential and the co-differential
operators of
$M$. By using \thetag{7.8} we find
$$
\delta(c^T)^{\sharp}=\lap h_c=\<\lap x ,c\>,\leqno(7.9)
$$
where $x$ is the position vector field of $M$ in ${\mathbb E}^{m}$. On the
other hand, the well-known Beltrami formula yields
$$
\lap x = -n H, \quad\quad n=\dim M,\leqno(7.10)
$$
where $H$ is the mean curvature vector of $M$ in ${\mathbb E}^{m}$. Hence, by
using \thetag{7.7}, \thetag{7.9} and \thetag{7.10}, we have
$$
A'(0)=\int_M H_c f dA,\leqno(7.11)
$$
where $H_c=\<H,c\>$.

If $M$ is a compact submanifold in ${\mathbb E}^{m}$ which satisfies the
variational minimal principle in the class $\CC_1$, then
\thetag{7.11} yields $$
(H_c,f)=\int_M H_c f dA=0.\leqno(7.12)
$$
Since \thetag{7.12} holds for any $f\in \sum_{t=1}^\infty V_t$ and
any $c\in {\mathbb E}^{m}$, the mean curvature vector field $H$ is a
constant vector. Therefore, 
can conclude that $H=0$, which is a 
compact. the Weingarten operator $A_H$ with respect to $H$
vanishes. Taking the trace of $A_H=0$, we obtain $H=0$, which
means that the submanifold would be minimal. This is impossible
since $M$ is compact. Consequently, $M$ cannot satisfy the
variational minimal principle in the class $\CC_1$. Hence, it
does not satisfy the variational minimal principle in the class
$\CC_0$. This proves (a).

Now we prove (b).
First, assume that $M$ is a compact submanifold of finite type. Then the
position vector field of $M$ in ${\mathbb E}^{m}$ has a finite spectral
decomposition:
$$
x= x_0 + x_{i_1} + \dots + x_{i_k},\leqno(7.13)
$$
where $x_0$ is a constant vector and $\lap x_{i_j}=\l_{i_j}x_{i_j}$ for
$j=1,\dots,k$. Assume $\l_{i_1}<\dots<\l_{i_k}$. From
\thetag{7.10} and \thetag{7.13} we obtain
$$
-nH= \l_{i_1}x_{i_1} + \dots + \l_{i_k} x_{i_k}.\leqno(7.14)
$$
Equation \thetag{7.14} implies that $H_c\in \sum_{j\le\ell} V_j$,
where $\ell$ is the upper order of  $M$. Put $q=\ell +1$. Then
\thetag{7.11} implies that $A'(0)=0$ for any deformation in
$\CC_q$.  Hence $M$ satisfies the variational minimal
principle in $\CC_q$. Obviously $q\geq 2$.

Conversely, if $M$ satisfies the variational minimal principal in $\CC_q$,
for some $q\geq 2$, then \thetag{7.12} implies that $H_c\in
\sum_{i<q} V_i$. Since this holds for any $c$,  $H$ has finite
spectral decomposition $H=H_0+H_1+\dots+H_{q-1}$, $\lap H_j=\l_j
H_j$. We put $P(u)= u \prod_{j=1}^{q-1}(u-\l_j)$. Then
$P(\lap)H=0$. Hence,  by applying the minimal polynomial criterion 
(Theorem 3.1), $M$ is of finite type. \qed

{\bf Theorem 7.2.}  {\it Every compact submanifold $M$ of finite
type in  a Euclidean space is stable in the class $\CC_q$ for any
$q\geq \hbox{u.o.}(M)+1$.}

{\bf Proof. }. Let $M$ be a compact submanifold of finite type in
${\mathbb E}^{m}$, then, from the proof of Theorem 7.1, we see that $M$
satisfies the variational minimal principal in the class $\CC_q$
for any $q\geq \hbox{u.o.}(M)+1$. On the other hand, from (7.6),
we have $$A''(0)=\int_M \< c^{\perp},c^{\perp}\>|\nabla f|^2
dA\geq 0.\leqno(7.15) $$
Therefore, $M$ is  stable in the class $\CC_q$
for any $q\geq \hbox{u.o.}(M)+1$.\qed

Now, we give some consequences of Theorem 7.1.

{\bf Corollary 7.1.}   {\it A compact
hypersurface  in a Euclidean space is a hypersphere if and only if
it satisfies the variational minimal principle in the class $\CC_2$.}
\sq

{\bf Corollary 7.2.}  {\it Circles in ${\mathbb E}^2$ are the only closed
planar curves which satisfy the variational minimal principle in
the class $\CC_q$ for some $q\geq 2$.} \sq

{\bf Corollary 7.3.}  {\it Every $k$-th standard immersion of a
compact homogeneous Riemannian manifold in ${\mathbb E}^{m}$ satisfies the
variational minimal principle in the class $\CC_{k+1}$.} \sq

For hypersurfaces in ${\mathbb E}^{n+1}$ we may prove the following.

{\bf Theorem 7.3.}  {\rm [CDVV4]}  {\it A compact embedded hypersurface $M$ in
${\mathbb E}^{n+1}$ has finite type Gauss map if and only if the volume
enclosed by $M$ is invariant
under the directional deformations in a class $\CC_q$, for some $q\geq 2$.}

{\bf Proof. } If $M$ is an embedded hypersurface in ${\mathbb E}^{n+1}$, $M$
can be regarded as the boundary $\partial\Omega$ of a domain
$\Omega\subset {\mathbb E}^{n+1}$.
 For  any fixed unit vector $c\in {\mathbb E}^{n+1}$, we
choose a Euclidean coordinate system
$(u_1,u_2,\ldots,u_{n+1})$ such that $c=(1,0,\ldots,0)$.
 Let $X=u_1c$. Then $X$ is a vector field defined on
${\mathbb E}^{n+1}$. According to the divergence theorem, we know the
volume $V$ enclosed by $M$ is given by
$$V=\int_\Omega
(\hbox{div}\;X)dV=\int_{M=\partial\Omega} \<\xi,X\> dA,\leqno(7.16)$$
where $\xi$ is the unit outward normal vector field of $M$ in
${\mathbb E}^{n+1}$ and $dA$ denotes the area element of $M$. Let
$e_1,\ldots,e_n$ be an positive-oriented orthonormal local
frame field of the tangent bundle of $M$. Then, by (7.16),  we
have  $$V=\int_M u_1\<e_1\times\cdots\times e_n,c\>dA,\leqno(7.17)$$
where $e_1\times\cdots\times e_n$ denotes the vector product of
$e_1,\ldots,e_n$ in ${\mathbb E}^{n+1}$.

For a given function $f$ on $M$, consider the directional
deformation given by $\phi_t(p)=x(p) +tf(p)c$, $t\in
(-\epsilon,-\epsilon)$. Let $(u_1)_t$ denote the $u_1$-component
of $M_t=\phi_t(M)$. Then  $$(u_1)_t=u_1+tf,\quad
(\phi_t)_*e_i=e_i+t(e_if)c. \leqno(7.18)$$ Let $V(t)$ denote the
volume enclosed by $M_t$. Then, from (7.17) and (7.18), we find
$$V(t)=V(0)+t\int_M<\xi,c>fdA.\leqno(7.19)$$ Formula (7.19) implies
$$V'(t)=\int_M<\xi,c>fdA.\leqno(7.20)$$
We recall that  the Gauss
map $\nu$ of $M$ in ${\mathbb E}^{n+1}$ is given by $\nu(p)=\xi(p)$ and
the Gauss map $\nu$ is mass-symmetric in the unit hypersphere
of ${\mathbb E}^{n+1}$ centered at the origin.

If the Gauss map $\nu$ of $M$  is of finite type, then $\nu$
has a finite spectral decomposition:
$$\nu=\nu_1+\cdots+\nu_k\leqno(7.21)$$
where $\nu_1,\ldots,\nu_k$ are ${\mathbb E}^{n+1}$-valued eigenfunctions
of $\Delta$. Let $\ell$ is the upper order of the Gauss map
$\mu$. Put $q=\ell+1$. Then (7.20) implies that $V'(0)=0$ for any
deformation in the class $\CC_q$. Thus, according to (7.19),
$V(t)=V(0)$ for any $t\in(-\epsilon,-\epsilon)$. This means that
the volume enclosed by $M$ is invariant under directional
deformations in the class $\CC_q$.

Conversely, if  the volume
enclosed by $M$ is invariant under directional deformations in
the class  $\CC_q$, for some $q\geq 2$, then $V(t)=V(0)$ for
$t\in(-\epsilon,-\epsilon)$ and $f\in \sum_{i\geqq} V_i$. Thus,
by applying (7.19), we get $\<\xi,c\>\in \sum_{i< q}V_i$. Since
this is true for any $c\in {\mathbb E}{n+1}$, we get
$\nu=\nu_1+\cdots +\nu_{q-1}$, $\Delta \nu_i=\lambda_i\nu_i$.
We put $P(t)=\prod_{j=1}^{q-1}(t-\lambda_j)$. Then
$P(\Delta)\nu=0$. Therefore,  by Theorem 2.2 of [CPi], the Gauss
map of $M$ is of finite type. \qed 

{\bf Theorem 7.4.} {\rm [CDVV4]}  {\it A compact embedded
hypersurface $M$ in a Euclidean space ${\mathbb E}^{n+1}$ is a
hypersphere if and only if the  volume enclosed by $M$ is invariant under
directional  deformations in  $\CC_2$.}

{\bf Proof. } If $M$ is a hypersphere of ${\mathbb E}^{n+1}$, then the
Gauss map of $M$ satisfies $\Delta\nu=\lambda_1\nu$, where
$\lambda_1$ is the first nonzero eigenvalue of $\Delta$. Thus,
from the proof of Theorem 7.3, we see that the volume enclosed
by $M$ is invariant under directional deformations in $\CC_2$.

Conversely, if the volume enclosed
by $M$ is invariant under directional deformations in $\CC_2$,
then $\nu=\nu_1,\;\Delta\nu_1=\lambda_1\nu_1$, since $\nu$ is
mass-symmetric. In particular, this implies that the Gauss map
is of 1-type. Therefore, by applying Theorem 4.2 of [CPi], $M$ is
a hypersphere in ${\mathbb E}^{n+1}$.
\qed

For maps of finite type we have the following results
[CDVV4].

 {\bf Theorem 7.5.}  {\it 
 A smooth map $\phi:M\to {\mathbb E}^{m}$ of a compact Riemannian
manifold  $M$ in ${\mathbb E}^{m}$ is of finite type if and only if it
satisfies the variational minimal
principle in the class $\CC_q$ for some $q\in \Bbb N$.} \sq

{\bf Theorem 7.6,}   {\it A smooth map $\phi:M\to {\mathbb E}^{m}$ of a
compact Riemannian manifold  $M$ in ${\mathbb E}^{m}$ is a constant map if
and only if it satisfies the variational minimal principle in the
class $\CC_1$.} \sq

\vfill\eject

\noindent {\bf Chapter III: LAPLACE MAPS OF SMALL RANK}

\vskip.2in
\noindent {\bf \S1. Curves in ${\mathbb E}^m$.}
\vskip.1in

In this section  we  discuss the Laplace map of a regular curve in a
Euclidean $m$--space. 

Consider a regular curve $\beta:I\rightarrow {\mathbb E}^m$
from an open interval $I$ into ${\mathbb E}^m$ parame\-tri\-zed by the arclength $s$. Then
the Laplacian operator $\Delta$ of $\beta$ is given by
$\Delta= -{{d^2}\over {ds^2}}$. Therefore, the Laplace map of $\beta$ is given by
$L(s)=\Delta \beta=-\beta''(s)$. Let $t=\beta'={{d\beta}\over{ds}}$ be the unit 
tangent vector of $\beta$. The differential of the Laplace map $L$ is given by
$L_*(t)=-\beta'''$. Recall that the curve $\beta$
 is said to be regular if $\beta'(s)$ is nowhere
zero on $I$. In this monograph, by a {\it 3--regular curve\/} we mean a 
regular curve $\beta:I\rightarrow {\mathbb E}^m$ such that  $\beta'''(s)$ is nowhere
zero.

If $\beta(s)$ is a regular plane curve, then there is a unique unit normal
vector field $n(s)$ such that $\{t(s),n(s)\}$ gives a right handed orthonormal
basis of $E^2$ for each $s$. The {\it plane curvature\/} $\kappa(s)$ of $\beta$
is defined by $\kappa(s)=<\beta''(s),n(s)>$, where $<\;,\; >$ is the 
Euclidean inner product on ${\mathbb E}^2$. It is easy to see that a regular plane
curve is 3--regular if the plane curvature and its derivative do not vanish
simultaneously. Clearly, there exist regular 3-regular plane curves whose
derivative of its planar curvature vanishes at a point; e.g. the plane
polynomial spirals (cf. [D1]) with curvature function $\kappa(s)=as^2 +bs
+c$ such that the discriminant $b^2-4ac\not=0$.

For curves in higher dimensional Euclidean spaces,
  we explain the Frenet curvatures and Frenet vectors 
of a regular curve $\beta:I\rightarrow E^m$ as folows.

Let $\beta_1=\beta'$ be the unit tangent vector and put $\kappa_1
=||{\tilde\nabla}_{\beta'}\beta_1 ||.$ If $\kappa_1$ vanishes on $I$, then
the curve $\beta$ is said to be of rank one. If $\beta_1$ is not identically
zero, then we define $\beta_2$ by 
$${\tilde\nabla}_{\beta'}\beta_1=\kappa_1\beta_2,\quad \hbox{on}\quad
I_1=\{s\in I:\kappa_1(s)\not=0\},\leqno(1.1)$$
$\beta_2$ is called the principaal normal vector of $\beta$.

Put
$$\kappa_2=|| {\tilde\nabla}_{\beta'}\beta_2 +\kappa_1\beta_1 ||.\leqno(1.2)$$
If $\kappa_2\equiv 0$ on $I_1$, then $\beta$ is said to be of rank 2. If
 $\kappa_2$ is not identically zero on $I_1$, then we define $\beta_3$ by
$${\tilde\nabla}_{\beta'}\beta_2=-\kappa_1\beta_1+\kappa_2\beta_3,\quad
\leqno (1.3)$$ on $I_2=\{s\in I:\kappa_2(s)\not=0\}$. For $m=3$, $\beta_3$ is
called the binormal vector of $\beta$. Inductively, we put 
$$\kappa_i=|| {\tilde\nabla}_{\beta'}\beta_i +\kappa_{i-1}\beta_{i-1}||\leqno (1.4)$$
and if $\kappa_i\equiv 0$ on $I_{i-1}$, then $\beta$ is said to be of rank $i$.
$\kappa_i$ is then called the $i$--th Frenet curvature of $\beta$ and
$\beta_i$ the $i$--th Frenet vector. 

A curve in ${\mathbb E}^m$ is called a {\it W--curve\/} if its Frenet curvatures are
constant.

{\bf Lemma 1.1.} {\it If $\beta $ is a regular curve in ${\mathbb E}^m$ whose
first Frenet curvature is nowhere zero, then
the Laplace map of $\beta$ is regular.}

\demo Let $s$ denote an arclength parameter of the curve
$\beta$. Then $t={{d\beta}\over{ds}}$ is a unit tangent vector field of
$\beta$.
 Denote by $\kappa_i$ and $\beta_i$ the $i$--th Frenet curvature and $i$--th
 Frenet normal vector of $\beta$, respectively.
Then we have $$L(s)=-\beta ''=-\kappa_1 \beta_2.$$
Therefore, 
$$L_*(t)=-\beta'''={\kappa_1(s)}^2 \beta_1
-({\kappa_1(s)})'\beta_2-\kappa_1(s)\kappa_2(s)\beta_3\leqno(1.5)$$
where the last term occurs only when $m\geq 3$.
From (1.5), we obtain the Lemma. \sq

 \vskip.05in
{\bf Proposition 1.2.} {\it Let $\beta (s)$ be a curve
parametrized by arclengh $s$ in ${\mathbb E}^m$. Then the Laplace
transformation of  $\beta $ is homothetic if and only if
(1) $\beta$ has nonzero first Frenet curvature function
$\kappa_1 (s)$ and (2) the first and the second Frenet
curvature functions satisfy the relation: $\kappa_{1}^4 +
(\kappa_{1} ')^2 + \kappa_{1}^2 \kappa_{2}^2 = c$ for some
positive constant c.}

{\bf Proof.} From (1.5) we have
$$\<dL(t),dL(t)\>=\kappa_{1}^4 +
(\kappa_{1}')^2 + \kappa_{1}^2 \kappa_{2}^2.\leqno (1.6)$$
From (1.6) we obtain the Proposition. \sq

In particular, if $\beta$ is a planar curve, this
Proposition  yields the following.

{\bf Corollary 1.3.} {\it A planar curve $\beta (s) $ has
homothetic Laplace transformation if and only if its
curvature function $\kappa$  satisfies
$\kappa^4 +  (\kappa ')^2  = c$ for some positive
constant. \sq}

{\bf Remark 1.1.} We will come back to this situation in Chapter IV, section
3. \sq

From the  Proposition 1.2 we also have the following

{\bf Corollary 1.4} {\it Every $W$--curve
 in ${\mathbb E}^m$ has homothetic Laplace transformation.
}

\vskip.2in
\noindent {\bf \S2. Laplace maps of small rank.}
\vskip.1in

\def\nt{\tilde \nabla}

Let $x : M^n \rightarrow {\mathbb E}^m$  be an isometric
immersion. Denote by $dL$  the
differential  of the Laplace map $L :
M^n \rightarrow {\mathbb E}^m$ associated with the isometric immersion $x$.

The main purpose of this section is to study isometric
immersions whose Laplace map has rank   $< n=dim\, M$.

First we give the following result which classifies submanifolds whose
Laplace map is of rank 0.

{\bf Proposition 2.1.} {\it Let $x : M^n \rightarrow {\mathbb E}^m$ be
an isometric immersion. Then the Laplace image $L(M^n
)$ of the immersion $x$ is a point if and only if 
 $x$ is a minimal immersion.
}

{\bf Proof.} If the Laplace image is a point, then the mean
curvature vector field $H$ is a constant vector, say $c$.
Thus, by the Weingarten formula, we have
$${\tilde \nabla}_X H = -A_H X+D_X H=0$$
for any vector $X$ tangent to $M^n$. Thus $A_H =0$ which
implies that $H=c=0$. The converse is trivial. \sq

For hypersurfaces with singular Laplace map we have the following.
 
{\bf Proposition 2.2.} {\it
 Let $x : M^n \rightarrow {\mathbb E}^{n+1}$ be
an isometric immersion. Then the Laplace map of $x$
satisfies $rank \,(dL) \equiv k,$ for some $k$ with $ 0<k <n,$ if and only
if $M^n$ is foliated by $(n-k)$-dimensional submanifolds such that
\begin{itemize}
\item [(a)]  the second fundamental form $h$ of $M^n$
in ${\mathbb E}^{n+1}$, restricted to each leaf $N^{n-k}$, vanishes
identically and
\item [(b)] the mean curvature vector $H$ of $x$ is constant along
each leaf $N^{n-k}$.
\end{itemize}}

{\bf Proof.} Assume that the Laplace map of $x$
satisfies $rank \,(dL) \equiv k, \, 0<k <n.$ Then, for
any $k+1$ vectors $X_1 ,\ldots,X_{k+1} $ tangent to $M^n$,
we have $${\tilde \nabla}_{X_{1}} H
\wedge\ldots\wedge {\tilde \nabla}_{X_{k+1}}
H=0.\leqno(2.1)$$
Let $e_1 ,\ldots,e_n$ be an orthonormal frame such that
$e_i, i=1,\ldots,n,$ are principal vectors of $M^n$ in ${\mathbb E}^{n+1}$ with
principal curvatures $\kappa_i,i=1,\ldots,n,$ respectively. Then we have $${\tilde
\nabla}_{e_{i}} H=-\alpha \kappa_i e_i +(e_i
\alpha)e_{n+1},$$ where  $H=\alpha e_{n+1}$. Thus from
(2.1) and the fact that $\alpha\not=0$ by the previous Proposition, we may
assume that the orthonormal frame \onf satisfies $$\kappa_1 \ldots \kappa_k
\kappa_{r}=0,\;\;\;  \kappa_1 \ldots
\kappa_{k}(e_{r}\alpha)=0,\;\;\; r=k+1,\ldots, n.$$
Consequently, also taking into account that $rank(dL)\equiv k$, without
loss of generality, we may assume that one of the following two cases occurs:
$$\kappa_1 ,\ldots,\kappa_k \not=0,\;\;\;\kappa_{k+1}=
\ldots =\kappa_{n}=0 ,\;\;\; e_{k+1}\alpha =\ldots =e_n
\alpha =0,\leqno(2.2)$$
or
$$\aligned\kappa_1 ,\ldots,\kappa_{k-1} \not=0,\;\;\;\kappa_{k}=
\ldots =\kappa_{n}=0 ,\\
e_k \alpha \not=0,\;\;\;
e_{k+1}\alpha =\ldots =e_n \alpha =0.
\endaligned
\leqno(2.3)$$

From (2.2) and (2.3) we have
$$\nt_{e_{k+1}}H= \ldots =\nt_{e_{n}}H =0.\leqno(2.4)$$
 Let ${\Cal D} = \{Z\in TM^n : {\tilde \nabla}_Z H=0\}$ and let $Z, W$ be any
two vector fields in ${\Cal D}$. Then from the Codazzi equation, we have
$$B([Z,W])=B(\nabla_Z W)-B(\nabla_W Z)= (\nabla_W B)Z
-(\nabla_Z B)W=0,\leqno(2.5)$$
where $B=A_{e_{n+1}}$.
From the definition of ${\Cal D}$ we also have
$$[Z,W]\alpha =ZW\alpha -WZ\alpha =0.\leqno(2.6)$$
From (2.5) and (2.6) we conclude that $\Cal D$ is
completely integrable. So $M^n$ is foliated by
$(n-k)$-dimensional submanifolds. Furthermore, from the
definition of $\Cal D$, we see that each leaf of
$\Cal D$ satisfies conditions (a) and (b). The converse is
clear. \sq 
 
In particular, if $rank (dL)\equiv 1$, we have the following.

{\bf Proposition 2.3.} {\it Let $x : M^n \rightarrow
{\mathbb E}^{n+1}$  be an isometric immersion. Then the
Laplace map of $x$ satisfies $rank (dL)
\equiv 1$ if and only if $M^n$ is locally a hypercylinder on a
3--regular curve.}

{\bf Proof.} Assume that $rank \, (L_*)\equiv 1$. Then, from the
proof of \prop 2.2, we see that there eixsts an orthonormal
frame $e_1 ,\ldots, e_n$ of principal vectors such that
their corresponding principal curvatures satisfy the
following conditions:
$$\kappa_1 =n\alpha,\;\; \kappa_2 =\ldots =\kappa_n
=0,\;\;\; e_2 \alpha=\ldots =e_n \alpha =0.\leqno(2.7)$$
Let ${\Cal D}=\{ Z\in TM :\nt_Z H=0 \}.$ Then $\Cal D$ is
integrable from the proof of \prop 3.2. Let $X$ be any
vector field in ${\Cal D}^{\perp}=Span \{ e_1 \}$ and $Z
\in {\Cal D}$. Then we have
$$
\nabla_Z (BX) =(\nabla_Z B)X+B(\nabla_Z
X)=(\nabla_X A)Z + B(\nabla_Z X)$$ $$ =B(\nabla_Z
X)-B(\nabla_X Z)=B([Z,X]).
$$
This implies that $\nabla_Z e_1 \in {\Cal D}$.
Consequentyly,  each leaf of $\Cal D$ is
totally geodesic in ${\mathbb E}^{n+1}$; and hence it is an open
portion of a linear subspace of ${\mathbb E}^{n+1}$. Since ${\Cal
D}^{\perp}$ is of rank one, ${\Cal D}^{\perp}$ is trivially
integrable. Because the second fundamental form of $M^n$
in ${\mathbb E}^{n+1}$ satisfies $h({\Cal D},{\Cal D}^{\perp})=\{
0\}$, a lemma of Moore implies that $M^n$ is locally a
hypercylinder on a planar curve. Moreover, since $rank(dL)\equiv 1$, this
curve is 3-regular.
The converse is easy to verify. \sq

For surfaces in a Euclidean space of higher codimension, we 
have the following.
\def\trace{\hbox{trace}}

{\bf Proposition 2.4.} {\it  Let $x : M^2
\rightarrow {\mathbb E}^m$ be an isometric immersion. If the Laplace
map of $x$ satisfies $rank \,(dL) \equiv 1$, then $M^2$ is
non-positively curved, i.e., the Gauss curvature $K$ of
$M^2$ is $\leq 0$.

In particular, if
$K=0$, then $M^2$ is a cylinder on a 3--regular curve.}

{\bf Proof.}  Let $x : M^2
\rightarrow {\mathbb E}^m$ be an isometric immersion whose Laplace
map satisfies $rank \,(dL) \equiv 1$. Then we have
$${\tilde \nabla}_{X_{1}} H
\wedge {\tilde \nabla}_{X_{2}}
H=0,\leqno(2.8)$$
for any vectors $X_1 ,X_2$ tangent to $M^2$. 

Let $U=\{p\in M^2:H(p)\not=0\}.$ Then $U$ is an open subset of $M^2$. 
It is easy to see from the Gauss equation that the Gauss curvature $K\leq 0$
on $M-U$. 

On $U$, we choose an
orthonormal frame $e_1 ,\ldots, e_m$, such that,
restricted to $U$, $e_1 $ and $e_2$ are tangent to
$M^2$ and $H=\alpha e_3$. Furthermore, we may also assume
that $$A_H e_1 =\kappa_1 e_1\quad \hbox{and}\quad A_H e_2 =\kappa_2 e_2.$$
From (2.8) we have 
$$(-\alpha \kappa_1 e_1 +(e_1 \alpha)e_3 +\alpha
D_{e_{1}} e_3 )\wedge (-\alpha\kappa_2 e_2 +(e_2
\alpha)e_3 +\alpha D_{e_{2}} e_3 )=0.$$
From this we may assume without loss of generality the
following:
$$\kappa_1 \not= 0,\;\;\kappa_2 =0,\;\;  e_2 \alpha
=  D_{e_{2}}e_3 =  0.\leqno(2.9)$$
Since $e_3$ is parallel to the mean curvature vector
$H$, we have $\,\trace\, (A_r ) = 0$ for $r=4,\ldots ,m$.
Thus, by (2.9),  the Gauss curvature $K$ of
$M^2$ is $\leq 0$. 

In particular, if $K \equiv 0$, then we have 
$$A_3 = \begin{pmatrix} 2\alpha & 0\\ 0 & 0 \end{pmatrix}
,\;\;\; A_4 =\ldots =A_m = 0.\leqno(2.10)$$ 

From
(2.10) have
$$h(e_1 ,e_1 )=2\alpha e_3\;\;\; h(e_1 ,e_2 )=h(e_2 ,e_2
)=0.\leqno(2.11)$$
From (2.9), (2.11)  we may
obtain
 $$({\bar \nabla}_{e_{1}} h)(e_1 ,e_2
)=-\omega^{1}_{2}(e_1 )h(e_1 ,e_1),\;\; 
({\bar \nabla}_{e_{2}} h)(e_1 ,e_1
)=0.\leqno(2.12)$$
By (3.10) and the equation of codazzi we obtain
$$\omega_{1}^{2} = f\omega^2,\leqno(2.13)$$ for some
function $f$ on $M^2$.
By taking the exterior derivative of $\omega_{2}^{3}=0$
and by applying (2.8) and (2.13), we may get $f=0$. This
implies $\nabla e_1 =\nabla e_{2}=0$. Hence, by
combining this with (2.11), we see that each integral
curve of $e_2$ is an open portion of a straight line in
$E^m$ and so $M^2$ is a cylinder on a 3--regular curve. \sq

{\bf Remark.} There do exist  negatively-curved surfaces
in ${\mathbb E}^m$  whose Laplace maps $L$ satisfy $rank(dL)\equiv
1$. \sq

From Proposition 2.2 and Proposition 2.4, we obtain the
following  

{\bf Corollary 2.5.} {\it Let   $x :
M^n \rightarrow {\mathbb E}^{m}$  be an isometric immersion of a
compact Riemannian manifold. Then
\begin{itemize}
\item[(1)] if $m=n+1$,  the Laplace map is
regular (i.e., it is of maximal rank) on a non-empty open
submanifold of $M^n$ and

\item[(2)]  if $n=2$ and $M^2$ is diffeomorphic to a
sphere, a real projective plane, a Klein bottle or a
torus, then the Laplace map is regular on a non-empty open
subset. \sq
\end{itemize}}

\vskip.2in 
\noindent {\bf \S3. Ruled surfaces in ${\mathbb E}^m$.}
\vskip.1in

In this section we study
the Laplace map of ruled surfaces in ${\mathbb E}^m$.

{\bf Theorem 3.1.} {\it Let $M^2$ be a ruled surface in a Euclidean
$m$--space ${\mathbb E}^m$. Then the restriction of $dL$ in the direction of
the rulings  vanishes if and only if either $M^2$ is an open portion of a
helicoid or it is an open portion of a cylinder over a curve.
}

{\bf Proof.}  We consider the two cases separately.

\noindent {\bf Case 1.} $M$ is a cylinder.

Suppose that the surface $M$ is a cylinder over a curve $\gamma$ in an
affine hyperplane ${\mathbb E}^{n-1}$, which we can choose to have the equation
$x_m=0$. Assume that $\gamma$ is parametrized by its arc length
$s$. Then a parametrization  of $M$ is given by
$$x(s,t)=\gamma(s)+te_m.\leqno (3.1)$$

The Laplacian $\Delta$ of $M$ is given in terms of $s$ and $t$ by
$$\Delta=-{{\partial^2}\over{\partial s^2}}-{{\partial^2}\over{\partial t^2}}
.\leqno (3.2)$$ 
Thus the Laplace map $L$ of the cylinder is given by 
$L(s,t)=-\gamma''(s)$; and hence $dL({{\partial}\over{\partial t}})=0$
{\it i.e.,} the restriction of  $dL$ in the direction of the rulings 
vanishes  identically.

\noindent {\bf Case 2:} $M$ is not cylindrical.

If the ruled surface $M$ is not cylindrical, we can decompose $M$ into
open pieces such that on each piece we can find a parametrization
$x$ of the form:
$$x(s,t)=\alpha(s)+t\beta(s)\leqno (3.3)$$
where $\alpha$ and $\beta$ are curves in ${\mathbb E}^m$ such that
$$\<\alpha',\beta\> =0,\quad \<\beta,\beta\>=1,\quad \<\beta',\beta'\>=1.$$

We have $x_s=\alpha'+t\beta'$ and $x_t=\beta$. We define functions
$q,u$ and $v$ by
$$q=||x_s||^2=t^2+2ut+v,\quad u=\<\alpha',\beta'\>,\quad
v=\<\alpha',\alpha'\>. \leqno (3.4)$$
The Laplacian $\Delta$ of $M$ can be expressed as follows:
$$\Delta= -{{\partial^2}\over{\partial t^2}} -
{1\over q} {{\partial^2}\over{\partial s^2}}
+{1\over 2}{{\partial q}\over{\partial s}}
{1\over {q^2}}{{\partial}\over{\partial s}}
-{1\over 2}{{\partial q}\over{\partial t}}
{1\over q}{{\partial}\over{\partial t}}. \leqno (3.5)$$

From (3.3) and (3.5) we see that the Laplace map is given by
$$L(s,t)={1\over {2q^2}}\{-2q\alpha''(s)+q_s
\alpha'(s)-2tq\beta''(s)+q_s\beta'(s)
-q_tq\beta(s)\}$$
where $q_s={{\partial q}\over {\partial s}},\,
 q_t={{\partial q}\over {\partial t}}$. This implies
 that $dL({{\partial q}\over {\partial t}})=0$ if and only if
$$\aligned 2qq_t\alpha''+(qq_{st}-2q_sq_t)\alpha'+2q(tq_t-q)\beta''\\
+(qq_{st}-2q_sq_t)\beta'+q(q_t^2-qq_{tt})\beta=0,\endaligned
\leqno (3.6)$$
where $q_{ss}={{\partial^2 q}\over{\partial s^2}},
q_{st}={{\partial^2 q}\over{\partial s\partial t}}$ and
$q_{tt}={{\partial^2 q}\over{\partial t^2}}.$ From (3.4) and (3.6)
we  obtain $$\aligned
(\beta+\beta'')t^4+2(2u\beta+u\beta''-\alpha'')t^3\hskip1in \\
-3(2u\alpha''-2u^2\beta+u'\alpha'+u'\beta')t^2\hskip1in\\
-\{2(v'+uu')\alpha'+2uv\beta''+2(v'+uu')\beta'\hskip.4in\\-4u^3\beta
+2(2u^2+v)\alpha''\}t\hskip.3in\\
-\{(2uv'-u'v)\alpha'+v^2\beta''+2uv\alpha''\hskip1.3in\\ +(2uv'-u'v)\beta'
+(v^2-2u^2v)\beta\}=0.\hskip.7in\endaligned\leqno (3.7)$$
From (3.7) we obtain
$$\beta ''+\beta =0,\leqno (3.8)$$
$$ \alpha ''=u\beta,\leqno (3.9)$$
$$u'(\alpha'+\beta ')=(uu'+v')(\alpha'+\beta')=0,\leqno (3.10)$$
$$(2uv'-u'v)(\alpha'+\beta')=0.\leqno (3.11)$$

If $\alpha'+\beta'=0$, then $\alpha(s)=-\beta(s)+c$ for some constant vector
$c\in {\mathbb E}^m$.
On the other hand, (3.8) implies that $\beta$ is a unit speed curve of 1--type.
Hence, $\beta$ is a plane circle. Consequently, by  
using $\alpha(s)=-\beta(s)+c$, we may conclude that $x(s,t)=\alpha(s)
+t\beta(s)$ is an open portion of a plane which is a special case of cylinder.

If $\alpha'+\beta'\not=0$, then (3.10) yields
$$v'=-uu',\quad u'v=-2uv'.$$
Thus $(2u^2+v)u'=0$. Since $v=\<\alpha',\alpha'\>>0$, 
$u$ is a constant. Hence 
$$\alpha''=-\lambda \beta\leqno (3.12)$$ for some nonzero constant $\lambda$.
From (3.8), (3.12) and $\<\beta,\beta\>=\<\beta',\beta\>=1$,
 we may choose a Euclidean
coordinate system on $E^m$ such that $\alpha$ and $\beta$ are given
respectively by
$$\beta(s)=(\cos s,\sin s,0,\ldots,0).\leqno (2.13)$$
$$\alpha(s)=(\lambda \cos s +c_1s,\lambda \sin s+c_2s,
c_3s,0,\ldots,0)\leqno(3.14)$$
for some constants $c_1,c_2,c_3$. From (3.13) and (3.14) we see that
the ruled surface $x(s,t)=\alpha(s)+t\beta(s)$ is an open portion of a
helicoid. Conversely, since a helicoid is a minimal surface, its Laplace map
$L$ is a constant map. Thus $dL=0$ identically. \sq

In the remaining part of this section we investigate the Laplace
map of flat ruled surfaces in ${\mathbb E}^m$. It is well known that
 flat ruled surfaces in ${\mathbb E}^m$ are ``in general'' cylinders,
cones or tangential developables of curves. As we already have seen,
the differential $dL$ of the Laplace map of a cylinder in ${\mathbb E}^m$  
has rank $\leq 1$. 

Now, we study the Laplace maps of cones and tangential developables.

{\bf Proposition 3.2.} {\it The Laplace map of a cone in ${\mathbb E}^m
(m\geq 3)$ is a cone.}

{\bf Proof.} Let $x:M\rightarrow {\mathbb E}^m$ be a cone in ${\mathbb E}^m$. Without
loss of generality, we may assume the cone has its vertex at
the orign and it is parametrized by
$$x(t,s)=t\beta(s),\quad |\beta|=|\beta'|=1.\leqno(3.15)$$
Put 
$$e_1={1\over t}{{\partial}\over{\partial s}}\quad
e_2={\partial\over{\partial t}}.\leqno(3.16)$$
Then, by direct computation, the second fundamental form $h$ of
$M$ in ${\mathbb E}^m$ satisfies
$$h(e_1,e_1)={1\over t}\beta''-{1\over t}\<\beta'',\beta\>\beta,
 \quad h(e_2,e_2)=0.\leqno(3.17)$$
From (3.17) we see that the Laplace map of $x$ is given by
$$L(t,s)={1\over t}(\<\beta''(s),\beta(s)\>\beta(s)-\beta''(s)).
$$
On the other hand, (3.15) yields $\<\beta'',\beta\>=-1$. Thus
$$L(t,s)=-{1\over t} (\beta(s)+\beta''(s)),\leqno(3.18)$$
which implies that the Laplace map $L$ is also a cone with
vertex at the origin. \sq

{\bf Proposition 3.3.} {\it The Laplace map of a tangential
developable surface in ${\mathbb E}^m (m\geq 3)$ is a cone.}

{\bf Proof.} Assume $M$ is a tangential developable surface in ${\mathbb E}^m$ given
by the tangent lines of a unit speed curve $\beta(s)$ in ${\mathbb E}^m$. Then
$M$ is parametrized by
$$x(s,t)=\beta(s)+t\beta'(s),\quad |\beta'(s)|=1.\leqno(3.19)$$
Let $\kappa_1$ denote the first Frenet curvature of $\beta$ in 
${\mathbb E}^m$.  

Put
$$e_1={1\over{t\kappa_1}}\Big({{\partial}\over{\partial s}}
-{{\partial}\over{\partial t}}\Big),\quad
e_2={\partial\over{\partial t}}.\leqno(3.20)$$
Then, by a direct computation, the second fundamental form $h$
of $M$ in ${\mathbb E}^m$ satisfies
$$h(e_1,e_1)={1\over{t\kappa_1}}\beta_3,\quad h(e_2,e_2)=0,\leqno(3.21)$$
where $\beta_3$ is the third Frenet vector.
(3.21) implies that the Laplace map of $M$ is given by
$$L(s,t)=-{1\over {t\kappa_1}}\beta_2(s).\leqno(3.22)$$
Therefore, the Laplace map is a cone with vertex at the origin. \sq

\vfill\eject

\noindent {\bf Chapter IV:  HOMOTHETIC LAPLACE
TRANSFORMATIONS.}
\vskip.2in

\noindent {\bf \S1. Some general results.}
\vskip.1in

Let $x : M \rightarrow {\mathbb E}^m$ be an isometric immersion
of an $n$-dimensional connected Riemannian manifold
$M$ into a Euclidean $m$-space. Denote by  $L :M^{n} \rightarrow {\mathbb E}^m$
 the {\it Laplace map\/} and by $L(M^n )$  the
{\it Laplace image\/} of the immersion $x$.  
Recall that  the transformation ${\Cal L} :M \rightarrow
L(M )$ from $M^n$ onto its Laplace image $L(M )$ is called 
 the {\it Laplace transformation\/} of the immersion  $x$.

The main purpose of this section is to give some general properties
concerning isometric immersions $x : M^n \rightarrow {\mathbb E}^m$ with
homothetic Laplace transformation. 

First we give the following general results.

{\bf Lemma 1.1.} {\it Let $x : M \rightarrow {\mathbb E}^m$ be an isometric immersion
of an n--dimensional Riemannian manifold into ${\mathbb E}^m$. Then the Laplace
transformation ${\Cal L}:M\rightarrow L(M)$ is homothetic 
if and only if
$$\<A_HX,A_HY\>+\<D_XH,D_XH\>=c^2\<X,Y\>\leqno (1.1)$$
holds for all vectors X,Y tangent to $M$, where $c$ is a positive constant.
}

{\bf Proof.} Because the Laplace map is given by $L(p)=(\Delta x)(p)=-n
H(p)$, the differential of the Laplace map satisfies
$$dL(X)=nA_HX-nD_XH.\leqno (1.2)$$
From (1.2) we obtain
$$\<dL(X),dL(Y)\>=\<A_HX,A_HY\>+\<D_XH,D_HY\>\leqno (1.3)$$
which implies the Lemma. \sq

In the following, by $S^{m-1}(r)$ we mean the hypersphere of ${\mathbb E}^m$ centered
at the origin and with radius $r$.

{\bf Lemma 1.2.} {\it
If  $$x:M\rightarrow S^{m-1}(r)\subset {\mathbb E}^m$$ is an isometric
immersion  which immerses $M$ into the hypersphere $S^{m-1}(r)$  as
a minimal submanifold, then $x$ has homothetic Laplace transformation.
}

{\bf Proof.}  If  $$x:M\rightarrow S^{m-1}(r)\subset {\mathbb E}^m$$
 is an isometric immersion 
which immerses $M$ into the hypersphere $S^{m-1}(r)$  as
a minimal submanifold, then
the Laplace map $L$ of $x$ is given by
$L={n\over {r^2}}x$. Thus, $dL(X)={n\over{r^2}}X$ for any vector $X$
tangent to $M$. Hence, $x$ has homothetic Laplace transformation. \sq

{\bf Lemma 1.3.}  Let $x : M \rightarrow {\mathbb E}^m$ and
$y : N \rightarrow {\mathbb E}^r$ be two isometric immersions whose
Laplace maps are given by $L_M : M \rightarrow {\mathbb E}^m$ and
$L_N :N\rightarrow {\mathbb E}^r$, respectively. Then
\begin{itemize}

\item[(1)] the Laplace map of
the product immersion $(x,y): M\times N \rightarrow
{\mathbb E}^{m+r}$ is given  by $(L_M ,L_N)
: M\times N \rightarrow {\mathbb E}^{m+r}$;

\item[(2)]  the Laplace image of
the product immersion $(x,y): M\times N \rightarrow
{\mathbb E}^{m+r}$ is given by $L_M (M) \times L_N (N)$; and 

\item[(3)] the
Laplace transformation ${\Cal L}_{M\times N} :M\times N
\rightarrow  L_{M\times N}(M\times N)$ of the
product immersion is homothetic if and only if the
Laplace transformations of $x$ and $y$ are 
homothetic transformations with the same constant homothetic
factor.
\end{itemize}
}

{\bf Proof.} This Lemma follows easily from the fact that the
Laplace operator of the Riemannian product $M\times N$ of
two Riemannian manifolds $M$ and $N$ is given by 
$$\Delta_{M\times N} =\Delta_M \times \Delta_N
.\quad \square\leqno(1.4)$$

 Lemma 1.2 and Lemma 1.3 show that there exist ample examples of
submanifolds with homothetic Laplace transformation.

 If $x : M \rightarrow {\mathbb E}^m$ is an isometric immersion of a Riemannian
manifold $(M ,g)$ with Riemannian metric $g$ and $c$ a positive
number, then the immersion $x^c$ defined by $(x^c )(p)
 = cx(p)$, for $p\in M$, is an isometric immersion from
the Riemannian manifold $(M ,c^2 g)$ into ${\mathbb E}^m$. 

{\bf Lemma 1.4.} {\it  Let  $x : M \rightarrow {\mathbb E}^m$
 be an isometric immersion. Then the Laplace
map $L^c :M \rightarrow {\mathbb E}^m$ of the isometric
immersion $x^c$ is given by $L^c (p)=c^{-2} L(p)$ for any
point $p\in M$.
}

{\bf Proof.} Follows from the fact that the Laplace operator
$\Delta^c$ of $(M ,c^2 g)$ and the Laplace operator
$\Delta$ of $(M ,g)$  are related by $$\Delta^c =c^{-2}
\Delta.\quad \square$$

Let $y_i : M \rightarrow {\mathbb E}^{m_{i}}$,
$i=1,\ldots,k,$ be $k$ maps of $M$ into ${\mathbb E}^{m_i}$, respectively, and
let $a_1,\ldots,a_k$ be $k$ positive numbers. Then 
 $$y(p) = (a_1 y_1 (p),\ldots, a_k y_k (p)),\quad p \in M,$$ 
 is a map of
$M$ into ${\mathbb E}^{m_1 +\ldots +m_k }$, which is called a
{\it diagonal map\/} of  $y_1, \ldots, y_k$.  
In particular, if $y_1,\ldots,y_k$ are isometric immersions and 
$a_1,\ldots,a_k$ satisfy $a_1^2+\ldots +a_k^2=1$, then the diagonal map
$y$ of $y_1,\ldots,y_k$ is an isometric immersion, called a {\it diagonal
immersion\/} of $y_1,\ldots,y_k$.

{\bf Lemma 1.5.} {\it Let $x : M \rightarrow {\mathbb E}^m$ be
a diagonal immersion of k isometric immersions. Then
 the Laplace map $L : M\rightarrow {\mathbb E}^m$ of $x$ is  a diagonal map.
}

{\bf Proof.}  From the definition of diagonal immersion, we have
$$(\Delta x)(p)=(a_1(\Delta x_1)(p),
\ldots,a_k(\Delta x_k)(p))\leqno (1.5)$$
which implies that the Laplace map $L$ of the diagonal immersion $x$
is related to the Laplace map $L_i$ of $x_i$ by
$$L(p)=(a_1L_1(p),\ldots,a_kL_k(p)), \quad p\in M.\leqno (1.6)$$
Thus, the Laplace map of a diagonal immersion is also a diagonal map.
\sq

If $M$ is a compact (connected) Riemannian homogeneous
manifold, let $G=I_{o}(M)$ be the identity component of
the group of all isometries of $M$. $G$ is a compact Lie
group acting on $M$ transitively.
A map $x :M \rightarrow E^m$  from $M$ into ${\mathbb E}^m$ is said
to be {\it equivariant\/} if there exists a Lie
homomorphism $\varphi$ of $G$ into the isometry group
$I(E^m)$ of $E^m$ such that
 $$x(\sigma (p))=\varphi(\sigma)(x(p))$$
 for any $\sigma \in G$ and $p \in M$.

For the Laplace map of an equivariant isometric immersion
of a compact homogeneous space, we have the following general result.

{\bf Lemma 1.6.} {\it Let $x :M \rightarrow {\mathbb E}^m$ be
an equivariant isometric immersion of a compact 
Riemannian homogeneous manifold $M$ into  ${\mathbb E}^m$. Then the
Laplace map $L : M \rightarrow {\mathbb E}^m$ of $x$ is also
equivariant.}

{\bf Proof.}  This lemma follows from the fact that
 the Laplacian operator $\Delta$ is a Riemannian
invariant. \sq

Lemma 1.6 implies the following. 

{\bf Lemma 1.7.} {\it  Let $x :M \rightarrow {\mathbb E}^m$ be an equivariant
isometric immersion of a compact Riemannian homogeneous
manifold $M$ into  ${\mathbb E}^m$. Then the immersion $x$ has
constant mean curvature function and the Laplace map $L :
M \rightarrow {\mathbb E}^m$ of $x$ is spherical. \sq}

For equivariant immersions we also have the following

{\bf Theorem 1.8.} {\it Let $x : M^n \rightarrow {\mathbb E}^m$ be
an equivariant isometric immersion of a compact
irreducible Riemannian homogeneous manifold. Then the
associated Laplace transformation ${\Cal L} : M^n
\rightarrow L(M^{n})$ is a homothetic transformation.}

{\bf Proof.} 
 Let $x:M\rightarrow {\mathbb E}^m $ be an
equivariant isometric immersion of a compact connected
 Riemannian homogeneous space $M$ into ${\mathbb E}^m$.
Without loss of generality we may assume the immersion is
full. From Lemma 1.7, the immersion $x$ is spherical.
Thus $x(M)$ is contained in a hypersphere $S^{m-1}$ of
${\mathbb E}^m$. Without loss of generality, we may also assume that
$S^{m-1}$ is centered at the origin of ${\mathbb E}^m$. Then
 there is a Lie homomorphism
$\phi : G\rightarrow SO({\mathbb E}^{m})$ such that 
$x(g(p))=\phi(g)(x(p))$ for every $g\in G$ and $p\in M$.
Because $(\phi,{\mathbb E}^{m})$ is a representation of the compact
Lie group $G$, $(\phi,{\mathbb E}^{m})$ is the direct sum of some
irreducible subrepresentations
$(\phi_{1},E_{1}),\ldots,(\phi_{k},E_{k})$ such that ${\mathbb E}^m$
is the Euclidean direct sum $E_{1}\oplus\ldots\oplus
E_{k}$ of $E_{1},\ldots,E_{k}$. Let $x_{i}$ denote the
$E_{i}$-component of $x$. Then we have
$$x_{i}(g(p))=\phi_{i}(g)(x_{i}(p)),\;\;\;g\in G,\;\;
p\in M,\;\;\;i=1, \ldots,k,\leqno(1.7)$$
where $E_{1},\ldots,E_{k}$ are mutually orthogonal in
${\mathbb E}^{m}$. 

  Now we  claim that each $x_{i}$ is a 1-type map, that is,
$\Delta x_{i}=\lambda_{i}x_{i}$, $i\in \{1,\ldots,k\},$
for some real numbers $\lambda_{i}$. In order to do so, we
choose a fixed point 
 $o\in M$. Denote by $K$ the
isotropy subgroup of $G$ at $o$. Then $M$ can be
identified with $G/K$ in a natural way. Consider a
biinvariant Riemannian metric on the compact Lie group $G$
such that the projection $\pi: G\rightarrow M=G/K$ is a
Riemannian submersion. Let $e_{1},\ldots,e_{N}$ be any
orthonormal basis of the Lie algebra 
${\frak g} =T_{e}G$ of $G$, where $e$ is the identity element
of $G$. 

 For each $i\in \{1,\ldots,k\}$,  denote also by
$\phi_{i}$ the  homomorphism  $ {\frak g} \rightarrow 
{\frak so}(E_{i})$ induced from $\phi_{i}:G\rightarrow
E_{i}$, where ${\frak so}(E_{i})$ is the Lie algebra of
$SO(E_{i})$. Then each $\phi_{i}(e_{a})$ is a
skew-symmetric linear transformation of $E_i$ and
$\sum_{a=1}^{N} \phi_{i}(e_{a})^2$ is a symmetric linear
transformation. Let $Ad: G\rightarrow GL({\frak g})$ be the
adjoint representation of $G$. Since $ad(g)(h)=ghg^{-1}$
and $\phi_{i}$ is a Lie homomorphism, we have
$$\phi_{i}(g)\phi_{i}(X)\phi_{i}(g^{-1})=\phi_{i}(Ad(g)X)$$
for any $X\in {\frak g}, \,g\in G$ and $i\in
\{ 1,\ldots,k\}$. Therefore we find
$$\phi_{i}(g)(\sum_{a=1}^{N}
\phi_{i}(e_{a})^{2})\phi_{i}(g^{-1}) =\sum_{a=1}^{N}
\phi_{i}(Ad(g)e_{a})^{2}=\sum_{a=1}^{N}
\phi_{i}(e_{a})^{2}$$
for any $g\in G$. This shows that $\sum_{a=1}^{N}
\phi_{i}(e_{a})^2$ lies in the centralizer of
$\phi_{i}(G)$. Since the representation $(\phi_{i},E_{i})$
is irreducible, Schur's lemma in representation theory
implies that 
$$\sum_{a=1}^{N} \phi_{i}(e_{a})^{2} = - \lambda_{i}
I_{i},\leqno(1.8)$$ for some constants $\lambda_{i}$, where
$I_i$ is the identity transformation on $E_i$. On the
other hand, it is known that the Laplacian of $x_{i}$ is
given by  $$\Delta
x_{i}(p)=-\sum_{a=1}^{N}{d^{2}\over
{dt^{2}}}x_{i}(exp\,\,te_{a}){
|}_{t=0}=-\sum_{a=1}^{N} \phi_{i}(e_{a})^{2}(x_{i}(p)).
\leqno(1.9)$$ Therefore, by (1.8) and (1.9), we have
$$\Delta x_i=\lambda_iI_i,\quad i=1,\ldots,k,\leqno (1.10)$$
for some constants $\lambda_1,\ldots,\lambda_k$.
In particular, this shows that $x$ is a diagonal immersion of 
1--type maps $x_i:M\rightarrow E_i,\,i=1,\ldots,k$. Therefore, the
Laplace map of the equivariant immersion $x$ is given by
$$L(p)=\lambda_1x_1(p)+\ldots+\lambda_kx_k(p).\leqno (1.11)$$
Since $E_1,\ldots,E_k$ are mutually orthogonal in ${\mathbb E}^m$, (1.11)
implies $$\<dL,dL\>=\lambda_1^2 \<dx_1,dx_1\>+\cdots
+\lambda_k^2\<dx_k,dx_k\>.\leqno (1.12)$$
If $M$ is an irreducible homogenous Riemannian manifold, the linear isotropy 
representation
is irreducible. So, in this case, each $\<dx_i,dx_i\>$ is a constant
multiple of the  original metric on $M$. Thus, by (1.12),  the Laplace
transformation of the equivariant immersion $x$ is a homothetic transformation.
\sq

For compact submanifolds with homothetic Laplace
transformation, we  have the following

{\bf Proposition 1.9.}  {\it Let $x
: M^n \rightarrow {\mathbb E}^{m}$ be an isometric immersion of a
compact Riemannian manifold $M^n$ into ${\mathbb E}^m$. Then we
have:
\begin{itemize}
\item[(1)] the Laplace map $L :M^n \rightarrow {\mathbb E}^m$
has center of gravity at the origin 0 of ${\mathbb E}^m$ and 

\item[(2)] if the Laplace transformation  ${\Cal L} : M^n \rightarrow L(M^{n})$
 of $x$ is homothetic, or more generally volume-element
preserving, then the center of gravity of the Laplace image
$L (M^n )$ (with respect to the induced metric) in ${\mathbb E}^m$ is
the origin 0 of ${\mathbb E}^m$.
\end{itemize}
}

{\bf Proof.}  Since $M$ is compact and $L=\Delta x$, Hopf's lemma implies that
the center of gravity of the Laplace map is the origin. This proves
statetement (1). Statement (2) follows from statement (1) easily.
\sq

{\bf Remark 1.1.}  Proposition 1.9 shows that not every submanifold 
in ${\mathbb E}^m$ can be realized as the Laplace image
of some submanifolds in ${\mathbb E}^m$, in particular, {\sl every compact
submanifold in ${\mathbb E}^m$ cannot be realized as the Laplace image of any
submanifold if its center of gravity differs from the origin.} \sq

The following results establish some relations between
Laplace transformations and the notion of submanifolds of
finite type.

{\bf Proposition 1.10.} {\it Let $x :M^n \rightarrow {\mathbb E}^{m}$  be an
isometric
 immersion. If $x$ has homothetic  Laplace transformation  ${\Cal L} : M^n
\rightarrow L(M^{n})$, then 
\begin{itemize}
\item[(1)] if  the immersion $x$ is of finite type, then
the Laplace map $L: M^n \rightarrow {\mathbb E}^m$
is of finite type;

\item[(2)] if $M$ is compact, then the Laplace map $L$ is of $k$--type if
and only if  the immersion $x$ is of $k$--type;

\item[(3)] if the immersion $x$ is of finite type, then
the Laplace map is of non-null finite type; in particular, if $x$ is of
non--null $k$--type, then $L$ is of non--null $k$--type; and if $x$ is
of null $k$--type, then $L$ is of non--null $(k-1)$--type.
\end{itemize}
}

{\bf Proof.}  If the Laplace transformation is homothetic, then the
Laplacian operator $\Delta_L$ of the Laplace image is related with the
Laplacian operator of $M$ by $\Delta_L=c\Delta$ for some positive constant
$c$. 

(1) If the immersion $x$ is of finite type, then we have
$$x=x_1+\cdots+x_k,\quad \Delta x_i=\lambda_i x_i,\quad i=1,\ldots,k\leqno(1.13)$$
for some natural number $k$ and real numbers $\lambda_1,\ldots,\lambda_k$
and non--constant maps $x_1,\ldots,x_k$.
From (1.10) we get
$$L=L_1+\cdots+L_k,\leqno (1.14)$$
where $L_i=\lambda_ix_i$. Because $\Delta_LL_i=c\Delta L_i
=c\lambda_i L_i$, (1.14) implies that the Laplace map $L$ of $x$ is also of
finite type.

(2) Suppose $M$ is compact and 
the Laplace map $L$ of $x$ is of $k$--type. Let $Q(t)$ be the minimal
polynomial of $L$. Then $\deg Q=k$ and $Q(\Delta_L)L=0$. 
Because $L=\Delta x$ and 
$\Delta_L=c\Delta$ for some positive constant $c$, there exists a
polynomial $\bar Q(t)$  such that $\deg \bar Q=\deg Q$ and
 $P(\Delta)x=0$ where $P(t)=t\bar Q$. Therefore, by Theorem 3.1 of Chapter 1, 
the immersion $x$ is of $\ell$--type with $\ell\leq k+1$. 

If $x$ is of $(k+1)$--type, then 
$$x=x_1+\cdots+x_{k+1},\quad \Delta x_i=\lambda_i x_i,\quad i=1,\ldots,k+1
\leqno(1.15)$$
for non--constant maps $x_1,\ldots,x_{k+1}$. Because, $M$ is compact, all of
$\lambda_1,\ldots,$ $\lambda_{k+1}$ are non--zero. This implies that $L$
 is of $(k+1)$--type which is a contradiction. Similarly, if $x$ 
is of $\ell$--type
with $\ell<k$, then $L$ is of $\ell$--type with $\ell<k$. 
The remaining part of statement (2) follows from the proof of statement (1).

(3)  If $x$ is of non--null $k$--type whose spectral decomposition
is given by (1.13), then (1.14) is the 
spectral decomposition of the Laplace map $L$. Since $\Delta_LL_i=c\lambda_i L_i$
and $c,\lambda_i$ are nonzero constants, (1.14) shows that $L$ is of non--null 
$k$--type. If $x$ is of null  $k$--type, then one of $\lambda_1,\ldots,\lambda_k$ is zero.
Without loss of generality, we may assume $\lambda_1=0$.
Then the spectral decompostion of the Laplace map is given by
$$L=L_2+\cdots+L_k,\quad \Delta_L L_i=c\lambda_iL_i,\quad i=2,\ldots,k$$
which implies that $L$ is non--null $(k-1)$--type. This proves statement (3).
\sq

 {\bf Proposition 1.11.} {\it Let $x : M \rightarrow {\mathbb E}^{m}$  be
an isometric immersion with homothetic  Laplace
transformation  ${\Cal L} : M \rightarrow L(M)$. Then
\begin{itemize}
\item[(1)] if $x$ is linearly independent, then the
Laplace map is also linearly independent; and
\item[(2)] if $x$ is orthogonal, then the Laplace
map is also orthogonal.
\end{itemize}}
{\bf Proof.}  Let $x: M
\rightarrow {\mathbb E}^m$ be an immersion of $k$-type and let $$x =
c+x_{1}+\ldots +x_{k},\hskip.3in \Delta
x_{i}=\lambda_{i}x_{i},\hskip.2in \lambda_{1}<\ldots
<\lambda_{k}\leqno(1.16)$$ be the spectral decomposition of
the immersion $x$, where $c$ is  a constant vector and
$x_{1},\ldots,x_{k}$ are non-constant maps. For each
$i\in \{1,\ldots,k\}$ we put $E_{i}=\hbox{Span}\{x_{i}(p):p\in
M\}.$ Then each $E_i$ is a linear subspace of ${\mathbb E}^m$.  

(1) If  the immersion $x$ is  linearly independent, then the
 subspaces $E_{1},\ldots,$ $E_{k}$ are linearly independent, that is,
the dimension of subspace spanned by all vectors in $E_{1}\cup
\ldots \cup E_{k}$ is equal to $\dim E_{1}+\ldots+\dim E_k$.
Since the Laplace transformation is homothetic, (1.16) yields
$$L=L_1+\ldots+L_k,\quad \Delta_L L_i=c\lambda_iL_i,\quad i=1,\ldots,k,
\leqno (1.17)$$
where $L_i=0$ when $\lambda_i=0$. Because $Span\{ L_i(p):p\in M\}=
Span \{x_i(p):p\in M\}=E_i$ for $L_i\not= 0$, (1.17) implies that
the Laplace map is also linearly independent. 

(2) If  the immersion $x$ is 
orthogonal, then the subspaces $E_{1},\ldots,E_{k}$ are
mutually orthogonal in ${\mathbb E}^m$.
Because $L$ has the spectral deocomposition given by (1.17),
statement (2) can be proved in a way similar to statement (1). \sq

{\bf Theorem 1.12.} {\it Let $x : M \rightarrow {\mathbb E}^m$ be
an equivariant isometric immersion of a compact
irreducible Riemannian homogeneous manifold. Then, with
respect to the induced metric, the Laplace map $L :
M \rightarrow {\mathbb E}^m$ is a homothetic immersion of finite type.
}

{\bf Proof.}  If $x : M \rightarrow {\mathbb E}^m$ be
an equivariant isometric immersion of a compact
irreducible Riemannian homogeneous manifold, then $x$ is of finite type
[C5, page 258]. On the other hand, by Theorem 1.8, we know that the 
Laplace transformation of $x$ is homothetic. Thus, by applying Proposition
1.10, the Laplace map $L$ is a homothetic immersion of finite type.
\sq

\vskip.2in
\noindent {\bf \S 2. Some classification theorems.}
\vskip.1in

In this section we assume $x : M^n \rightarrow {\mathbb E}^{m}$ 
is an isometric immersion whose  Laplace
transformation is given by ${\Cal L} : M^n \rightarrow L(M^{n})$.
If $\Cal L$ is homothetic, then the Laplacian operator $\Delta_L$ of
the Laplace image is related with the Laplace operator of $M$ by
$\Delta_L=c\Delta$ for some positive number $c$.

First we give  the following.

{\bf Theorem 2.1.} {\it Let $x : M^n \rightarrow {\mathbb E}^{m}$ 
be an isometric immersion with homothetic  Laplace
transformation  ${\Cal L} : M^n \rightarrow L(M^{n})$.
Then the Laplace image  lies in a hypersphere of
${\mathbb E}^m$ as  a minimal submanifold if and
only if either 

{\rm (1)} $M^n$ is
 a minimal submanifold of a hypersphere of ${\mathbb E}^m$ or 
 
{\rm (2)} $x$ is of null 2-type.}

{\bf Proof.}  Let $x : M^n \rightarrow {\mathbb E}^{m}$ 
be an isometric immersion with homothetic  Laplace
transformation. If the Laplace image lies in a hypersphere 
$S^{m-1}$ of ${\mathbb E}^m$ centered at the origin
as a minimal submanifold, then the Laplace map $L$ of $x$ is of
non--null 1--type by Theorem 2.4 of Chapter II. Thus
$$\Delta^2 x =\Delta L={1\over c}\Delta_LL={{\lambda}\over c}
 L={{\lambda}\over c} \Delta x,\leqno (2.1)$$
for some positive number $\lambda$. Put
$$x_1={c\over{\lambda}}x,\quad x_0=x-x_1.$$
Then, by (2.1), we have
$$\Delta x_1={c\over{\lambda}} \Delta^2 x=\Delta x={{\lambda}\over c}x_1,
\quad \Delta x_0=\Delta x-\Delta x_1=0.$$
Thus, $x$ is either of 1--type (when $x_0=0$) or of null 2--type (when
$x_0\not= 0$). If the first case occurs, $M^n$ is immersed by $x$ as a
minimal submanifold of a hypersphere of ${\mathbb E}^m$. \sq

{\bf Proposition 2.2.} {\it  Let $x : M^n \rightarrow {\mathbb E}^m$ be
an equivariant isometric immersion of a compact
irreducible Riemannian homogeneous manifold. Then
\begin{itemize}
\item[(1)] the immersion $x$ is of 1-type if and only
if the Laplace image $L (M^n )$ is a minimal submanifold
of a hypersphere; and
\item[(2)] if the immersion $x$ is not of 1-type,
then the Laplace image is a minimal submanifold in some
hyperquadric of ${\mathbb E}^m$.
\end{itemize}}

{\bf Proof.}  If $x : M^n \rightarrow {\mathbb E}^m$ is an equivariant isometric
immersion of a compact irreducible Riemannian homogeneous manifold,
then $x$ has homothetic Laplace transformation according to Theorem
1.8. Thus, by applying Theorem 2.1 we obtain statement (1), since
$M^n$ is compact.

In general, Lemma 1.4 implies that the Laplace map $L$ of $x$ is
an  equivariant homothetic
immersion. By multiplying the metric on $M$ by a suitable constant, $L$
becomes an  equivariant isometric immersion. Therefore, by applying Theorem 6.3
of Chapter II, we see that the Laplace image of $M$ is a minimal
submanifold of some hyperquadric of ${\mathbb E}^m$. \sq

The following result classifies all submanifolds with homothetic Laplace 
transformation and with parallel mean curvature vector.

{\bf Theorem 2.3.} {\it Let $x : M \rightarrow {\mathbb E}^{m}$  be
an isometric immersion with parallel mean curvature
vector.  Then the Laplace transformation  $ {\cal L}: M
\rightarrow L(M)$ of $x$ is homothetic if and only if 
$M$ is immersed as a minimal submanifold of
a hypersphere of ${\mathbb E}^m$ via $x$.
}

{\bf Proof.}  Assume the immersion  $x : M \rightarrow {\mathbb E}^{m}$ has
parallel mean curvature and homothetic Laplace transformation.
Then by Lemma 1.1 we have
$$\<A_HX,A_HY\>=c^2\<X,Y\>,\leqno (2.2)$$
for vectors $X,Y$ tangent to $M$, where $c$ is a  positive number.
 Put $$U=\{p\in M:H(p)\not=0\}.$$ Then
$U$ is a dense open subset of $M$, since $\cal L$ is homothetic.
On $U$ we choose a local orthonormal frame field $e_1,\ldots,e_n$
of the tangent bundle and a local orthormal  frame field
 $e_{n+1},\ldots, e_m$ of the normal
bundle such that
$e_1,\ldots,e_n$ are eigenvectors of $A_H$
and $e_{n+1}$ is in the 
direction of the mean curvature vector $H$ on $U$. Hence we have
$$A_{n+1}e_i=\kappa_ie_i,\quad i=1,\ldots,n,\leqno(2.3)$$
where $A_{r}=A_{e_r}$ is the Weingarten map of $M$ and $\kappa_1,\ldots,
\kappa_n$ are eigenvalues of $A_{n+1}$. From (2.2) and (2.3) we may
assume
$$\kappa_1=\ldots=\kappa_k=\lambda,\quad \kappa_{k+1}=\ldots=\kappa_n
=-\lambda,\leqno (2.4)$$
where $\lambda =n\alpha/(2k-n)$ and $\alpha^2=\<H,H\>$ are constant.

Assume $0<k<n$.
We put $${\Cal D}_1=\{X\in TM:A_{n+1}X=\lambda X\},$$
$$
{\cal D}_2=\{X\in TM:A_{n+1}X=-\lambda X\}.\leqno(2.5)$$
Because $DH=0$,  the equation of Codazzi implies that ${\cal D}_1,
{\cal D}_2$ are integrable distributions on $U$. Furthermore, by applying
Codazzi's equation, it is easy to show that leaves of ${\cal D}_1$ and
${\cal D}_2$ are totally geodesic submanifolds of $U$. Therefore, locally
$U$ is the Riemannian product $M_1\times M_2$ of some integrable submanifolds
of ${\cal D}_1$ and ${\cal D}_2$.  

From (2.3) and (2.4) we know that $A_{n+1}$ takes the following form:
$$A_{n+1}=\begin{pmatrix} \lambda I_k & 0\\
0 & -\lambda I_{n-k}\end{pmatrix}.\leqno (2.6)$$
On the other hand, since the mean curvature vector $H$ is parallel, the equation of Ricci yields
$[A_{n+1},A_r]=0.$ Therefore, by using (2.6), the second fundamental form
$h$ of $U$ in ${\mathbb E}^m$ satisfies
$h({\cal D}_1,{\cal D}_2)=\{0\}.$
Consequently, by applying a lemma of Moore, we know that $M$ is locally a
product submanifolds, say 
$$x=(y,z): M_1\times M_2\rightarrow {\mathbb E}^{m_1}\times
{\mathbb E}^{m-m_1}={\mathbb E}^m.\leqno (2.7)$$  Let $\bar e_1,\ldots,\bar e_n,\bar
e_{n+1},\ldots,\bar e_m$ be a  local orthonormal frame field such
that $\bar e_{n+1}$ is in the direction of the mean curvature
 vector $H_1$ of
$M_1$ in ${\mathbb E}^{m_1}$, $\bar e_{n+2}$ is in the direction of the mean
curvature  vector $H_2$ of $M_2$ in ${\mathbb E}^{m-m_1}$, $\bar e_1,\ldots,\bar
e_k$ are tangent to $M_1$ and $\bar e_{k+1},\ldots,\bar e_n$ are
tangent to $M_2$, respectively. Then we have
$$H={1\over n}\{(\hbox{trace} A_{\bar e_{n+1}})\bar e_{n+1}
+(\hbox{trace} A_{\bar e_{n+2}})\bar e_{n+2}\}.\leqno (2.8)$$
From (2.8) we get
$$A_H={1\over n}\{(\hbox{trace} A_{\bar e_{n+1}})A_{\bar e_{n+1}}
+(\hbox{trace} A_{\bar e_{n+2}})A_{\bar e_{n+2}}\}.\leqno (2.9)$$
In particular, if  $\bar e_1,\ldots,\bar e_k$ are tangent to $M_1$
and  eigenvectors
of $A_{\bar e_{n+1}}$ and $\bar e_{k+1},\ldots,$ $\bar e_{n}$ are 
tangent to $M_2$ and are eigenvectors of $A_{\bar e_{n+2}}$, then we
have $$A_H={1\over n}\begin{pmatrix} (\hbox{trace} A_{\bar e_{n+1}}) B &
0\\ 0 & (\hbox{trace} A_{\bar e_{n+2}})C\end{pmatrix},\leqno (2.10)$$
where $B$ and $C$ are diagonal matrices given by
$$B=\hbox{Diag}(\delta_1,\ldots,\delta_k),\quad
C=\hbox{Diag}(\beta_{k+1},\ldots,\beta_n).\leqno (2.11)$$
Since $H=\alpha e_{n+1}$, (2.6), (2.10) and (2.11) imply
$$\delta_1=\ldots=\delta_k=\delta,\quad \beta_{k+1}=\ldots=\beta_n=\beta.$$
Therefore, we find $n\alpha\lambda=k\delta^2$ and $-n\alpha\lambda
=(n-k)\beta^2$ which is impossible, since $\alpha$ and $\lambda$ are non--zero
on $U$. Consequently, either $k=0$ or $k=n$. 
Thus, $U$ is a pseudo--umbilical submanifolds with parallel mean
curvature vector. Hence, $U=M$.  Therefore, by applying a result of
[CY] (i.e. Theorem 1.1 of Chapter II), we conclude
 that $M$ is immersed by $x$ as
a minimal submanifold of a hypersphere of ${\mathbb E}^m$.

The converse of this is given by Lemma 1.2. \sq

In the following we study hypersurfaces with homothetic Laplace 
transformation. For simpicity, here 
we classify hypersurfaces of dimenson
2 and 3 only. 

{\bf Theorem 2.4.} {\it Let $x : M \rightarrow {\mathbb E}^3$ be a
surface in ${\mathbb E}^3$. Then the Laplace transformation
${\Cal L} : M \rightarrow L(M)$ of the surface is homothetic if and only if $M$ is an open portion
of an ordinary sphere in ${\mathbb E}^3$.} 

{\bf Proof.}  If $M$ is an open protion of an ordinary sphere in
${\mathbb E}^3$,  it is easy to see that it has homothetic Laplace
transformation. 

Conversely, suppose $M$ is a surface of ${\mathbb E}^3$ with homothetic Laplace
transformation. It suffices to show that $M$ has constant mean
curvature  according to Theorem 2.3.

 Let $U=\{p\in M:d\alpha^2\not=0\; \hbox{at}\; p\,\}.$ Then
$U$ is an open subset of $M$. 
In order to prove that $U$ is an empty set, we
choose a local orthonormal frame field
 $e_1,e_2,e_3$   such that $e_1,e_2$ are tangent vector
fields of $U$ which are eigenvectors of
the Weingarten
map $A_{3}=A_{e_{3}}$. So, we have
$$A_{3}e_i=\kappa_i e_i,\quad i=1,2,\leqno(2.12)$$
where $\kappa_1,\kappa_2$ are principal curvatures.

From Lemma 1.1 and (2.12) we get
$$(e_1\alpha)(e_2\alpha)=0.\leqno(2.13)$$
Therefore, without loss of generality, we may choose $e_1$  such that
$e_1$ is in the direction of  the gradient of $\alpha,\;\nabla\alpha$.
So we have
$$e_2\alpha=0.\leqno (2.14)$$
 From Lemma 1.1 and (2.12) we also have
$$\alpha^2\kappa^2_i+(e_i\alpha)^2=c^2,\leqno(2.15)$$
where $c$ is a positive number.
 Thus  (2.14) yields
$$\alpha^2\kappa^2_2=c^2.\leqno(2.16)$$
So, by choosing suitable $e_3$, we have
$$\kappa_2={c\over{\alpha}}>0.\leqno(2.17)$$ 

From (2.15) and (2.17), we have
$$\omega^3_1=(2\alpha-{c\over{\alpha}})\omega^1,\quad
\omega^3_2={c\over{\alpha}}\omega^2,
\leqno(2.18)$$
$$(e_1\alpha)^2=4\alpha^2(c-\alpha^2).\leqno(2.19)$$
Put
$$\omega^2_1=f_1\omega^1+f_2\omega^2.\leqno(2.20)$$
Then we have
$$d\omega^1=f_1\omega^1\wedge\omega^2,\quad d\omega^2=f_2
\omega^1\wedge\omega^2.\leqno(2.21)$$

By taking exterior derivative of (2.18)  and using (2.21) and the
structure equations, we obtain
$$\alpha^2=cf_1,\quad e_1\alpha=2(\alpha-{{\alpha^3}\over c})f_2.
\leqno(2.22)$$
Combining (2.19) and (2.22) we find
$$(c-\alpha^2)f_2^2=c^2.\leqno(2.23)$$
(2.22) and (2.23) yield
$$e_2f_1=e_2f_2=0.\leqno(2.24)$$
Taking the exterior derivative of (2.20), we find
$$f_1^2+f_2^2={{c(c-2\alpha^2)}\over{\alpha^2}}-e_1f_2.\leqno(2.25)$$
On the other hand, (2.22) and (2.23) imply
$$e_1f_2={{c^5-4c^4\alpha^2+2c^3\alpha^4-c\alpha^6+\alpha^8}\over
{c^2\alpha^2(c-\alpha^2)}}.\leqno(2.26)$$
From (2.23) we get
$$f_2(e_1f_2)={{c^2\alpha(e_1\alpha)}\over{(c-\alpha^2)^2}}.\leqno(2.27)$$
Using (2.19), (2.23) and (2.27) we find
$$(e_1f_2)^2={{4c^2\alpha^4}\over{(c-\alpha^2)^2}}.\leqno (2.28)$$
Combining (2.26) and (2.28), we conclude that  $U$ is an empty
set. Thus, $M$ has constant mean curvature and it
 is an open portion of an ordinary sphere in ${\mathbb E}^3$. \sq

{\bf Theorem 2.5.} {\it Let $x : M \rightarrow {\mathbb E}^{4}$ be a
 hypersurface of ${\mathbb E}^4$. Then the Laplace transformation
${\Cal L} : M \rightarrow L(M)$ of the hypersurface is homothetic if
and only if $M$ is an open portion
of a hypersphere of ${\mathbb E}^4$.} 

{\bf Proof.} 
Suppose $M$ is a hypersurface of ${\mathbb E}^4$ with homothetic Laplace
transformation. It suffices to show that $M$ has constant mean
curvature  according to Theorem 2.3.

 Let $U=\{p\in M:d\alpha^2\not=0 \;\;\hbox{at}\; p\,\}.$ Then
$U$ is an open subset of $M$. 
In order to prove that $U$ is an empty set, we
choose a local orthonormal frame field
 $e_1,e_2,e_3,e_4$   such that $e_1,e_2,e_3$ are tangent vector
fields of $U$ which are eigenvectors of
the Weingarten
map $A_4=A_{e_4}$. So, we have
$$A_{n+1}e_i=\kappa_i e_i,\quad i=1,2,3,\leqno(2.29)$$
where $\kappa_1,\kappa_2,\kappa_3$ are principal curvatures.

From Lemma 1.1 and (2.29) we get
$$(e_i\alpha)(e_j\alpha)=0,\quad i\not=j.\leqno(2.30)$$
Therefore, without loss of generality, we may choose $e_1$  such
 that $e_1$ is in the direction of the
gradient of $\alpha,\;\nabla\alpha$.
So we have
$$e_2\alpha=e_3\alpha=0.\leqno (2.31)$$
 From Lemma 1.1
and (2.29) we have
$$\alpha^2\kappa^2_i+(e_i\alpha)^2=c^2,\leqno(2.32)$$
where $c$ is a positive number.
 Thus (2.31) yields
$$\alpha^2\kappa^2_2=\alpha^2\kappa_3^2=c^2.\leqno(2.33)$$
So, by choosing a suitable $e_4$, we have 
either
$$\kappa_1=3\alpha,\quad \kappa_2={c\over {\alpha}},\quad
\kappa_3=-{c\over{\alpha}},\leqno(2.34)$$
or
$$\kappa_1=3\alpha-{{2c}\over{\alpha}},\quad \kappa_2=
\kappa_3={c\over{\alpha}},\leqno(2.35)$$

We treat these two cases separately.

{\bf Case 1.}  (2.34) holds. In this case, we put $\mu={c\over {\alpha}}$.
From (2.32) and (2.34) we have
$$(e_1\alpha)^2=c^2-9\alpha^4.\leqno (2.36)$$
$$\omega^4_1=3\alpha \omega^1, \quad \omega^4_2=\mu\omega^2, \quad
\omega^4_3=-\mu\omega^3.\leqno(2.37)$$
By taking the  exterior derivative of (2.37) and applying (2.37)
and the structure equations, we obtain
$$\omega_1^2(e_1)=\omega^3_1(e_1)=0.\leqno(2.38)$$
on $U$. Hence, integral curves of $e_1$ in $U$ are geodesics 
of $M$. 

Since $\alpha\mu=c$, (2.31) yields
$$e_2\mu=e_2\mu=0.\leqno(2.39)$$

By taking exterior derivatives of the last two equations of (2.37) and
using (2.37), (2.39) and the structure equations, we obtain
$$e_1\mu=(3\alpha-\mu)\omega^2_1(e_2)=-(3\alpha+\mu)
\omega^3_1(e_3),\leqno(2.40)$$
$$\omega^2_1(e_1)=\omega^3_1(e_1)=\omega^3_2(e_2)=\omega^3_2(e_3)=0,
\leqno (2.41_{}$$
$$(3\alpha-\mu)\omega_1^2(e_3)=2\mu\omega_2^3(e_1),
\quad (3\alpha+\mu)\omega^3_1(e_2)=2\mu\omega^3_2(e_2)
.\leqno(2.42)$$

From (2.41) we have
$$\omega^2_1=f_2\omega^2+f_3\omega^3,\leqno(2.43)$$
$$ \omega^3_1=g_2
\omega^2+g_3\omega^3,\leqno(2.44)$$
$$ \omega^3_2=h\omega^1,\leqno (2.45)$$
for some functions $f_2,f_3,g_2,g_3,h$.
Moreover, by (2.40), and (2.42), we have
$$e_1\alpha={\alpha\over c}(c-3\alpha^2)f_2=
{\alpha\over c}(c+3\alpha^2)g_3,\leqno (2.46)$$
$$(3\alpha^2+c)g_2=2c h,\quad (3\alpha^2-c)f_3=2c h.\leqno(2.47)$$
Also, from (2.46) and (2.47), we have
$$g_3=({{c-3\alpha^2}\over{c+3\alpha^2}}) f_2,
\quad
g_2=-({{c-3\alpha^2}\over{c+3\alpha^2}}) f_3.\leqno(2.48)$$
By taking exterior derivative of (2.34) we obtain
$$e_2h=e_3h=0,\quad dh=(e_1h)\omega^1,\leqno (2.49)$$
$$\alpha^2h(g_2-f_3)=c^2+\alpha^2(g_2f_3-g_3f_2).\leqno (2.50)$$
Similarly, by taking exterior derivative of (2.43) and (2.44), we find
$$e_1f_2=hg_2-3c-f^2_2+f_3(h-g_2),\leqno (2.51)$$
$$e_1f_3=hg_3-f_3g_3-(f_3+h)f_2,\leqno(2.52)$$
$$e_3f_2=e_2f_3,\leqno2.53$$
$$e_1g_2=-hf_2-f_2g_2+(h-g_2)g_3,\leqno(2.54)$$
$$e_1g_3=-hf_3+3c-(f_3+h)g_2-g_3^2,\leqno(2.55)$$
$$e_3g_2=e_2g_3.\leqno(2.56)$$
From (2.31) and (2.46) we have
$$e_2e_1\alpha={\alpha\over c}(c+3\alpha^2)(e_2g_3),\quad
e_1e_2\alpha=0.$$
Therefore, using (2.41), we get
$${\alpha\over c}(c+3\alpha^2)(e_2g_3)=[e_2,e_1]\alpha=-
(\nabla_{e_1}e_2)\alpha=0.$$
Hence, by (2.56), we get
$$e_2g_3=e_3g_2=0.\leqno (2.57)$$
Similarly by using
 $$e_3e_1\alpha={\alpha\over c}(c+3\alpha^2)(e_3g_3),\quad
e_1e_3\alpha=0,$$
we obtain 
$$e_3g_3=0.\leqno(2.58)$$ 

Similarly using $e_1\alpha={\alpha\over c}(c-3\alpha^2)f_2$,
we also obtain
$$e_2f_2=e_3f_2=0,\quad e_2f_3=0.\leqno (2.59)$$
From (2.47) and (2.57), we also find
$$e_2h=e_3h=0.\leqno (2.60)$$
Consequently, from (2.31), (2.47) and (2.60), we obtain
$$e_ig_j=e_if_j=e_ih=0,\quad i,j=2,3.\leqno (2.61)$$
By using (2.37) and a long computation we derive
$$e_1g_2={{f_2f_3(c-3\alpha^2)}\over{c+3\alpha^2}}.\leqno (2.62)$$

On the other hand, substituting (2.47) and (2.48) into
(2.54)  yields
$$e_1g_2={{f_2f_3(c-3\alpha^2)}\over{c(c+3\alpha^2)^2}}
(2c^2+3c\alpha^2+9\alpha^4).\leqno (2.63)$$
Combining (2.51) and (2.52), we obtain $f_2f_3=0$.

 If $f_2=0$, then $g_3=0$ according to (2.48). Thus, 
(2.46) and  $e_1\alpha=0$ imply $M$ has constant mean curvature.

If $f_3=0$, (2.47) and (2.48) yield $h=g_2=0$. Thus, (2.51)
implies 
$$e_1f_2=-3c-f^2_2.\leqno(2.64)$$
On the other hand, by taking derivative of $\omega^2_1=f_2\omega^2$ and
using (2.32),  (2.43), (2.45) together with $h=f_3=0$, we have
$$e_1f_2={{c^2}\over{\alpha^2}}-f^2_2.\leqno (2.65)$$
Combining (2.64) and (2.65) we see that $U$ is an empty set which
implies that $M$ has constant mean curvature $M$ is an 
open portion of a hypersphere. 

{\bf Case 2.}  (2.35) holds. In this case, we have
$$\omega^4_1=(3\alpha-{{2c}\over{\alpha}})\omega^1,\quad
\omega^4_2={c\over\alpha}\omega^2,\quad \omega^4_3={c\over\alpha}
\omega^3,\leqno(2.66)$$
$$(e_1\alpha)^2=3(c-\alpha^2)(c-3\alpha^2).\leqno(2.67)$$

By taking exterior derivatives of the three equations in (2.66), we
may prove after a long computation that 
$$\omega^2_1=\omega^3_1=0.\leqno(2.68)$$ 
Thus, $d\omega^1=0$. Hence, if we put ${\cal D}_1=\hbox{Span}
\{e_1\}$ and ${\cal D}_2=\hbox{Span}\{e_2,e_3\}$ on $U$, then
${\cal D}_1,{\cal D}_2$ are totally geodesic integrable distributions
on $U$. Hence, $U$ is locally the Riemannian product of a curve and
a surface $N$. Moreover, since the second fundamental form $h$ satisfies
$h({\cal D}_1,{\cal D}_2)=\{0\}$, a lemma of Moore implies that $U$ 
is locally a product hypersurface of a line and a surface in a affine
3--subspace ${\mathbb E}^3$ of ${\mathbb E}^4$. This is a contradiction since such a 
hypersurface does not have a homothetic Laplace transformation. 
Consequently, $U$ is an empty set. Thus, $M$ is an open portion of
a hypersphere of ${\mathbb E}^4$.

The converse of this is trivial. \sq

{\bf Remark 2.1.} Similar but much more complicated arguments 
can be used to classify
hypersurfaces of a Euclidean space of any dimension with homothetic
Laplace  transformation. \sq

\vskip.2in
\noindent{\bf \S3. Surfaces with homothetic Laplace tranformation.}
\vskip.1in

The purpose of this section is to study surfaces with
homothetic Laplace transformation.

In the following, a submanifold $M$ of 
a Riemannian manifold is said to
have {\it parallel normalized mean curvature vector\/} if the
mean curvature vector $H$ is nonzero and the unit vector
field in the direction of $H$ is parallel in the normal
bundle. It is clear that a non-minimal hypersurface 
has parallel normalized mean curvature vector.
For submanifolds with parallel normalized mean curvature
vector, we have the following.

{\bf Theorem 3.1.} {\it Let $x: M \rightarrow {\mathbb E}^{m}$ be an 
isometric immersion with
parallel normalized mean curvature vector of a surface $M$. Then the
Laplace transformation 
 ${\Cal L} : M \rightarrow L(M)$ of
$x$ is homothetic if and only if $M$ is immersed as a
minimal surface in a hypersphere of ${\mathbb E}^m$ via $x$.
}

{\bf Proof.} Assume $M$ is a surface in $E^m$ with parallel 
normalized mean curvture vector and homothetic Laplace 
transformation. Denote by $\alpha$ the mean curvature of $M$
in $E^m$. We choose an orthonormal local frame field
$e_1,e_2,e_3,$ $\ldots,e_m$ such that $e_1,e_2$ are eigenvectors
of $A_H$ and $H=\alpha e_3$. We put $A_3e_i=\kappa_i e_i,i=1,2$.
Then we have $(e_1\alpha)(e_2\alpha)=0$. So, we may choose
$e_1,e_2$ such that $e_1$ is in the direction of the gradient
of $\alpha^2$. Let $U$ be the open subset of $M$ on which
the gradient of $\alpha^2$ is non--zero. Thus, we have
$$e_2\alpha=0.\leqno(3.1)$$
Hence,  from Lemma 1.1 and (3.1), we have
$$\kappa_1=2\alpha-{c\over\alpha},\quad \kappa_2={c\over\alpha},
\quad \omega^3_1=\kappa_1\omega^1,\quad \omega^3_2=\kappa_2\omega^2,\leqno(3.2)$$
$$(e_1\alpha)^2=4\alpha^2(c-\alpha^2),\leqno (3.3)$$
where $c$ is a positive constant. (3.3) implies
$c>\alpha^2$ on $U$. Becaus $M$ has parallel
normalized mean curvature vector, we have $De_3=0$; and hence
$$\omega^3_4=\cdots=\omega^3_m=0.\leqno(3.4)$$ 
From the equation of Ricci and (3.4), we get
$$[A_3,A_r]=0,\quad r=4,\ldots,m.\leqno(3.5)$$
Since the gradient of $\alpha^2$ is nowhere zero on $U$,
$c\not=\alpha^2$. Therefore, (3.5) implies that each  $A_r$
is diagonalized with respect to $e_1,e_2$. Because $e_3$ is 
parallel to $H$, this implies that one may choose $e_4,\ldots,
e_m$ in such a way that $A_4,\ldots,A_m$ take the following
forms:
$$A_4=\begin{pmatrix} \beta & 0\\0 & -\beta\end{pmatrix},\quad
A_5=\cdots=A_m=0.\leqno(3.6)$$
By taking exterior derivatives of $\omega^3_1$ and $\omega^3_2$
and using the structure equations and (3.2), we find
$$(2\alpha-{{2c}\over\alpha})f_1=0,\leqno(3.7)$$
$$e_1\alpha={{2\alpha}\over c}(c-\alpha^2)f_2,\leqno(3.8)$$
where $\omega^2_1=f_1\omega^1+f_2\omega^2.$ Since $\alpha$ is 
nowhere constant on $U$, (3.7) yields $f_1=0$. Put $f=f_2$. 
Then $\omega^2_1=f\omega^2.$

Taking exterior derivative of $\omega^4_1=\beta\omega^1$ and
$\omega^4_2=-\beta\omega^2$ we find
$$e_1\beta=-2\beta f,\quad e_2\beta=0.\leqno(3.9)$$

Combining (3.3) and (3.8), we get
$$\omega^2_1=f\omega^2=\pm{{c}\over{\sqrt{c-\alpha^2}}}
\omega^2.\leqno(3.10)$$
By taking exterior derivative of (3.10) and applying (3.8) and
(3.10), we obtain
$$\beta^2=2c-{{c^2}\over{\alpha^2}}+{{2c\alpha^2+c^2}\over
{c-\alpha^2}}.\leqno(3.11)$$

Taking exterior derivative of (3.11) and using (3.8) we get
$$\beta(e_1\beta)=\pm2c\alpha\sqrt{c-\alpha^2}
\{{{2c}\over{\alpha^3}}+{{6c\alpha}\over
{(c-\alpha^2)^2}}\}.\leqno(3.12)$$
From (3.9), (3.11) and (3.12), we may conclude that $U$ is an
empty set. Thus, $M$ has constant mean curvature. Because $M$
has parallel normalized mean curvature vector, $M$ has paralllel
mean curvature vector. Therefore, by Theorem 2.3, $M$ is 
immersed as a minimal surface in a hypersphere of ${\mathbb E}^m$. The
converse follows from Lemma 1.2. \sq

{\bf Theorem 3.2.}  {\it Let $x : M \rightarrow {\mathbb E}^{4}$ be
a surface in ${\mathbb E}^4$ with constant mean curvature. Then $M$
has homothetic Laplace transformation  if and
only if $M$ is  a minimal surface in a
hypersphere of ${\mathbb E}^4$.}

{\bf Proof.}  Assume $M$ has constant mean curvatur and homothetic
Laplace transformation. Then the mean curvature is nonzero.
If $M$ has parallel mean curvature vector, then $M$ is a minimal
surface in a hypersphere of ${\mathbb E}^4$ by Theorem 3.1. So, we assume
that $U=\{p\in M:DH\not= 0\;\hbox{at}\; p\,\}$ is non--empty. It is
clear that $U$ is an open subset of $M$.

On $U$, we choose a local orthonormal frame field $e_1,e_2,e_3,e_4$
such that $e_3$ is in the direction of $H$ and $e_1,e_2$ are 
eigenvectors of $A_3$ with eigenvalues $\kappa_1,\kappa_2$,
respectively. Since $DH=\alpha De_3=\alpha \omega_3^4$,
Lemma 1.1 implies
$\omega_3^4(e_1)\omega^4_3(e_2)=0.$ Without loss of generality,
we may assume $\omega^4_3(e_1)=0$. Then we have
$$D_{e_1}H=0,\quad \omega^4_3=\mu\omega^2,\leqno(3.14)$$ 
for some local function $\mu$. So, by Lemma 1.1, we have
$$\alpha^2\kappa_1^2=c^2,\quad \alpha^2\kappa_2^2 +
\alpha^2\mu^2=c^2,\leqno (3.15)$$
where $c$ is a positive constant. Because $M$ has constant
mean curvature, (3.15) implies $\kappa_1, \kappa_2$ and
$\mu$ are constant. Put
$$A_4=\begin{pmatrix} \gamma & \delta\\ \delta &-\gamma
\end{pmatrix}.\leqno(3.16)$$
Taking exterior derivative of $\omega^3_1=\kappa_1\omega^1$
and $\omega^3_2=\kappa_2\omega^2$, respectively, we obtain
$$(\kappa_2-\kappa_1)f_1=-\mu\gamma,\leqno(3.17)$$
$$e_1\delta+2\gamma f_1=\kappa_1\mu+e_2\gamma-2f_2\delta,
\leqno(3.18)$$
$$(\kappa_1-\kappa_2)f_2=\delta\mu,\leqno(3.19)$$
$$e_2\delta+2f_2\gamma=-e_1\gamma+2f_1\delta,\leqno(3.20)$$
where $\omega^2_1=f_1\omega^1+f_2\omega^2$.
From (3.17) and (3.19) we find
$$f_1\delta+\gamma f_2=0.\leqno(3.21)$$
Taking exterior derivative of $\omega^4_3=\mu\omega^2$, we find
$$\mu f_2=\delta(\kappa_2-\kappa_1).\leqno(3.22)$$
From (3.19) and (3.22) we get
$$\{\mu^2+(\kappa_1-\kappa_2)^2\}f_2=0.
\leqno(3.23)$$

If $\kappa_1\equiv\kappa_2$ on $U$, then $U$ is a pseudo--umbilical
surface with non--zero constant  mean curvature in $E^4$. Thus,
by applying a result of [C1] (cf. Theorem 1.2 of Chapter II), we know
that  $U$ has parallel mean curvature vector. This is a contradiction. 
Therefore,  $V=\{ p\in U:\kappa_1\not=\kappa_2
\;\hbox{at}\; p\,\}$ is a non--empty open subset of $U$.  On $V$,
(3.19) and (3.23) yield $f_2=0$ and $\delta\mu=0$. Since $DH\not=0$
on $U$, we obtain $\delta=0$.

Since $f_2=0$, we have $\omega^2_1=f_1\omega^1$. So, by taking
exterior derivative of $\omega^2_1$ and applying $\delta=0$, 
we find
$$e_2f_1=f_1^2-\gamma^2+\kappa_1\kappa_2.\leqno(3.24)$$
On the other hand, (3.17) implies
$$e_2f_1=({\mu\over{\kappa_1-\kappa_2}})(2\gamma f_1-\kappa_1\mu).
\leqno(3.25)$$
Combining (3.17), (3.24) and (3.25) we conclude that $\gamma$ is constant. Thus, 
by (3.17) and (3.22), we conclude that $f_1$ is also a constant and
$$f_1^2=\kappa_1\kappa_2-\gamma^2.\leqno (3.26)$$
Moreover, using (3.17) and (3.18), we find
$$f_1^2={{\kappa_1\mu^2}\over{2(\kappa_1-\kappa_2)}}\leqno(3.27)$$
Because $\kappa_1\not=0$, (3.26) and (3.27) give
$$\mu^2=(\kappa_1-\kappa_2)(3\kappa_2-\kappa_1).\leqno(3.28)$$
From (3.15) and (3.28) we obtain $\kappa_1=\kappa_2$ on $V$
which is a contradiction. Consequently, $M$ has parallel mean curvature
vector. Hence, $M$ is immersed as a minimal surface of a 
hypersphere of ${\mathbb E}^4$. 

The converse of this is clear. \sq

It is interesting to point out that Theorem 3.1 and Theorem 3.2
 are best possible. In fact we have  the following
results.

{\bf Proposition 3.3.} {\it Let $C_1$ and $C_2$ be two
planar curves  parametrized by
arclength. Then the product surface $M=C_1\times C_2$
in ${\mathbb E}^4$ has homothetic Laplace transformation if and only if
the curvature functions of $C_1, C_2$ are of the following form:
$$\kappa(s) =c(1+c^4 a e^{-8c^2s})^{-1\over4}\leqno(3.29)$$ for some
 constants $a$ and $c>0$. 
}

{\bf Proof.}  Assume $C_1\times C_2$ is given  by
$$X(u,v)=(x(u),y(u),z(v),w(v)),\leqno(3.30)$$
where $u,v$ are arclength parametrizations of $C_1$ and $C_2$,
respectively. Then the mean curvature vector $H$ of $M$ in $E^4$ 
is given by
$$H={1\over 2}(x''(u),y''(u),z''(v),w''(v)).\leqno(3.31)$$

Denote  $e_1={\partial\over{\partial u}},
e_2={\partial\over{\partial v}}$. Then we have
$$\aligned L_*(e_1)=(-x'''(u),-y'''(u),0,0),\\
L_*(e_2)=(0,0,-z'''(v),-w'''(v)).\endaligned \leqno(3.32)$$
Therefore, by Lemma 1.1, the plane curvatures $\kappa$
of $C_1,C_2$ satisfy the following differential equation:
$$(\kappa'(s))^2+\kappa(s)^4=c^2,\leqno(3.33)$$
for some positive constant $c$, where $s$ is 
the arc--length parametrization of the curve. (3.32) implies
$c\geq\kappa^2.$ Solving equation (3.33), we obtain
$$\kappa(s)^4={{c^2}\over{1+c^2ae^{-8c^2s}}},\leqno(3.34)$$
where $a$ is a constant.

 The converse of this is easy to verify.
\sq
 
{\bf Remark 3.1.} Proposition 3.3 still holds if we
replace $C_1\times C_2$ by the product of $k$ plane curves 
in ${\mathbb E}^{2k}$, for any $k=1,3,4,\cdots,$ (see also Corollary 1.3 of
Chapter III). \sq

{\bf Remark 3.2.} The plane curve whose curvature function
satisfies (3.29) 
is an open portion of a circle if and only if $a=0$. When $a\not=0$,
then a plane curve  whose curvature function satisfies
condition (3.29) is congruent to a curve given by
$$\Bigg(\int \cos\theta(s)ds,\int \sin\theta(s)ds\Bigg),\leqno (3.35)$$
where$$\theta(s)=c_1-{1\over {8c^{3\over2}}}\{\ln (\ell(s)-1)
-\ln(\ell(s)+1)+2
\tan^{-1}\ell(s)\},\quad c_1\in {\Bbb R},\leqno(3.36)$$
$$\ell(s)=(1+c^2ae^{-8c^2s})^{1\over 4} . \leqno (3.37)$$

{\bf Remark 3.3.} From Proposition 3.3 and Remark 3.2 it follows
that there exist amples examples of surfaces which do not
lie in any hypersphere of  ${\mathbb E}^4$ but with
homothetic Laplace transformation. \sq

The following two Propositions also show that Theorems 3.1 and 
3.3 are best possible.

{\bf Proposition 3.4.} {\it  Let $M^2$ be a surface in
${\mathbb E}^5$ defined by
$$x(u,v)=\Bigg(au, b^2 \cos u, b^2 \sin u, {1\over b}(a^2
+b^4)^{3/4}\cos v, {1\over b}(a^2
+b^4)^{3/4}\sin v\Bigg)$$
for some nonzero constants $a$ and $b$. Then the surface
$M^2$ in ${\mathbb E}^5$ satisfies the following properties:
\begin{itemize}

\item[(1)] the Laplace transformation is homothetic;

\item[(2)] the mean curvature function is a nonzero constant;

\item[(3)] $M^2$ has
non-parallel normalized mean curvature vector field; and

\item[(4)] $M^2$ is not contained in any hypersphere of ${\mathbb E}^5$.
\end{itemize}}

{\bf Proposition 3.5.} {\it Let $M^2$ be the product of
two circular helices  in ${\mathbb E}^6$ defined by
$$x(u,v)=(au, av, c\cos u, c\sin u, c\cos v, c\sin v)$$
for some nonzero constants $a$ and $c$. Then the surface
satisfies the following properties:
\begin{itemize}
\item[(1)] the Laplace transformation is homothetic;

\item[(2)] the mean curvature function is a nonzero constant;

\item[(3)] the surface is pseudo-umbilical;

\item[(4)] $M^2$ has non-parallel normalized mean curvature vector field;  

\item[(5)] $M^2$ is not contained in any hypersphere of ${\mathbb E}^6$; and

\item[(6)] the immersion $x$ is of null 2-type. 
\end{itemize} }

Propositions 3.4 and 3.5 can be proved by direct computation.

\vfill\eject

\noindent{\bf Chapter V: CONFORMAL LAPLACE TRANSFORMATION}

\vskip.2in
\noindent{\bf   \S1. Some general results. } \vskip.1in

Let $x : M \rightarrow \E^m$ be an isometric immersion
of an $n$-dimensional  connected Riemannian manifold
$M$ into a Euclidean $m$-space. Denote by  $L :M^{n} \rightarrow E^m$
 the {\it Laplace map\/} and by $L(M^n )$  the
{\it Laplace image\/} of the immersion $x$. 
As before, we denote the Laplace trasformation by
 ${\Cal L} :M \rightarrow L(M )$. The Laplace transformation
 $\Cal L$ is said to be {\it conformal\/} (respectively, 
{\it weakly conformal\/}) if there exists a function $\rho>0$
 (respectively, $\rho\geq 0$) such that
 $\<d{\Cal L}(X),d{\Cal L}(Y)\>=\rho^2 \<X,Y\>$ for all vectors
$X,Y$ tangent to $M$.

The main purpose of this section to give some general results
concerning submanifolds with conformal Laplace transformation.

First we give the following general lemma.

{\bf Lemma 1.1.} {\it Let $x : M \rightarrow \E^m$ be an isometric
immersion of an n--dimensional Riemannian manifold into $\E^m$.
Then the Laplace transformation ${\Cal L}:M\rightarrow L(M)$ is
conformal (respectively, weakly conformal) if and only if
$$\<A_HX,A_HY\>+\<D_XH,D_YH\>=\rho^2\<X,Y\>,\leqno (1.1)$$
holds for vectors X,Y tangent to $M$, where $\rho$ is a
strictly positive  function (respectively, positive function) on
$M$. }

\demo Because the Laplace map is given by $L(p)=(\Delta x)(p)=-n
H(p)$, the differential of the Laplace map satisfies
$$dL(X)=nA_HX-nD_XH.\leqno (1.2)$$
From (1.2) we obtain
$$\<dL(X),dL(Y)\>=\<A_HX,A_HY\>+\<D_XH,D_HY\>\leqno (1.3)$$
which implies the Lemma. \sq

{\bf Theorem 1.2.} {\it  Let $M$ be a submanifold in
$\E^m$. Then $M$ is a minimal submanifold of a hypersphere
of $\E^m$ if and only if  $M$ has confomal
Laplace transformation and the  mean curvature
vector field of $M$ in $\E^m$ is parallel in the normal
bundle.}

\demo If $M$ has parallel mean curvature vector and 
conformal Laplace transformation, then Lemma 1.1 yields
$$\<A_HX,A_HY\>=\rho^2\<X,Y\>,\leqno(1.4)$$
for some strictly positive function $\rho$. It is clear that $H$
is nowhere zero.

Let $e_1,\ldots,e_n$ be a local orthonormal frame field given
 by eigenvectors of $A_H$ satisfying
$A_He_i=\mu_ie_i,i=1,\ldots,n$. Then (1.4) implies
$\mu_i^2=\rho^2.$ Therefore, $A_H$ has at most two distinct
eigenvalues given by $\rho,-\rho$ with multiplicities, say
$k,n-k$, respectively. 

If $k=0$ or $k=n$, then $M$ is a pseudo--umbilical
submanifold with para\-llel mean curvature. Thus, by applying
Theorem 1.1 of Chapter II, $M$ is a minimal submanifold of a 
hypersphere of $\E^m$. 
 
If $0<k<n$, then, without loss of generality, we may put
$$\mu_1=\cdots=\mu_k=\rho,\quad \mu_{k+1}=\cdots=\mu_n=-\rho.
\leqno(1.5)$$

Because $M$ has parallel mean curvature vector, the mean 
curvature of $M$ is constant. Thus, (1.5) implies that
$\rho$ is a positive constant. Therefore, $M$ has homothetic Laplace
transformation. Hence, by Theorem 2.3 of Chapter IV, $M$ is a minimal
submanifold of a hypersphere of $\E^m$. 

The converse of this follows immediately from Theorem 2.3 of
Chapter IV. \sq

Theorem 1.2  implies the following

{\bf Corollary 1.3.}  {\it Let $M$ be a hypersurface of
 $\E^{n+1}$. Then $M$ is an open part of a
hypersphere in $\E^{n+1}$ if and only if $M$ has constant
mean curvature function and  conformal Laplace
transformation. \sq} 

It is easy to see that  the Laplace image of
 a minimal submanifold $M$ of a hypersphere of $\E^m$ is
 a minimal submanifold of a hypersphere of $\E^m$
centered at the origin. So it is natural to ask the
following problem. 

{\bf Problem 1.1.} {\it When is the Laplace image (with the induced
metric) of an isometric immersion $x: M \rightarrow \E^m$
 a minimal submanifold of a hypersphere of $\E^m$ 
centered at the origin? \sq
}

The following result gives a  solution to Problem 1.1.

{\bf Theorem 1.4.} {\it Let $x:M\rightarrow \E^m$ be an isometric
immersion of an $n$--dimensional Riemannian manifold into $\E^m$
with conformal Laplace transformation. Then the Laplace
image $L(M)$ (with respect to its induced metric) is a minimal
submanifold of a hypersphere of $\E^m$ centered at the origin
if and only if $\Delta H=fH$ for some nonzero function $f$ on $M$.
}

\demo 
Assume the Laplace transformation ${\cal L}:M
\rightarrow L(M)$ is conformal. Then the Laplace operator
$\bar \Delta$ of $L(M)$(with respect to its induced metric)
is given by $\bar\Delta=\rho^{-2}\Delta$
for some positive function $\rho$. Hence the mean curvature vector
$\bar H$ of $L(M)$ in $\E^m$ is given by
$\bar H=\rho^{-2}\Delta H.$

If $L(M)$ is immersed as a minimal submanifold of a hypersphere
of $\E^m$ centered at the origin, then the mean curvature vector
$\bar H$ of $L(M)$ in $\E^m$ satifies $\bar H=c\bar x$, where
$\bar x$ denotes the position vector field of $L(M)$ in $\E^m$ and
$c$ is  a nonzero constant. Since
$\bar x=Lx=-nH$, we have
$\Delta H=\rho^2\bar H=c\rho^2 \bar x=f H$,
where $f=-nc\rho^2\not=0$.

Conversely, if $\Delta H=fH$ for some non--zero function $f$, then 
we get 
$\bar H=\rho^{-2}\Delta H={f\over{\rho^2}}H= -{f\over{n\rho^2}}
\bar x$. This shows that the position vector field of $L(M)$ in
$E^m$ is normal to $L(M)$. Hence, $\<\bar x,\bar x\>$ is a
constant function on $L(M)$. Thus, $L(M)$ is contained in a
hypersphere $S^{m-1}(r)$ of $\E^m$ with radius $r$ and centered at
the origin. Let $\hat H$ denote the mean curvature vector of
$L(M)$ in $S^{m-1}(r)$. Then we have
$\bar H=\hat H-{1\over{r^2}}\bar x$. Because $\bar H$ is
parallel to $\bar x$, we get $\hat H=0$, {\it i.e.,} $L(M)$ is a
minimal submanifold of $S^{m-1}(r)$. \sq

Combining Lemma 1.3 of Chapter II and Theorem 1.4, we obtain 
the following.

{\bf Corollary 1.5.} {\it 
Let $x: M \rightarrow \E^m$
be an isometric immersion with conformal Laplace
transformation. Then the Laplace image
$L (M)$  is a minimal submanifold of a hypersphere of
$\E^m$ centered at the origin if and only if 
\begin{itemize}
\item[(1)]
$\hbox{trace}(\bar\nabla A_H)=0$ and 
\item[(2)] $\Delta^D H + a(H)$ is parallel
to $H$, where $a(H)$ denotes the allied mean curvature
vector field of $x$, $D$ denotes the normal connection
and $\Delta^D$ the Laplace operator of the normal
bundle. \sq
\end{itemize}
}

From Theorem 1.4 we also have the following.

{\bf Corollary 1.6.} {\it If $x : M \rightarrow \E^m$
is a null 2-type isometric immersion with conformal Laplace
transformation, then the Laplace image lies in a hypersphere
of $\E^m$ as a minimal surface.}

\demo This follows from the fact that every null 2--type
isometric immersion satisfies 
$\Delta H=\lambda H$ for some nonzero constant.  \sq

Now, we give the following

{\bf Lemma 1.7.} {\it Let $M$ be a hypersurface of
$\E^{n+1}$. If $M$ has conformal Laplace transformation,
then 
\begin{itemize}
\item[(1)]  the gradient
of the square of mean curvature function,
$\hbox{grad}\, \alpha^2$, is a
principal direction on the open subset 
$U=\{p\in M:\nabla\alpha^2\not=0 \,\hbox{at}\, p\}$ and 
\item[(2)] $M$ has at most {\rm 3} distinct principal 
curvatures of the following forms:
 $$\nu,\quad\sqrt{\nu^2 + |grad\,({\ln}\,\alpha
)|^2},\quad -\sqrt{\nu^2 + |grad\,({\ln}\,\alpha
)|^2},$$  
\end{itemize}
where the multiplity of $\nu$ is one on the open subset $U$.}

\demo If $M$ is a hypersurface of $\E^{n+1}$ with conformal
Laplace transformation, then Lemma 1.1 implies
$$\<A_HX,A_HY\>+(X\alpha)(Y\alpha)=\rho^2\<X,Y\>,\leqno(1.6)$$
for  $X,Y$ tangent to $M$, where $\rho$ is a positive function. 
Let $e_1,\ldots, e_n$ be a local orthonormal frame field given by 
eigenvectors of $A_H$ with eigenvalues $\mu_1,\ldots,\mu_n$,
respectively. Then, from (1.6), we have
$(e_i\alpha)( e_j\alpha)=0,\, i\not= j$. Without loss of generality,
we may assume 
$$e_2\alpha=\cdots=e_n\alpha=0.\leqno(1.7)$$
Then $e_1$ is parallel to the gradient of $\alpha^2$, which 
proves statement (1).

(2) If $M$ has constant mean curvature, then Corollary 1.3 implies
that $M$ has exactly one principal curvature. In this case statement
(2) holds. So, we assume that the mean curvature is not constant, 
{i.e.,} $U\not=\emptyset$.
From (1.6) we have
$$\mu_i^2+(e_i\alpha)^2=\mu_j^2+(e_j\alpha)^2,\quad i\not= j.\leqno
(1.8)$$
Combining (1.7) and (1.8), we obtain 
$$\mu_j^2=\mu^2+(e_1\alpha)^2,\quad \mu=\mu_1,
\quad j=2,\ldots,n.\leqno1.9$$
Let $\kappa_1,\ldots,\kappa_n$ be the principal curvatures of $M$
in $E^{n+1}$. Then $\mu_i=\alpha\kappa_i,i=1,\ldots,n$. Thus, 
from (1.9), we have statement (2).
\sq

If $M$ is a conformally flat hypersurface of dimension
$n > 3$ in $\E^{n+1}$, then, according to a well-known
result of \'E. Cartan, the shape operator of $M$ in
$E^{n+1}$ has an eigenvalue $\beta$ of multiplicity $\geq
n-1$. Let $\gamma$ denote the other eigenvalue of the
shape operator at the points where the multiplicity of
$\beta$ is $n-1$ and let $\gamma = \beta$ when the
multiplicity is $n$. Then the mean curvature function of
$M$ is given by $n\alpha = (n-1)\beta + \gamma$.

 The following result characterizes the class of conformally
flat hypersurfaces with conformal Laplace transformation.

{\bf Theorem 1.8.} {\it  Let $M$ 
be an $n$--dimensional $(n>3)$ conformally  flat hypersurface 
of $\E^{n+1}$. Then
$M$ has conformal Laplace transformation if and only if
\begin{itemize}
\item[(1)] the two eigenvalues of
the shape operator satisfy the following relation:
$$\alpha^2 \beta^2 = \alpha^2 \gamma^2 + | grad\, \alpha |^2$$ 
on $U=\{p\in M:\nabla\alpha^2\not=0 \,\hbox{at}\, p\}$,  and
\item[(2)] the gradients of the eigenvalues $\beta$
and $\gamma$ are parallel.
\end{itemize}
}

\demo (1) follows easily from statement (2) of Lemma 1.6.

(2) Let $e_1,\ldots,e_n$ be a local orthonormal frame field of
$M$ given by principal directions with 
$$A_{n+1}e_1=\gamma e_1,\quad A_{n+1}e_i=\beta e_i,\quad
i=2,\ldots, n.\leqno 1.10$$
Denote by $S$ and $\tau$ the Ricci tensor and the scalar curvature
of $M$, respectively. Then we have
$$S(X,Y)=\sum_{i=1}^n<R(e_i,X)Y,e_i>,\leqno(1.11)$$
$$\tau=\sum_{i=1}^n S(e_i,e_i),\leqno(1.12)$$
where $R$ is the Riemann curvature tensor of $M$.

We put
$$L(X,Y)=-{1\over{n-2}}S(X,Y)+{\tau\over{2(n-1)(n-2)}}g(X,Y).
\leqno(1.13)$$
Let $\hat L$ denote the (1,1)--tensor associated with $L$.
It is well--known that if $M$ is a conformally flat manifold
of dimension $n>3$, then we have (cf. for instance, [C2])
$$(\nabla_X\hat L)Y=(\nabla_Y\hat L)X,\leqno(1.14)$$
for $X,Y$ tangent to $M$.

From (1.10), it follows that the Weingarten map of $M$ in
$\E^{n+1}$ is given by
$$A_{n+1}=\beta I+(\gamma-\beta)\omega^1\otimes e_1,\leqno(1.15)$$
where $I$ is the identity map.

Moreover, (1.10) also implies
$$\hat L=-{{\beta^2}\over 2}I+\beta(\beta-\gamma)\omega^1\otimes
e_1.\leqno(1.16)$$
From the Codazzi equation, we have
$$(\nabla_X A_{n+1})Y=(\nabla_Y A_{n+1})X.\leqno(1.17)$$
(1.15) and (1.17) yields$$\hat L +\beta A_{n+1}={{\beta^2}\over
2}I.\leqno)1.18)$$ From (1.14), (1.15), (1.17) and (1.18), we obtain
$$(X\beta)A_{n+1}Y-\beta(X\beta )Y=(Y\beta)A_{n+1} X
-\beta(Y\beta)X.\leqno(1.19)$$ Setting $X=e_1, Y=e_i$ with
$i=2,\ldots,n$ in (1.19), we find $e_2\beta=\cdots=e_n\beta=0$.
Therefore, the 
gradient of $\beta$ is parallel to $e_1$. Consequently, by statement
(1), we  conclude that the gradients of the eigenvalues
$\beta$ and $\gamma$ are parallel. \sq

{\bf Remark 1.1.} Hypersurfaces of revolution $M$ in $\E^{n+1}$
automatically satisfy condition (2) of Theorem 1.8,
whereas condition (1) amounts to an ordinary differential
equation of order 3 to be satisfied by the planar profile
curve of $M$; hence, many hypersurfaces of revolution do
have conformal Laplace transformations. On the other
hand, for instance, the hypercatenoids $M$ in
$\E^{n+1}$, {i.e.} the hypersurfaces of revolution
whose profile curve is a catenary curve, do not satisfy
condition (1); hence, the hypercatenoids have
non-conformal Laplace transformations. \sq

By using Lemma 1.7 we obtain the following

{\bf Theorem 1.9.} {\it  Let $M$  be a
hypersurface of $\E^{n+1}\, (n>3)$. If $M$ has conformal Laplace
transformation and positive semi--definite  Ricci
tensor, then $M$ is conformally flat.
}

\demo If $M$ has conformal Laplace transformation, 
Lemma 1.4 implies that $M$ has at most 3 distinct principal
curvatures of the following forms:
 $$\nu,\delta, -\delta,\leqno 1.20$$
with mutliplicities $1, k$ and $n-k-1$, respec\-tively,
where $$\delta=\sqrt{\nu^2 + |grad\,({\ln}\,\alpha
)|^2}.$$   Without
loss of generality, we may assume $\nu>0$. 

If $k=0$ or $k=n-1$,  $M$ is a quasi--umbilical hypersurface.
In this case, $M$ is conformally flat. 

If $0<k<n-1$, then (1.20) implies that 
 the Ricci tensor $S$ of $M$ satisfies
$$ S(e_1,e_1)=(2k+1-n)\nu\beta,\quad 
S(e_n,e_n)=-\nu\beta+(n-2k-2)\beta^2.\leqno1.21$$
Since $M$ is assumed to have positive semi--definite
Ricci tensor,  (1.21) gives
$2k+2<n\leq 2k+1$ which is impossible.
\sq 

\vskip.2in
 \noindent{\bf  \S2. Surfaces in $\E^m$ with conformal Laplace
transformation. } \vskip.1in

In this section we investigate surfaces of general 
codimension with conformal Laplace transformation.

First we give the following  result.

{\bf Proposition 2.1.} {\it  Let $x: M \rightarrow \E^m$
be an isometric immersion of a Riemannian surface into
$E^m$. If the Laplace transformation ${\cal L} :M
\rightarrow L (M)$ is conformal, then the Laplace
image $L (M)$ (with respect to its induced
metric) is a minimal surface of $\E^m$ via $L$ if and
only if the immersion $x$ is biharmonic.
}

\demo Assume the Laplace transformation ${\cal L}:M
\rightarrow L(M)$ is conformal. Let $g$ and $\bar g$ denote
the Riemannian metrics of $M$ and $L(M)$, respectively. Then we have
${\cal L}^*\bar g=\rho^2g$ for some positive function $\rho$
on $M$. Denote the Laplace operator of $(M,g)$ and $(M,\bar g)$,
 respectively, by $\Delta$ and $\bar \Delta$. Then we have
$\bar\Delta=\rho^{-2}\Delta$. With respect to $\bar g$,
the Laplace map $L:(M,\bar g)\rightarrow L(M)\subset 
E^m$ is isometric. Hence, we have
$$-2\bar H=\rho^{-2}\Delta H,\leqno (2.1)$$
where $H$ and $\bar H$ denote the mean curvature vector of $x$
and of $L$, respectively. From (2.1) we conclude that $L(M)$ is a
minimal surface in $E^m$ if and only if the original immersion $x$
is biharmonic.
\sq

For flat ruled surfaces  in $\E^m$ we have the following. 

{\bf Proposition 2.2.} {\it  If $M$ is a flat ruled
surface in $\E^m$, then its Laplace transformation is not
conformal.
}

\demo As we mentioned in Chapter III,
 flat ruled surfaces in $\E^m$ are ``in general'' cylinders,
cones or tangential developables of curves. If $M$ is a cylinder
in $\E^m$,
the differential $dL$ (or $L_*$) of its Laplace map   
has rank $\leq 1$. So, its Laplace transformation 
is not conformal. 

Let $M$ be a cone in $\E^m$. Without
loss of generality, we may assume the vertex of the
cone is at
the orign and the cone is parametrized by
$$x(t,s)=t\beta(s),\quad |\beta|=|\beta'|=1.\leqno(2.2)$$ 
From the proof of Proposition 3.2 of Chapter III, we know
that the Laplace map of $M$ is a also cone given by
$$L(t,s)={1\over t}(\<\beta''(s),\beta(s)\>\beta(s)-\beta''(s)),
\leqno(2.3)$$
which implies
$$dL({{\partial}\over{\partial t}})=
\<\beta''(s),\beta(s)\>\beta(s)-\beta''(s),\leqno(2.4)$$
$$dL({{\partial}\over{\partial s}})
={1\over t}(\<\beta'''(s),\beta(s)\>\beta(s)
+\<\beta'',\beta\>\beta'-\beta'''(s)).\leqno(2.5)$$
(2.4) and (2.5) show that the metric tensor of $L(M)$ is given
by
$$\bar g_{11}=\<\beta'',\beta\>(1-2\<\beta'',\beta\>)+
|\beta'' |^2, $$
$$\bar g_{12}={1\over t}(\<\beta'',\beta'''\>-\<\beta'',\beta\>
\<\beta''',\beta\>),\leqno(2.6)$$
$$\bar g_{22}={1\over
{t^2}}(\<\beta'',\beta\>^2-\<\beta''',\beta\>^2
+|\beta'''|^2-\<\beta'',\beta\>\<\beta',\beta'''\>).$$

On the other hand, (2.2) implies that the metric tensor of $M$ is given
by
$$g_{11}=1,\quad g_{12}=0,\quad g_{22}=t^2.\leqno(2.7)$$
Comparing (2.6) and (2.7) we conclude that the Laplace transformation
of $M$ is not conformal.

If $M$ is a   tangential developable surface in $\E^m$ given by
the tangent lines of a unit speed curve $\beta(s)$ in $\E^m$, then
$M$ can be parametrized by
$$x(s,t)=\beta(s)+t\beta'(s),\quad |\beta'(s)|=1.\leqno(2.8)$$
Let $\kappa_1$ denote the first Frenet curvature and $\beta_2$
the second Frenet vector of $\beta$ in 
$\E^m$.  From the proof of Proposition 3.3 of Chapter III we know
that the Laplace map of $M$ is given by
$$L(s,t)=-{1\over {t\kappa_1}}\beta_2(s).\leqno(2.9)$$
By computing the metric tensors of $M$ and $L(M)$ and using (2.8)
and (2.9), we conclude that the Laplace transformation of 
a tangential developable surface is also not conformal.
\sq

 The following result shows that 
if a submanifold $M$  of $\E^m$
has a conformal Laplace transformation and if its
Laplace image $L (M)$ is a minimal surface
of a hypersphere of $\E^m$, then $M$ is not necessary a minimal
submanifold of a hypersphere of $\E^m$.

{\bf Proposition 2.3.} {\it  If $M$ is a surface in
$\E^6$ defined by
$$x(u,v) =(au,av, b\cos u,b\sin u, b\cos
v, b\sin v )$$
for some nonzero constants $a$ and $b$, then 
\begin{itemize}
\item[(1)] $M$ has homothetic Laplace transformation, 
\item[(2)] the Laplace image $L(M)$ is a minimal submanifold of a
hypersphere of $\E^6$ centered at the origin and, 
\item[(3)] the surface $M$ is not spherical.
\end{itemize}
}

\demo This follows from straight-forward computations.
\sq

As an application of Corollary 1.5  we give the
following

{\bf Theorem 2.3.} {\it Let $x: M \rightarrow \E^4$
be an isometric immersion with conformal Laplace
transformation. Then the Laplace image $L(M)$ of
 $x$ is a minimal surface of a hypersphere of
$\E^4$ centered at the origin if and only if $M$ is a
minimal surface of a hypersphere.
}

\demo  Assume $M$ is a minimal surface of a hypersphere
$S^3(r)$ of $\E^4$. Then without loss of generality we may assume
$S^3(r)$ is centered at the origin. In this case, we have
$L={2\over {r^2}}x$. This implies that the Laplace image
$L(M)$ lies in the hypersphere $S^3({2\over r})$ as a minimal
surface.

Conversely, assume that the Laplace transformation is conformal and 
the Laplace image $L(M)$ lies in a hypersphere centered at the
origin as a minimal surface. 
Then $M$ has constant mean curvature
in $\E^4$, since $L=-2H.$

 Let $e_1,e_2,e_3,e_4$ be a local orthonormal
frame field defined on $M$ such that $H=\alpha e_3$ and $e_1,e_2$ are eigenvectors
of $A_3$ with $A_3e_1=\kappa_1e_1,A_3e_2=\kappa_2e_2$. If
the Laplace transformation is conformal, then Lemma 1.1 implies
$(D_{e_1}e_3)(D_{e_2}e_3)=0$. Since $e_3 $ is perpendicular to
both $D_{e_1}e_3,D_{e_2}e_3$, at least one of 
$D_{e_1}e_3,D_{e_2}e_3$ vanishes. Without loss of generality, we
may assume $D_{e_1}e_3=0$. From Lemma 1.1 we find
$$ \alpha^2\kappa_1^2=\alpha^2\kappa_2^2+\phi^2,\quad \phi=
\omega^4_3(e_2).\leqno(2.10)$$

On the other hand, since $M$ has constant mean curvature, 
Corollary 1.5 implies
$$0=\hbox{trace}(\bar\nabla A_H)=\hbox{trace}(A_{DH})
=\phi A_4e_2.$$
Therefore, either $\phi=0$ or $A_4e_2=0$.

If $\phi=0$, then $DH=0$ and (2.10) implies $\kappa_1=\pm 
\kappa_2$. Since the Laplace transformation is conformal,
$\alpha$ is  nonzero. Thus, $\kappa_1=\kappa_2$. Hence, $M$ is
a pseudo--umbilical surface with parallel mean curvature vector.
By applying Theorem 1.1 of Chapter II, $M$ is contained in a 
hypersphere of $E^4$ as a minimal surface.

If $\phi\not=0$, then $A_4e_2=0$. Hence, by $\hbox{trace}\,A_4=0$, we get $A_4=0$.
So, by applying Codazzi's equation, we obtain
$\kappa_1=0$. Combining this with (2.10) we find
$\phi=0$, which is a contradiction. 
\sq
 
 \vskip.2in
\noindent{\bf  \S 3. Surfaces in $\E^3$ with conformal Laplace
transformation.} \vskip.1in

In this section we investigate surfaces in $\E^3$ with conformal
Laplace transformation. Since  such surfaces are abundant,
the complete classification of such surfaces is formidable.

{\bf Proposition 3.1.} {\it Let $M$ be a
surface in $\E^3$ with nonzero mean curvature function.
Then $M$ has conformal Laplace transformation if and
only if the following three conditions hold:
\begin{itemize}

\item[(1)]  the gradient
$\nabla \alpha^2$ of the mean curvature function is a
principal direction on the open subset $U$ of $M$ where
$\nabla\alpha\not= 0$;

\item[(2)] the Gauss curvature $K$ of
$M$ is given by $K={\alpha}^2 -{1\over{16\alpha^6}}
|\nabla\alpha |^4 $; and 

\item[(3)] the mean curvature function satisfies
the following equation: 
$$ \Delta \alpha^2 = 4\alpha^4 - 5|\nabla\alpha |^2 .$$
\end{itemize}
}

\demo If $M$ is a surface in $\E^3$ with conformal Laplace
transformation, then Lemma 1.7 implies that the gradient
$\nabla\alpha^2$ is a principal direction on $U=\{p\in M:
\nabla\alpha^2\not=0\;\hbox{at}\,p\}.$ Moreover, from
 Lemma 1.7,  $A_H$ also satisfies
$$A_He_1=\mu e_1,\quad A_He_2=\pm\sqrt{\mu^2+|\nabla\alpha |^2}e_2,
\leqno(3.1)$$
where $e_1$ is choosen in the direction of $\nabla\alpha^2$.
(3.1) yields
$$2\alpha^2=\mu \pm\sqrt{\mu^2+|\nabla\alpha |^2}.$$
 Thus we obtain
$$\mu=\alpha^2-{{|\nabla\alpha |^2}\over{4\alpha^2}},\leqno(3.2)$$
$$\mu^2+|\nabla\alpha |^2=(\alpha^2+
{{|\nabla\alpha |^2}\over{4\alpha^2}})^2.\leqno(3.3)$$
From (3.2) and (3.3) we get
$$A_3e_1=\Bigg(\alpha-{{|\nabla\alpha |^2}\over{4\alpha^3}}\Bigg)e_1,\quad
A_3e_2=\pm \Bigg(\alpha+{{|\nabla\alpha |^2}\over{4\alpha^3}}\Bigg)e_2.
$$
Because $\alpha\not=0$ and the two eigenvalues of $A_3$ satisfy
$2\alpha=\kappa_1+\kappa_2$, we must have
$$\kappa_1=\alpha-{{|\nabla\alpha |^2}\over{4\alpha^3}},\quad
\kappa_2=\alpha+{{|\nabla\alpha |^2}\over{4\alpha^3}}.
\leqno (3.4)$$

Therefore, the Gauss curvature $K$ of $M$ is given by
$$K={\alpha}^2 -{1\over{16\alpha^6}}
|\nabla\alpha |^4.\leqno(3.5) $$
Put
$\omega^2_1=f_1\omega^1+f_2\omega^2.$ Then, by
Codazzi's equation, (3.5), $e_2\alpha=0$, and a direct computation,
we find
 $$e_2e_1\alpha=f_1e_1\alpha,\leqno(3.6)$$
 $$e_1e_1\alpha={3\over{2\alpha}}(e_1\alpha)^2-f_2e_1\alpha-2\alpha^3.
 \leqno(3.7)$$
 From (3.6), (3.7) and $e_2\alpha=0$, we get
 $$\Delta\alpha=-e_1e_1\alpha-f_2e_1 \alpha=2\alpha^3-{3\over{2\alpha}}
 (e_1\alpha)^2.\leqno(3.8)$$
 Therefore
 $$\Delta\alpha^2=2\alpha\Delta\alpha-2 (e_1\alpha)^2=4\alpha^4
 -5|\nabla\alpha |^2.$$
 
 Conversely, assume $M$ is a surface in $\E^3$ satisfying
conditions
 (1) and (2) of the Proposition. Then we have
 $$e_2\alpha=0$$ on $U$ and 
 $$\kappa_1\kappa_2=\alpha^2 -{{(e_1\alpha )^4}\over{16\alpha^6}}
.\leqno(3.9)$$
Because $\kappa_1=2\alpha-\kappa_1$, (3.9) yields
$$\kappa_1=\alpha\pm {{(e_1\alpha)^2}\over{4\alpha^3}},\leqno(3.10)$$
$$\kappa_2=\alpha\mp {{(e_1\alpha)^2}\over{4\alpha^3}}.\leqno(3.11)$$
If  $\kappa_1=\alpha- {{(e_1\alpha)^2}\over{4\alpha^3}}$
$\kappa_2=\alpha+ {{(e_1\alpha)^2}\over{4\alpha^3}},$ then it is easy
to see that $\alpha^2\kappa_1^2+(e_1\alpha)^2=
\alpha^2\kappa_2^2$. Thus, by using $e_2\alpha=0$, we conclude
that $M$ has conformal Laplace transformation by Lemma 1.1.
If $\kappa_1=\alpha+ {{(e_1\alpha)^2}\over{4\alpha^3}}$
$\kappa_2=\alpha-{{(e_1\alpha)^2}\over{4\alpha^3}},$ then, by using
Codazzi's equation and $e_2\alpha=0$, we conclude that the mean
curvature does not satisfy condition (3) of the Proposition.
Consequently, we conclude that if conditions (1), (2) and (3) of
the Proposition hold, then $M$ has conformal Laplace transformation.
\sq

If $M$ is a compact surface, we have the
following characterization theorem.

 {\bf Theorem 3.2.} {\it   Let $M$ be a compact surface in
$\E^3$ whose  mean curvature function is nonzero. Then $M$ has
conformal Laplace transformation if and only if we have:
\begin{itemize}
\item[(1)]  $\nabla \alpha^2$  is a
principal direction on the open subset $U$ of $M$ where
$\nabla\alpha^2\not= 0$, and

\item[(2)] the Gauss curvature $K$ of
$M$ is given by $K={\alpha}^2 - {1\over{16\alpha^6 }}
|\nabla\alpha |^4 $.
\end{itemize}
} 

\demo Let $M$ be a compact surface in $\E^3$ whose mean
curvature function is nonzero. If $M$ has conformal Laplace
transformation, Proposition 3.1 implies (1) and (2) hold.

Conversely, if $M$ is compact and (1) and (2) hold, then from
the proof of Proposition 3.1 we know that either $M$ has
conformal Laplace transformation or $\alpha, \kappa_1,\kappa_2$
satisfy
$$\kappa_1=\alpha+ {{(e_1\alpha)^2}\over{4\alpha^3}},
\quad\kappa_2=\alpha-{{(e_1\alpha)^2}\over{4\alpha^3}},\quad
e_2\alpha=0.\leqno(3.12)$$

Put $\omega^2_1=f_1\omega^1+f_2\omega^2$ as before. Then, by (3.12),
Codazzi's equation and  direct computation, we obtain
$$e_2e_1\alpha=f_1e_1\alpha,\quad e_1e_1\alpha=2\alpha^3-f_2
e_1\alpha-{3\over 2}{{(e_1\alpha)^2}\over \alpha}.\leqno(3.13)$$
Thus, we get
$$\Delta\alpha=2\alpha^3+{3\over 2}{{(e_1\alpha)^2}\over \alpha}.
\leqno(3.14)$$
Consequently, we find
$$\Delta\alpha^2=2\alpha\Delta\alpha-2(e_1\alpha)^2
=4\alpha^4+|\nabla\alpha|^2.$$
Because $M$ is assumed to be compact, this implies $\alpha\equiv
0$ which is a contradiction.
\sq

The following result shows that examples of 
surfaces in $\E^3$ with conformal Laplace transformation
are abundant.

{\bf Proposition 3.3.}  {\it  There exist infinitely 
many surfaces of
revolution in $\E^3$ whose Laplace transformations are
conformal, but not homothetic.
}

\demo Let $M$ be a surface of revolution parametrized by
$$x(t,\theta)=(t,f(t)\cos\theta,f(t)\sin\theta).\leqno(3.15)$$
Then
$${\partial\over{\partial t}}=(1,f'(t)\cos\theta,f'(t)\sin\theta),
\quad 
{\partial\over{\partial \theta}}=(0,-f(t)\sin\theta,f(t)\cos\theta)
.\leqno(3.16)$$
From (3.15) and (3.16) we obtain
$$L(t,\theta)={{1+f'{}^2-ff'{}'}\over{f(1+f'{}^2)^2}}(-f',\cos\theta,
\sin\theta).\leqno(3.17)$$
By using (3.17) we can prove that the Laplace transformation is
conformal if and only if $f=f(t)$ satisfies the following third
order ordinary differential equation:
$$Af'''{}^2+Bf'''+C=0,\leqno(3.18)$$
where
$$A=f^4(1+f'{}^2)^3f'''{}^2,$$
$$B=2f^2(1+f'{}^2)\{(1+f'{}^2)(2ff'f''+f'(1+f'{}^2))-4f^2f'f''$$
$$+f'{}^3(1+f'{}^2)^2+ff'f''(f'{}^4+ff''-1-3ff'{}^2f'')\}$$
$$C=\{(1+f'{}^2)(2ff'f''+f'(1+f'{}^2))-4f^2f'f''{}^2\}^2\leqno(3.19)$$
$$+\{f'{}^2(1+f'{}^2)^2+ff''(f'{}^4+ff''-1-3ff'{}^2f'')\}^2$$
$$-(1+f'{}^2)^3(1+f'{}^2-ff'')^2.$$

It is easy to verify that, for infinitely many values of
 $a=f(0), b=f'(0), c=f''(0)$, the value of the 
discriminant $D=D(f,f',f'')=B^2-4AC$ at $t=0$ is positive.

Consider the third  order ordinary differential equation:
$$f'''=F(f,f',f'')={1\over {2A}}(-B+\sqrt{B^2-4AC}).\leqno(3.20)$$
For  $a=f(0), b=f'(0), c=f''(0)$ with $a\not=0$
and $D(a,b,c)>0$, the initial value problem:
$$f'''=F(f,f',f''), \quad f(0)=a,\quad f'(0)=b,\quad f''(0)=c$$
has a unique solution $f=f(t), $ in an open neighborhood of $0$.
 Such a solution $f=f(t)$ of the initial value problem gives
rise to a surface of revolution in $E^3$ with conformal 
Laplace transformation. Generically, such surfaces of revolution
do not have homothetic Laplace transformation .
\sq

Finally, we give the following.

{\bf Proposition 3.4.}  {\it  The Laplace image of a surface
$M$ in $\E^3$ with conformal Laplace transformation is not
a minimal surface in $\E^3$ (with respect to its induced
metric).
}

\demo Assume $M$ is a surface in $\E^3$ with conformal
Laplace transformation. Then the Laplace operator of the Laplace
image $L(M)$ satisfies $\bar\Delta=\rho^{-2}\Delta$ for some
positive function $\rho$ and the mean curvature vector
$\bar H$ of $L(M)$ in $\E^3$ is given by $\bar H=\rho^{-2}\Delta
H$. Thus, the Laplace image $L(M)$ is a minimal surface of $\E^3$
if and only if $M$ is a biharmonic surface. Since the only
biharmonic surfaces in $\E^3$ are minimal surfaces and
minimal surfaces in $\E^3$ do not have conformal Laplace
transformation, the Laplace image $L(M)$ cannot be a minimal
surface in $\E^3$. \sq

\vfill\eject

\noindent {\bf Chapter VI: GEOMETRY OF LAPLAE IMAGES}
\vskip.2in

\noindent {\bf   \S1. Immersions whose Laplace images lie in a cone. }
\vskip.1in

From the proof of Propositions 3.2 and 3.3 of Chapter III, we
know if $M$ is a cone in $\E^m$ with vertex at the origin or
a tangential developable surface, then the Laplace image of $M$ lies
in a cone in $\E^m$ with vertex at the origin. 
So, it is natural to ask the following problem:

{\bf Problem 1.1.}  {\it  ``When  the
Laplace image of a surface lies in a cone?''
} 

The purpose of this section is to investigate this Problem.

First,  we prove the following

{\bf Theorem 1.1.} {\it Let  $M$ be a surfce in $\E^3$ with 
regular Laplace map. Then the Laplace image
of $M$ lies in a cone with vertex at the origin if and
only if $M$ is locally a cone with vertex at the origin or
a tangential developable surface.
}

\demo Let $M$ be a surface in $\E^3$ with regular Laplace
map. If the Laplace image of $M$ lies in a cone with vertex at 
the origin. Then the Laplace image can be locally reparametrized 
by $$L(s,t)=t\beta(s).\leqno (1.1)$$
Since $L$ is regular, $s, t$ can be considered as
 local coordinates of $M$, too.  With respect to $(s,t)$ we have
$$ dL\Big({{\partial}\over{\partial s}}\Big)=2A_H\Big(
{{\partial}\over{\partial s}}\Big)-2{{\partial\alpha}\over{\partial s}}
e_3=t\beta'(s),$$ $$
dL\Big({{\partial}\over{\partial t}}\Big)=2A_H \Big(
{{\partial}\over{\partial t}}\Big)-2{{\partial\alpha}\over{\partial t}}
e_3=\beta(s),
\leqno (1.2)$$
where $H=\alpha e_3$ and $e_3$ is a unit normal vector of $M$
in $\E^3$. Because $L=-2H$, (1.1) implies that $\beta$
is normal to the surface $M$ in $\E^3$. Thus, (1.2) yields
$$A_H\Big({{\partial}\over{\partial t}}\Big)=0.\leqno (1.3)$$

Since $L$ is assumed to be regular, $H\not=0$. Thus, (1.3) implies
 $M$ is a flat surface in $\E^3$. Hence, $M$ 
is locally a cylinder, a cone or a tangential developable
surface. Because $L$ is regular, $M$ cannot be a cylinder. Therefore,
$M$ is locally a tangential developable surface or a cone.
If it is a cone with vertex at a point $p$, then from the
proof of Proposition 3.2 of Chapter III we see that the Laplace
image lies in a cone with vertex also at  $p$. So,
in this case, $M$ is a cone with vertex at the origin.

The converse follows from Propositions 3.2 and 3.3 of Chapter
III.
\sq

{\bf Remark 1.1.} In view of Theorem 1.1, it is interesting to point out
that there exist {\it non-flat\/} surfaces in $\E^3$ whose
Laplace images are contained in a cone with vertex at a
point other than the origin. In fact, by applying the
fundamental theorem of ordinary differential equations,
we can prove that there exist surfaces of revolution
which are non-flat and whose Laplace images are cones
with vertices at points other than the origin of $\E^3$. \sq

The following theorem yields a general result for
surfaces in $\E^m$ whose Laplace images lie in cone with
vertex at the origin.

{\bf Theorem 1.2.}  {\it  Let $x : M
\rightarrow \E^{m}$ be a surface in $\E^m$ with regular Laplace
map. If the Laplace image of $M$ lies in a cone with vertex at the
origin, then $M$ is non-positively curved, i.e., the
Gauss curvature $K$ of $M$ is $\leq 0$.
 Furthermore, 
 if $K=0$, then  locally $M$ is either a cone with vertex at the origin
or a tangential developable surface.
}

\demo Let $M$ be a surface in $\E^m$ with regular Laplace
map. If the Laplace image of $M$ lies in a cone with vertex at 
the origin, then the Laplace image can be locally parametrized 
by $$L(s,t)=t\beta(s), \quad |\beta|=|\beta'|=1.\leqno (1.4)$$
From (1.4) we may choose $e_3=\beta$ with $H=\alpha e_3$ and
$t=-2\alpha.$

Since $L$ is regular, $s, t$ can be considered as
 local coordinates of $M$, too.  With respect to $(s,t)$ we have
$$\aligned dL\Big({{\partial}\over{\partial s}}\Big)=2A_H\Big(
{{\partial}\over{\partial s}}\Big)+t D_{{\partial}
\over{\partial s}}
e_3=t\beta'(s),\\
dL\Big({{\partial}\over{\partial t}}\Big)=2A_H \Big(
{{\partial}\over{\partial t}}\Big)+e_3+tD_{{\partial}\over{\partial t}}
e_3=\beta(s).\endaligned
\leqno (1.5)$$
 (1.5) yields
$$A_H\Big({{\partial}\over{\partial t}}\Big)=0,\quad 
D_{{\partial\alpha}\over{\partial t}}e_3=0,
\leqno (1.6)$$
$$\beta'=-A_3\Big({\partial\over{\partial s}}\Big)+ D_{{\partial}
\over{\partial s}} e_3. \leqno (1.7)$$
By choosing $e_2$ in the direction of ${{\partial}\over{\partial t}}$,
and $e_1$ a unit tangent vector of $M$ perpendicular to $e_2$, 
we have
$$A_3 e_1=-t e_1,\quad A_3e_2=0.\leqno (1.8)$$
Since $e_3$ is parallel to the mean curvatur vector, we have
$$\hbox{trace}\,A_4=\cdots=\hbox{trace}\,A_m=0,\leqno (1.9)$$
where $A_r=A_{e_r}$ and $e_3,\ldots,e_m$ is a local orthonormal
frame field of the normal bundle of $M$ in $\E^m$. Thus, 
$\det(A_r)\leq 0,\,r=4,\ldots,m$.
Hence, by (1.8) and (1.9), we conclude that the Gauss curvature
$K$ of $M$ satisfies $K\leq 0$.

 If $K=0$, then (1.8) and (1.9) yield
 $$A_4=\cdots=A_m=0.\leqno (1.10)$$

From (1.6), (1.8) and (1.10) we have
$$\omega^3_1=-t\omega^1,\quad \omega^3_2=0,\quad \omega^r_i=0,
\quad i=1,2,\quad r=4,\ldots,m.\leqno (1.11)$$
Taking exterior derivative of $\omega^3_2=0$ and applying
(1.11) and the structure equations, we find
$\omega^2_1(e_2)=0$. Thus, $t$--curves on $M$ are geodesics.
Hence, by using (1.8) and (1.10), we conclude that $t$--curves
on $M$ are straight lines in $\E^m$. Hence, in this case, $M$
is a flat ruled surface. Therefore, locally, $M$ is a cylinder,
a cone, or a tangential developable surface. Because the Laplace
map is assumed to be regular, $M$ cannot be a cylinder. It $M$ 
is a cone, then its vertex must be at the origin, since the
Laplace image lies in a cone with vertex at the origin.
\sq

In view of Theorem 1.2, it is interesting to give the
following. 

{\bf Proposition 1.3.}  {\it There exist infinitely many
negatively curved surfaces in $\E^4$ whose Laplace
images lie in a cone with vertex at the
origin.   }

\demo Let $b$ and $c$ be real numbers such that $b>0$
and $c>4b^2$. Let
$D$ be the unit disk of $\E^2$ centered at $(0,1)$ with Riemannian
metric given by
$$g={1\over{4v^2}}du^2+{1\over{\mu^2v^2}}dv^2,\quad
\mu^2=v^2(c-b^2v^2).\leqno (1.12)$$ Put
$$e_1=-2v{{\partial}\over{\partial u}},\quad e_2=\mu
v{{\partial} \over{\partial v}}.\leqno (1.13)$$
Denote by $\omega^1,\omega^2$ the dual frame of $e_1,e_2$. Then
we have
$$\omega^1=-{1\over{2v}}du,\quad \omega^2={1
\over{\mu v}}dv,\quad \omega^2_1=\mu\omega^1.\leqno (1.14)$$
Let $F=D\times \E^2$ be the rank 2 trivial bundle over $D$. With
the usual Euclidean metric on fibres, $F$ is a Riemannian
vector bundle with usual connection $D$. Let $e_3,e_4$ be the
canonical orthonormal frame of the fibre $\E^2$. 

We define a bilinear map $h$ by
$$h(e_1,e_1)=2ve_3+bv^2e_4,\quad h(e_1,e_2)=0,\quad
h(e_2,e_2)=-bv^2e_4.\leqno (1.15)$$
Then, by direct computation, we can prove that $(D,g)$ together
with $D,h$ satisfy the equations of Gauss, Codazzi and Ricci.
Hence, by the fundamental theorem of submanifolds, we conclude
that there exists an isometric immersion of $(D,g)$ into $\E^4$
with $F$ as its normal bundle and $h$ as its second fundamental
form. 

Condition (1.15) says that the Gauss curvature of $D$ is
given by $K=-b^2v^2<0$. 

Now, we claim that the Laplace image of such a surface lies in
a cone with vertex at the origin. This can be seen as follows.

From (1.15) we know that the Laplace map of such a surface is
given by
$$L(u,v)=-2ve_3.\leqno (1.16)$$
From (1.15) and $De_3=0$ we get ${{\partial e_3}\over{\partial
v}}=0$, which shows that $e_3$ is a function of $u$. Therefore,
(1.16) implies that the Laplace image is contained in a cone
with vertex at the origin.
\sq

{\bf Proposition 1.4.}  {\it Every surface constructed in 
Proposition 1.3 is the locus of a planar curve $\gamma(s)$ moving 
along a space curve; moreover the planar curves are
congruent to  curves of the following form:
$$\Big({b\over c}\,\sinh^{-1}({{cs}\over b}), -{1\over
c}\sqrt{b^{2}+c^{2}s}\Big).\leqno (1.17)$$
where $b$ and $c$ are nonzero constants.
}

\demo First we observe from (1.13) and $\omega^2_1=\mu\omega^1$
in (1.14) that $v$--curves of the surface are geodesics.
Let $\gamma(v)$ be a $v$--curve with $s$ as its arclength
parameter.

Denote by $\tilde\nabla$ the Euclidean connection of $\E^4$.
From (1.14) and (1.15) we get
$$\tilde\nabla_{e_2}e_2=-bv^2e_4,\quad \tilde\nabla_{e_2}e_3=0,
\quad
\tilde\nabla_{e_2}e_4 =bv^2e_2.\leqno (1.17)$$
From (1.17) it follows that the $v$--curve $\gamma$
is a plane curve whose curvature function is given by
$$\kappa=-bv^2.\leqno (1.18)$$
From (1.12) and (1.13) we see that $s$ and $v$ are related by
$$s=\int {{dv}\over{v^2\sqrt{c-b^2v}}}=-{{\sqrt{c-b^2v}}\over{cv}}
+a.\leqno (1.19)$$
where $a$ is an integration constant. Without loss of generality, we
may assume $a=0$. From (1.18) and (1.19) we find
$$\kappa=-{{bc}\over{b^2+c^2s^2}}.\leqno (1.20)$$
Put $$\gamma\,'(s)=(\cos\theta(s),\sin\theta(s)).\leqno (1.21)$$
Then
$${{d\theta}\over{ds}}=-{{bc}\over{b^2+c^2s^2}}.\leqno (1.22)$$
So we obtain
$$\theta(s)=-\tan^{-1}({{cs}\over b})+c_2,\leqno (1.23)$$
where $c_2$ is a constant. Without loss of generality, we may assume
$c_2=0$. From (1.23) we have
$$\gamma\,'(s)=\Big({b\over{\sqrt{b^2+c^2s^2}}},-{{cs}\over{\sqrt{b^2+
c^2s^2}}}\,\Big).\leqno (1.24)$$
By taking integration of $\gamma\,'(s)$ with respect to $s$,
we obtain the Proposition.
\sq
The final result of this section is the following
theorem.

{\bf Theorem 1.5.}  {\it  Let $M$ be a non--flat surface in $\E^{m}$
with regular Laplace map and parallel normalized
mean curvature vector. If the Laplace image of $M$ lies in a
cone with vertex at the origin,   then $M$ is a locus of planar
curves moving along a space curve, moreover, the planar curves are
congruent to  curves of the following form:
$$\Big({b\over c}\,\sinh^{-1}({{cs}\over b}), -{1\over
c}\sqrt{b^{2}+c^{2}s}\,\Big),$$
where $b$ and $c$ are real numbers with  $b>0$.
}

\demo
Let $M$ be a surface in $\E^m$ with regular Laplace
map and parallel normalized mean curvature vector. If the Laplace
image of $M$ lies in a cone with vertex at  the origin, then the
Laplace image can be locally parametrized  by
$$L(s,t)=t\beta(s), \quad |\beta|=|\beta'|=1.\leqno (1.25)$$ From
(1.25), we may choose $e_3=\beta$ with $H=\alpha e_3$ and
$t=-2\alpha.$ Since $M$ has parallel normalized mean curvature
vector, $De_3=0$.

Since $L$ is regular, $s, t$ can be considered as
 local coordinates of $M$, too.  With respect to $(s,t)$ we have
$$\aligned dL\Big({{\partial}\over{\partial s}}\Big)=2A_H\Big(
{{\partial}\over{\partial s}}\Big)=t\beta'(s),\\
dL\Big({{\partial}\over{\partial t}}\Big)=2A_H\Big(
{{\partial}\over{\partial t}}\Big)+e_3=\beta(s),\endaligned
\leqno (1.26)$$
The first equation in (1.26) shows that $\beta'$ is tangent to
$M$ in $\E^m$. The second equation in (1.26) yields
$$A_H\Big({{\partial}\over{\partial t}}\Big)=0.\leqno (1.27)$$
By choosing $e_2$ in the direction of ${{\partial}\over{\partial
t}}$, and $e_1$ a unit tangent vector of $M$ perpendicular to
$e_2$,  we have
$$A_3 e_1=-t e_1,\quad A_3e_2=0.\leqno (1.28)$$

Since $\beta'(s)$ is tangent to $M$ and ${{\partial\beta}\over
{\partial t}}=0$, the equation of Gauss yields
$$\nabla_{{\partial}\over{\partial t}}\beta'=0,\quad
h(\beta',e_2)=0. \leqno (1.29)$$
Put ${{\partial }\over{\partial
s}}=a_1e_1+a_2e_2.$ Then
$$t\beta'=2A_H\Big({{\partial }\over{\partial
s}}\Big)=2a_1A_He_1+2a_2A_He_2=4a_1\alpha^2e_1.$$
Hence
$${{\partial}\over{\partial
s}}={1\over{t}}e_1+a_2e_2, \quad
e_1=\beta'.\leqno (1.30)$$ 
(1.28) and $De_3=0$ imply that $M$ has flat normal
connection, {\it i.e.,} $R^D=0$. Hence, we may choose
$e_3,\ldots,e_m$ such that $$\omega^4_1=\eta \omega^1,\quad
\omega^4_2=-\eta\omega^2, \quad A_5=\cdots=A_m=0.\leqno (1.31)$$
From (1.28), (1.29) and (1.31) we have
$$h(e_1,e_1)=2\alpha e_3+\eta e_4,\quad h(e_1,e_2)=0, \quad 
h(e_2,e_2)=-\eta e_4.\leqno (1.32)$$
Also, from (1.28) and (1.30), we obtain $\omega^2_1(e_2)=0$,
which implies $t$--curves of $M$ are geodesics.  We put
$$\omega^2_1=f\omega^1.\leqno (1.33)$$ 
From (1.32), (1.33), $De_3=0$ and the equation of Codazzi, we
 obtain
$$e_1\eta=0, \quad e_2\eta=2f\eta, \quad e_2\alpha=f\alpha,\leqno (1.34)$$
$$\eta D_{e_2}e_4=\eta D_{e_1}e_4=0.\leqno (1.35)$$
Because $M$ is assumed to be non--flat, $\eta\not=0$. So, (1.35)
implies
$$De_4=0.\leqno (1.36)$$

By using  (1.32), (1.33) and (1.36) we  obtain
$$\tilde\nabla_{e_2}e_1=0,\quad \tilde\nabla_{e_2}e_2= -\eta
e_4, \quad \tilde\nabla_{e_2}e_3=0,\quad \tilde\nabla_{e_2}
e_4=\eta e_2.\leqno (1.37)$$

Because $d\omega^2=-\omega^2_1\wedge\omega^1=0$ by (1.33),
$\omega^2=d\sigma$ for some function $\sigma$. We put
$$d\sigma=\sigma_s ds+\sigma_tdt.\leqno (1.38)$$
Since $\omega^2(e_1)=0$, (1.30) yields
$a_2=\sigma_s$. Therefore, (1.30) gives
$${{\partial}\over{\partial s}}={1\over t}e_1+\sigma_se_2.\leqno (1.39)$$ Similarly, we may also obtain
$${{\partial}\over{\partial t}}=\sigma_te_2.\leqno (1.40)$$
Because $[{{\partial}\over{\partial s}},
{{\partial}\over{\partial t}}]=0$, (1.39) and (1.40) imply
$\sigma_t={1\over{tf}}$. Consequently, $e_2=tf
{{\partial}\over{\partial t}}$. Therefore, the second equation
in (1.34) implies
$\eta=bt^2$ for some nonzero constant $b$ along a $t$--curve
in $M$. Combining this with (1.37) we may conclude as in the
proof of Proposition 1.4 that $t$--curves are plane curves
which are congruent to the curve mention in the Proposition.
\sq

\vskip.2in
\noindent {\bf  \S2. Laplace images of  surfaces of revolution.}
\vskip.1in

In this section we study the Laplace images
of surfaces of revolution.

First we make the following obveration.

{\bf Proposition 2.1.}  {\it Let $M$ be a surface of
revolution in $\E^3$ about an axis $L$. If $M$ is
non--minimal, then the Laplace map of $M$ is also a surface
of revolution about the same axis.
}

\demo Let $M$ be  a surface of revolution in $\E^3$. 
Without loss of generality, we may assume the axis of $M$ is
the $x$--axis. So, $M$ can be parametrized by  
$$x(t,\theta)=(t,f(t)\cos\theta,f(t)\sin\theta).\leqno(2.1)$$
By a direct computation, we see that the Laplace operator  is
given by 
$$\Delta = -{1\over {1+{f'}^2 }}{{\partial}^{2}\over
{\partial t^2}} - {{{f'}+{f'}^3 -ff'f''}\over {f{(1+{f'}^2
)}^2 }} {\partial \over {\partial t}} - {1\over {f^2}}
{{\partial}^2 \over {\partial {{\theta}^2 }}}.\leqno(2.2)$$ 
By  (2.2) we may conclude that the Laplace map
of the surface in $\E^3$ is given by 
 $$L(t,\theta)=\Bigg({{1+{f'}^2
-ff''}\over {f{(1+{f'}^2 )}^2}}\Bigg)(-f',\cos\theta,
\sin\theta).\leqno(2.3)$$ 
From (2.3) it follows that the
surface of revolution $M^2$ in $\E^3$ is minimal if and
only if $f$ satisfies the differential equation $
1+{f'}^2-ff''=0$. If $M^2$ is non-minimal, then,  by applying
(2.3), we see that the Laplace map of $M$ is also a surface
of revolution with the same axis. \sq 

The following result classifies surfaces of revolution
whose Laplace image lies in a cylinder.

{\bf Proposition 2.2.}  {\it  Let $M$ be a surface of
revolution in $\E^3$ with the $x$--axis as its axis. Then the
Laplace image of $M$ in $\E^3$ lies in a  cylinder if and only
if, after mutliplying the Euclidean coordinates of $\E^3$ by
a suitable constant, $M$ is of the following form:
$$x(t,\theta)=\Big(-{t\over 2}\sqrt{a^2t^{2}-1} + {1\over
{2a}} {\ln}(at+\sqrt{a^2t^{2}-1})-b, t\cos
\theta, t\sin \theta\Big)\leqno (2.4)$$ 
for some suitable constant $a,b$.
}

\demo Assume $M$ is a surface of revolution
parametrized by 
$$x(u,\theta)=(u,f(u)\cos\theta,f(u)\sin\theta).\leqno (2.5)$$
Then the Laplace map of $M$ is given by
 $$L(u,\theta)=\Big({{1+{f'}^2
-ff''}\over {f{(1+{f'}^2 )}^2}}\Big)(-f',\cos\theta,
\sin\theta).\leqno(2.6)$$ 
Formula (2.6) implies that the Laplace image of $M$ lies in
a cylinder if and only if
$$1+{f'}^2-ff''=cf(1+{f'}^2)^2\leqno (2.7)$$
where $c$ is a nonzero constant.  By applying a suitable
constant to each coordinate, we may obtain $c=1$ in the
follwing. So, we may assume $c=1$. Thus, (2.7) yields
$$ff''+{f'}^2(1+{f'}^2)=0.\leqno (2.8)$$
By solving  this second order ordinary differential
equation, we obtain the solution of (2.8) which is given by
$$u=-{f\over 2}\sqrt{a^2f^2-1}+{1\over{2a}}
\ln(af+\sqrt{a^2f^2-1})-b,\leqno (2.9)$$
where $a,b$ are integration constants.
(2.4) now follows from (2.5) and (2.9) by letting
 $t=f(u)$. 

The converse can be verified by direct computation.
\sq

From the definition of Laplace map, it is easy to
see that the Laplace image of a surface in $\E^3$ with
non--zero constant mean curvature lies in a sphere. So it is
natural to ask the following

{\bf Problem 2.1.}  {\it  ``Beside surfaces of constant mean
curvature, do there exist  surfaces in $E^3$ whose
 Laplace image  lies in some sphere?'' }
 
 Concerning this problem we have the following solution.

{\bf Proposition 2.3.}  {\it There exist infinitely 
many surfaces of revolution in $\E^3$ with non--constant mean
curvature whose Laplace images are contained in a sphere
which passes through the origin of $\E^3$. }

\demo Let $M$ be a surface of revolution parametrized
by (2.5). Then  the Laplace map of $M$
is given by (2.6) which
 implies that if the Laplace image of $M$ lies in a 
sphere, then the center of the sphere $S^2$ must be a point
on the $x$--axis. Assume $S^2$ is centered at $(r,0,0)$ and
with radius $R$. Then we have
$$\left( {{f'(1+f'{}^2-ff'')}\over{f(1+f'{}^2)^2}}-r\right)^2
+\left( {{1+f'{}^2-ff''}\over{f(1+f'{}^2)^2}}\right)^2=R^2,$$
which is equivalent to
$$(1+f'{}^2-ff'')^2-4rff'(1+f'{}^2)(1+f'{}^2-ff'')
+4(r^2-R^2)f^2 (1+f'{}^2)^3=0.$$
In particular, if the sphere passes through the origin, then
we have
$$(1+f'{}^2-ff'')(1+f'{}^2-ff''-4rff'(1+f'{}^2))=0.\leqno (2.10)$$
Since $M$ is minimal if and only if $1-f'{}^2-ff''=0$,
 the hypothesis yields
$$1+f'{}^2-ff''-4rff'(1+f'{}^2)=0.\leqno (2.11)$$ Put
$$F(u,f,v)={1\over f}(1+v^2-4rfv(1+v^2)).$$
It is easy to see that  $F,{{\partial
F}\over{\partial f}},{{\partial F}\over {\partial v}}$ are
continuous functions on any open set of the $(u,f,v)$--space
on which $f\not=0$. 

By applying the existence theorem
of second order ordinary differential equations, we know that
for any given point $x_0$ with $f(x_0)\not=0$, 
differential equation (2.11) together with  initial
conditions: $f(x_0)=f_0, f'(x_0)=f'_0$ has a unique local
solution about $x_0$. Generically, for such a solution
$f=f(u)$, the corresponding surface of revolution defined by
(2.5) does not have constant mean curvature.
Moreover, the Laplace image of such a surface lies in the
sphere centered at $(r,0,0)$ with radius $r$. \sq

 Concerning totally geodesic Laplace images,
we have the following  

{\bf Theorem 2.4.}  {\it Let $M$ be a surface of revolution in
$\E^3$. Then the Laplace image of $M$ lies in a plane if and
only if, up to rigid motions of $\E^3$, the surface is either
a circular cylinder or a surface of revolution given by
 $$x(t,\theta)=\Big({1\over a}\ln |at+\sqrt{a^2t^2-1}|,
t\cos\theta,t\sin\theta\Big),$$
where $a$ is a positive constant.
 }

\demo Assume $M$ is a surface of revolution parametrized by
(2.5). Then the Laplace map of $M$ is given by
 $$L(u,\theta)=\Bigg({{1+{f'}^2
-ff''}\over {f{(1+{f'}^2 )}^2}})(-f',\cos\theta,
\sin\theta\Bigg).\leqno(2.12)$$ 
From (2.12) we see that the Laplace image of $M$ lie in a
plane $P$ if and only if
$$f'({1+{f'}^2-ff''})=cf(1+{f'}^2 )^2\leqno (2.13)$$
for some constant $c$.

Without loss of generality, 
we may assume that the plane $P$ is the $yz$--plane and hence 
$c=0$. Then, from (2.13), we have
$$f'(ff''-1-f'{}^2)=0.\leqno (2.14)$$

If $f'=0$, then $M$ is a circular cylinder. In this case,
the Laplace image is a circle in the plane $P$. 

If $f'\not=0$, then (2.14) yields
$$ff''=1+f'{}^2.\leqno (2.15)$$
From (2.15) we get
$$f'{}^2=a^2f^2-1\leqno (2.16)$$
for some constant $a>0$. Solving (2.16) we get
$$u=\pm {1\over a}\ln |af+\sqrt{a^2f^2-1}|+b\leqno (2.17)$$
where $b$ is a constant. 
Without loss of generality, we may assume
$$u= {1\over a}\ln |af+\sqrt{a^2f^2-1}|.\leqno (2.18)$$
By putting $f(u)=t$, (2.5) and (2.18) show that the
surface of revolution is given by
 $$x(t,\theta)=\Big({1\over a}\ln |at+\sqrt{a^2t^2-1}|,
t\cos\theta,t\sin\theta\Big).$$

The converse follows from direct computations.
\sq

\vskip.2in

\noindent {\bf  \S3. Laplace images of curves.}

In this section we study  Laplace image of
curves.

It is clear that the Laplace image of a circle is a circle and
the Laplace image of a line is a point. So, it is natural to ask
the following.

{\bf Problem 3.1.}   {\it ``Besides circles and lines in $\E^2$, do
there exist other planar curves whose Laplace images lie in a
circle?'' 
}

The following result gives an affirmative answer to this
Problem.

{\bf Proposition 3.1.}  {\it The Laplace image of a
unit speed planar curve $\gamma(s)$ lies in a circle if and
only if locally it is either a line, a circle, or, up to
similarity transformations of $\E^2$, it is a curve given by
$$\gamma(s)=\Big(s-{1\over 2}\ln(1+c^2e^{4s}),\tan^{-1}(ce^{2s})\Big),
\leqno ( 3.1)$$
  where $c$ is a positive constant. }

\demo Let $\gamma(s)=(x(s),y(s))$ be a unit speed curve
in the $xy$--plane and $\theta(s)$ be the angle between
$\gamma\,'(s)$ and the positive direction of the $x$--axis.
Then we have
$$\gamma\,'(s)=(\cos\theta(s),\sin\theta(s)).\leqno (3.2)$$
Because the Laplace operator of $\gamma$ is given by
$\Delta=-{{d^2}\over{ds^2}}$, the Laplace map of $\gamma$ is
given by 
$$L(s)={{d\theta}\over{ds}}(\sin\theta(s),-\cos\theta(s)).
\leqno (3.3)$$

Assume the Laplace image lies in a circle. By applying a
suitable
similarity transformation on the plane,  the center of the
circle can be chosen to be $(1,0)$ and its radius to
be one.  Therefore, by using (3.3) we obtain 
$$(\theta'(s))^2-2(\sin\theta)\theta'(s)=0.\leqno (3.4)$$

If $\gamma$ is not contained in a line or in a circle, then
$\theta'(s)\not=0$. Thus (3.4) yields $\theta'(s)=2\sin\theta$.
By solving this differential equation, we obtain
$$\theta(s)=2\tan^{-1}(ce^{2s})\leqno (3.5)$$
for some positive constant $c$. Combining (3.2) and (3.5),
we may obtain 
$$\gamma\,'(s)=\Big({{1-c^2e^{4s}}\over{1+c^2e^{4s}}},{{2ce^{2s}}
\over{1+c^2e^{4s}}}\Big),$$
from which we obtain (3.1).

The converse is easy to verify.
\sq

Similar to Problem 3.1, we would like to know whether there
exist curves whose Laplace images lie in a line.

{\bf Proposition 3.2.}   {\it Let $\beta (s)$ be a curve in
$\E^m$ parametrized by arclength.  Then the Laplace image of
$\beta$ lies in a line if and only if either 
\begin{itemize}
\item[(1)] $\beta$ is a
planar curve whose curvature function $\kappa (s)$
satisfies the differential equation: 
$$\kappa {\kappa}'' -{\kappa}^4 =3({\kappa}' )^2\leqno (3.6)$$
or
\item[(2)] $\beta$ is a helix in an
affine 3-space whose first and second Frenet curvature
functions $\kappa_1$ and $\kappa_2$ satisfy the equations:
$$\kappa_2 = c\kappa_1 ,\quad \kappa_1 {\kappa_1}'' -
(1+c^2 ){\kappa_1}^4 =3({\kappa_1}')^2,\leqno (3.7)$$
where $c$ is a nonzero constant.\end{itemize}
}

\demo The Laplace map of the unit speed curve $\beta$ is
given by  $L(s)=-\beta''(s)$.

(1) If $\beta$ is a planar curve, then 
$$L'(s)=\kappa^2\beta_1-\kappa'\beta_2,\quad
L''(s)=3\kappa\kappa'\beta_1+(\kappa^3-\kappa'')
\beta_2,\leqno (3.8)$$ 
where $\kappa$ is the plane curvature and
$\beta_1=\beta'$ and $\beta_2$ is a unit normal vector field.
Because the Laplace image of $\beta$ is contained in a line if
and only if $L'(s)$ and $L''(s)$ are linearly dependent, 
(3.8) implies  statement (1).

(2)  Let $\beta$ be a unit speed space curve in $\E^m$ with
$m\geq 3$. If $m=3$, we may consider $\E^3$ as a linear subspace
of $\E^4$. Thus, without loss of generality, we may assume
$m\geq 4$.
Denote the $i$--the Frenet curvature by $\kappa_i$ and the
$i$--the Frenet vector by $\beta_i$. Then we have
$$L'(s)=\kappa_1^2\beta_1-\kappa'_1\beta_2-\kappa_1\kappa_2
\beta_3,\leqno (3.9)$$
$$\aligned
L''(s)=3\kappa_1\kappa'_1\beta_1-(\kappa''_1-\kappa_1^3
-\kappa_1\kappa_2^2)\beta_2 \\
-(2\kappa'_1\kappa_2-\kappa_1\kappa'_2)\beta_3
-\kappa_1\kappa_2\kappa_3\beta_4.\endaligned\leqno (3.10)$$
From (3.9) and (3.10) we conclude that the Laplace image of
the space curve $\beta$ lies in a line if and only if the
following four conditions hold:
$$\kappa_1\kappa''_1-\kappa_1^4-\kappa_1^2\kappa_2^2
=3(\kappa'_1)^2,\leqno (3.11)$$
$$\kappa'_1\kappa_2=\kappa_1\kappa'_2,\leqno (3.12)$$
$$2\kappa'_1{}^2\kappa_2+\kappa_1\kappa'_1\kappa'_2 +\kappa_1^4
\kappa_2+\kappa_1^2\kappa_2^3=\kappa''_1\kappa_1
\kappa_2,\leqno (3.13)$$ and 
$$\kappa_1^2\kappa_2\kappa_3=\kappa_1^2\kappa'_1
\kappa_2\kappa_3
=\kappa_1^2\kappa_2^2\kappa_3=0.\leqno (3.14)$$
 (3.12) implies $\kappa_2=c\kappa_1$ for some constant $c$. Thus,
from (3.11) as well as from (3.13), we obtain (3.7). From
(3.14), we obtain $\kappa_3=0$. Thus $\beta$ is a curve in an
affine 3--space with $\kappa_2=c\kappa_1$. Hence, $\beta$ is a
helix. 

Conversely, if $\beta$ is a unit speed helix in $\E^3$ whose
Frenet curvatures satisfy (3.7), then conditions (3.11) and
(3.12) hold. Condition (3.14) holds automatically since $\beta$
is a curve in $\E^3$. Furthermore, because $\kappa_2=c\kappa_1$,
equation (3.13) is nothing but the second equation in (3.7).
Thus, the Laplace image of $\beta$ lies in a line.
\sq

\vskip.2in
\noindent {\bf  \S4. Laplace images of totally real submanifolds.}
\vskip.1in

Let $x:M\rightarrow {\Bbb C}^m$ be a map from
an $n$--dimensional Riemannian manifold $M$ into the complex
number $m$--space ${\Bbb C}^m$ equiped with the standard flat
Kaehler metric. Denote by $J$ the almost complex structure
of ${\Bbb C}^m$. The map $x$ is said to be a {\it complex map
\/} if each tangent space of $M$ is invariant
under the action of $J$, {\it i.e.,}
$$J(dx(T_pM))\subset dx(T_pM),\quad p\in M.$$ 
The map $x$ is said to be a {\it totally real map\/} if
$$J(dx(T_pM))\subset (dx(T_pM))^{\perp},\,p\in M,$$
 where $(dx(T_pM))^{\perp}$ is the orthogonal complementary
subspace  of $dx(T_pM)$ in $\Bbb C^m$.

For a totally real isometric immersion  from $M$ into  ${\Bbb
C}^m$, it is natural to ask the following question:

{\bf Problem 4.1.}  {\it ``When is the Laplace image of a
totally real submanifold  totally real?''
}

Concerning this question we have the following results.

{\bf Lemma 4.1.}   {\it  Let $x: M \rightarrow {\Bbb C}^n$ be a
totally real isometric immersion of an $n$--dimensional
Riemannian manifold $M$ into ${\Bbb C}^n$. Then the Laplace map
$L:M\rightarrow {\Bbb C}^n$
 of the immersion $x$ is totally
real  if and only if the shape operator
$A_H$ in the direction of $H$ and the normal connection $D$ of
$x$ satisfy the condition: 
$$\< A_{H}X,JD_Y H\>=\< A_H Y,JD_XH\>,\leqno (4.1)$$ 
for any $X,Y$ tangent to $M^n$. 
}

\demo The Laplace map of the isometric totally real
immersion $x$ is given by $L=-nH$.  For any vector $X$ 
tangent to $M$, we have  $$dL(X)=nA_HX-nD_XH.$$
Thus the Laplace map $L$ is totally real if and only if
$$\<A_HX-D_XH,JA_HY-JD_YH\>=0\leqno (4.2)$$ for all $X,Y$ tangent to
$M$. Since the immersion $x$ is totally real and $m=n$, (4.2) is
equivalent to condition (4.1).
\sq

Recall that a submanifold $M$  is said to have
parallel normalized mean curvatue vector if $M$ has nonzero
mean curvature and the unit normal vector field in the
direction of the mean curvature vector field is a parallel
normal vector field.

{\bf Theorem 4.2.}  {\it  Let $x: M \rightarrow {\Bbb C}^m$ be a
totally real isometric immersion of an $n$--dimensional
Riemannian manifold into ${\Bbb C}^m$. Then  
\begin{itemize} 
\item[(1)] if the mean curvature vector field $H$ of
$x$ is parallel in the normal bundle, then the Laplace map is 
totally real; and 
\item[(2)] if  $M$ has parallel
normalized mean curvature vector field and $m=n$, then the
Laplace map is totally real if and only if 
$A_{H}JH$ is parallel to the gradient of the mean
curvature function, $\nabla\,\alpha.$
\end{itemize} 
}

\demo Let $x: M \rightarrow {\Bbb C}^m$ be a totally
real isometric immersion of an $n$--dimensional Riemannian
manifold into ${\Bbb C}^m$ for which $H$ is parallel in the
normal bundle. Then the Laplace map satisfies $dL(X)=nA_HX$.
Thus, for any $X$ tangent to $M$, $dL(X)$ is also a tangent
vector. Since $x$ is a totally real isometric immersion, this
implies $\<JdL(X),dL(Y)\>$ $=0$ for any $X,Y$ tangent to $M$.
Hence, by definition, the Laplace map is totally real.

(2) Assume $x$ has parallel normalized mean curvature vector.
Then the mean curvature function is nonzero, {\it
i.e.,} $\alpha\not=0$, and $D\xi=0$, where $\xi$ is the unit
normal vector field in the direction of $H$. Thus, for any
vector $X$ tangent to $M$, we have
$$dL(X)=nA_HX-n(X\alpha)\xi.\leqno (4.3)$$ Hence, the Laplace map is
a totally real map if and only if 
$$(Y\alpha)\<A_HX,J\xi\>=(X\alpha)\<A_HY,J\xi\>,\leqno (4.4)$$ for
any $X,Y$ tangent to $M$.

If we choose a local orthonormal frame field $e_1,\ldots,e_n$
of $M$ such that $e_1$ is in the direction of $\nabla\alpha$,
then  $e_2\alpha=\cdots=e_n\alpha=0$. Let $U$ be the
open subset of $M$ on which $\nabla\alpha\not=0$. Then (4.4)
implies $$\<A_He_i,J\xi\>=0,\,i=2,\ldots,n,$$ which is
equivalent to $$\<e_2,A_H(J\xi)\>=\cdots=\<e_n,A_H(J\xi)\>=0.$$
Therefore, $A_H(JH)$ is parallel to $\nabla\alpha$.

Conversely, if $A_H(JH)$ is parallel to $\nabla\alpha$, then
$$\<A_He_i,J\xi\>=0,\,i=2,\ldots,n,\leqno (4.5)$$
where $e_1$ is in the direction of $\nabla\alpha$. Because
$e_2\alpha=\cdots=e_n\alpha=0$, (4.5) implies (4.4) 
for any $X,Y$ tangent to $M$. Since $M$ is assumed to have
parallel normalized mean curvature vector, from (4.4) we see
that  condition (4.1)  holds, too. Therefore, by applying Lemma
4.1, we may conclude that the Laplace map of $x$ is totally real.
\sq

{\bf Remark 4.1.} Statement (1) of Theorem 4.2 implies that
there exist infinitely many totally real submanifolds in ${\Bbb
C}^m$ whose Laplace map are also totally real. \sq

\vfill\eject

\noindent {\bf Chapter{VII}: SUBMANIFOLDS WITH HARMONIC
 LAPLACE MAPS AND TRANSFORMATIONS}

\vskip.2in

\noindent{\bf   \S1. Submanifolds with harmonic Laplace map. } \vskip.1in

Let $\phi :M\rightarrow N$ be a differentiable map between
Riemannian manifolds $M$ and $N$. Denote by $\nabla$ and
$\tilde\nabla$ the Levi--Civita connections of $M$ and $N$,
respectively. Then the {\it second fundamental form\/} $h_\phi$
and  the {\it energy density\/}  $e(\phi)$ of the map $\phi$ are
given respectively  by
$$h_\phi(X,Y)=\tilde\nabla_X(\phi_*Y)-\phi_*(\nabla_XY),\leqno (1.1)$$
$$e(\phi) = {1\over 2} ||d\phi ||^2 = {1\over 2} trace({\phi}^{*}
g'),\leqno(1.2)$$ where $X,Y$ are tangent vectors of $M$, $g'$ is
the  Riemannian metric on $N$, $\phi_*=d\phi$, and $\phi^*$ is the
induced map of $\phi$.

 If $M$ is compact, the
{\it energy\/} $E(\phi)$ of $\phi$ is defined by $$E(\phi) =\int_M
e(\phi)*1.$$ The Euler-Lagrange operator
associated with $E$ will be written
$\tau(\phi)=div(d\phi)$ and called the {\it tension
field\/} of $\phi$. A map is said to be {\it harmonic\/}
if its tension field vanishes identically. For a
differentiable map $\psi : M^n \rightarrow \E^m$, one has
(cf. [EL1,2]) $$\Delta \psi =-\tau(\psi).\leqno(1.3)$$
Furthermore, let $i: N \rightarrow P$ be an isometric
immersion and $\Phi :M \rightarrow P$ be the composition
of $\phi$ and $i$. Then the tension field $\tau(\phi)$
is the orthogonal projection of $\tau(\Phi)$ onto the
tangent bundle $T(N)$. More precisely, we have (cf.
[EL1,EL2]) $$\tau(\Phi) =\tau(\phi) + \hbox{trace}\,h(d\phi,d\phi
),\leqno(1.4)$$
where $h$ is the second fundamental form of $N$ in $P$.
In particular, $\phi : M \rightarrow N$ is harmonic if
and only if $\tau(\Phi)$ is perpendicular to $N$. 

Concerning harmonic maps, it is natural to ask the
following problem:

{\bf Problem 1.1.} {\it ``When is the Laplace map  $L :M^n
\rightarrow \E^m$ of an isomtric immersion  $x:
M^n\rightarrow \E^m$  harmonic?'' \sq
}

From Beltrami's formula and formula
(1.3), it follows that the Laplace map $L :M^n
\rightarrow \E^m$ of the immersion $x$ is a harmonic map if
and only if the immersion $x$ is biharmonic. Biharmonic
submanifolds were first studied by the first author
as one of the off--springs of the theory of finite type. In
fact, the first author obtained in 1985
(unpublished then, see [Di1] for details) that if $x : M^2
\rightarrow \E^3$ is an isometric immersion, then the mean
curvature vector field $H$ is harmonic, {\it i.e.,\/} $\Delta
H =0$, (or equivalently, $x$ is biharmonic: $\Delta^2 x =0$)
if and only if the immersion is minimal.   In his 1989 
doctoral thesis at Michigan State
University, I. Dimitric  generalized this result to some
more general classes of submanifolds (cf.
[Di1]). Furthermore,  it is shown in [CI1] that space--like
biharmonic surfaces in the Minkowski space--time $\E^3_1$ are
minimal surfaces, too.

In terms of the Laplace map, this provides us the
following first solution  to Problem 1.1.

{\bf Theorem 1.1.} {\it Let $M$ be a
surface in $\E^3$. Then the Laplace map $L
:M^2 \rightarrow \E^3$ is a harmonic map if and only if
$M^2$ is a minimal surface.}

\demo  Let $x:M \rightarrow \E^3$ be a surface whose
Laplace map is a harmonic map, then the isometric immersion
$x$ is biharmonic. We put $H=\alpha e_3$.
Then, from Theorem 1.3 of Chapter II, we know that
$$\Delta H=\Delta^DH+||h||^2H+\nabla \alpha^2+2A_3(\nabla
\alpha),\leqno (1.5)$$
where $||h||$ is the length of the second fundamental form
$h$.
Since $M$ is biharmonic, we have 
$$\Delta\alpha+||h||^{2}=0,\leqno(1.6)$$
$$A_{3}(grad\,\alpha)=-\alpha (grad \,\alpha).\leqno(1.7)$$
 We choose $e_{1},e_2$ which diagonalize $A_3$.
Thus we have $$h(e_{1},e_{1})=\beta e_{3},\,\,\,
h(e_{1},e_{2})=0,\,\,\,h(e_{2},e_{2})=\gamma
e_{3}\leqno(1.8)$$ for some functions $\beta, \gamma$. From
the equation  of Codazzi and (1.8) we have $$e_{2}\beta =
(\beta - \gamma )\omega_{1}^{2}(e_{1}),\,\,\,e_{1}\gamma =
(\beta - \gamma)\omega_{1}^{2}(e_{2}).\leqno(1.9)$$ 
Let $U=\{\, p\in M\, \vert \,\nabla\alpha^{2} \not= 0\,\, at\,
\,p\,\}.$
 Then $U$ is an open subset of $M$. Assume $U\not=
\emptyset .$ Then, by (1.9), $\nabla\alpha$ is an eigenvector
of $A_3$ with eigenvalue $-\alpha$ on $U$. We
choose $e_1$ in the direction of $\nabla\alpha$ on $U$. Then
we have $e_{2}\alpha = 0$ on $U$. Moreover, the Weingarten map
$A_3$ satisfies 
$$\beta = -\alpha ,\quad \gamma = 3\alpha,
\quad\ ||h||^{2}=10\alpha^{2}.\leqno(1.10)$$
Since $e_{2}\alpha = 0,$ (1.9) and (1.10) imply
$$\omega_{1}^{2}(e_{1})=0,\,\,\,
d\omega^{1}=0.\leqno(1.11)$$ Thus, locally,
$\omega^{1}=du$ for some function $u$. Since
$d\alpha\wedge\omega^{1}=d\alpha\wedge du = 0,\, \alpha\,$
is a function of $u$. We denote by $\alpha'$ and
$\alpha''$ the first two derivatives of $\alpha$ with
respect to $u$. From (1.9), (1.10), and (1.11), we have
$$4\alpha\omega_{1}^{2}=-3\alpha'\omega^{2}.\leqno(1.12)$$
Also from $e_{2}\alpha=0$ and (1.12), we have
$$4\alpha\Delta\alpha =
3(\alpha')^{2}-4\alpha\alpha''.\leqno(1.13)$$ By using
(1.6), (1.10), and (1.13), we get
$$4\alpha\alpha''-3(\alpha')^{2}-40\alpha^{4}=0.\leqno(1.14)$$
From (1.14) we may obtain
$$(\alpha')^{2}=8\alpha^{4}+C\alpha^{3/2}\leqno(1.15)$$
for some constant $C$. On the other hand, the equation  of
Gauss, (1.8) and (1.12) imply $$\alpha\alpha'' -{7\over 4}
(\alpha')^{2}+4\alpha^{4}=0.\leqno(1.16)$$
From (1.14) and (1.16) we get
$$(\alpha')^{2}=16\alpha^{4}.\leqno(1.17)$$
Combining (1.15) and (1.17) we conclude that $\alpha$ is
constant on $U$ which is a contradiction. Therefore, $U$
is empty. So, $M$ has constant mean curvature. Thus, by
applying (1.6), we obtain $\alpha = 0$. Hence
the immersion $x$ is a minimal immersion. 
\sq

Next, we consider the following problem.

{\bf Problem 1.2.} {\it ``When is the Laplace
transformation  ${\cal L} :M \rightarrow L(M )$ of an
isometric immersion  $x: M\rightarrow \E^m$  harmonic?''}

Since a (weakly) conformal map between two
Riemannian surfaces is a harmonic map,  Proposition 3.3 of
Chapter V shows that there exist many surfaces of revolution 
which have harmonic Laplace transformations and which have
non-constant mean curvature functions.
Therefore, this seems to be a difficult problem in
general. However, for curves in $\E^m$, this Problem is much
easier to solve.

First, it is easy to see that lines and circles in a plane
have harmonic Laplace transformations. If the curve is
neither a line or a circle, we have the following results.

{\bf Proposition 1.2.} {\it  Let $\beta(s)$ be a  curve
parametrized by arclength $s$. If $\beta$ is neither a line
or a circle, then the Laplace transformation of $\beta$ is a
harmonic map if and only if $\beta$ is a curve whose first and
second Frenet curvature functions satisfy the  
equation:
$${\kappa_1}'(2{\kappa_1}^3 +{\kappa_1}'' )+
{\kappa_1} {\kappa_2}({\kappa_1}' {\kappa_2}+{\kappa_1}
{\kappa_2}') = 0.\leqno (1.18)$$

 In particular, if $\beta$ is a planar curve, then
the Laplace transformation of $\beta$ is a harmonic
map if and only if the curvature function $\kappa$
of $\beta$ satisfies the differential equation: $\kappa
''=-2\kappa^3$.}

\demo Let $\beta(s)$ be a unit speed curve. Then the
Laplace map $L$ is given by $L(s)=-\beta''(s).$
Thus 
$$\aligned\Delta L=-3\kappa_1\kappa'_1\beta_1-(\kappa^3-\kappa''_1
+\kappa_1\kappa^2_2)\beta_2\\
+(2\kappa'_1\kappa_2+\kappa_1
\kappa'_2)\beta_3+\kappa_1\kappa_2\kappa_3\beta_4,\endaligned\leqno (1.19)$$
where $\beta_i$ is the $i$--the Frenet vector of $\beta$. 
On the other hand, we have
$$dL({d\over{ds}})=\kappa_1^2\beta_1-\kappa_1'\beta_2
-\kappa_1\kappa_2\beta_3.\leqno (1.20)$$
It is known that the Laplace transformation $\cal L$ of
$\beta$ is a harmonic map if and only if the component of
$\Delta L$ in the direction of $dL({d\over{ds}})$ vanishes.
Therefore, by using (1.19) and (1.20), we may conclude that
the Laplace transformation ${\cal L}$  is a harmonic
map if and only if 
$$-3\kappa^3_1\kappa'_1+\kappa'_1(\kappa_1^3-
\kappa''_1+\kappa_1\kappa^2_2)+ \kappa_1\kappa_2(-2
\kappa'_1\kappa_2-\kappa_1\kappa'_2)=0.\leqno (1.21)$$
It is easy to verify that (1.21) is equivalent to (1.18).

If $\beta$ is a planar curve, then (1.18) is equivalent to
$\kappa'(2\kappa^3+\kappa'')=0$ where $\kappa$ is the plane
curvature of $\beta$. Because $\beta$ is neither a line or a
circle, this yields $\kappa''=-2\kappa^3.$ 

The converse is easy to verify.
\sq
 
From this Proposition we know that the class of planar
curves whose Laplace transformations are harmonic maps
depends on two parameters. Moreover, by
combining Proposition 3.3 and Remark 3.1 of Chapter IV and the
last Proposition, we have the following

{\bf Corollary 1.3.} {\it  The Laplace transformation
of a planar curve $\beta$ is homothetic if and only if
the curve $\beta$ has nonzero curvature function and it has
harmonic Laplace transformation.}

\demo From Proposition 3.3 and Remark 3.1 of Chapter
IV we know that a unit speed planar curve $\beta(s)$ has
homothetic Laplace transformation if and only if the plane
curvature of $\beta$ satisfies
$$(\kappa'(s))^2+\kappa(s)^4=c^2,\leqno (1.22)$$
for some positive constant $c$. By taking the derivative of
(1.22) with resepct to $s$, we obtain
$$\kappa'(\kappa''+2\kappa^3)=0.\leqno (1.23)$$
If $\kappa'=0$, then (1.22) implies that $\beta$ is a circle
which has harmonic Laplace transformation. If $\beta$ is not
a circle, then (1.23) yields $\kappa''=-2\kappa^3$. Thus, by
applying Proposition 1.2, we conclude that the Laplace
transformation of $\beta$ is a harmonic map.

Conversely, if $\beta$ is a planar curve which has nonzero
curvature function and harmonic Laplace transformation, then
either $\beta$ is a circle or the plane curvature of $\beta$
satisfies $\kappa''=-2\kappa^3$. By solving this differential
equation, we may obtain (1.22), which implies $\beta$ has
homothetic Laplace transformation.
\sq

For submanifolds with parallel mean curvature vector field,
we have the following result.

 {\bf Theorem 1.4.} {\it  Let $x: M
\rightarrow \E^{m}$ be an isometric immersion from an
$n$--dimensional Riemannian manifold into $\E^m$ with parallel
mean curvature vector field. Then the Laplace transformation
of $x$ is a harmonic map. }

\demo From Theorem 1.3 of Chapter II we have the following
fundamental formula:
 $$\Delta H=\Delta^{D}H + \sum_{i=1}^{n}
h(A_{H}e_{i},e_{i}) +  {{n}\over {2}}grad<H,H> + 2\,
trace\,A_{DH},\leqno(1.24)$$ where
 $\Delta^D$ is the Laplacian
associated with the normal connection $D$.
If $x$ has parallel mean curvature vector field, then (1.24)
 yields
$$\Delta L= -n\sum_{i=1}^{n} h(A_{H}e_{i},e_{i}).\leqno (1.25)$$
On the other hand, from $DH=0$, we have
$$L_*(X)=nA_H(X),\leqno (1.26)$$ for any vector $X$ tangent to $M$. 
Comparing (1.25) and (1.26) we conclude that the Laplace
transformation of $x$ is a harmonic map.
\sq

\vskip.2in

\noindent{\bf  \S2. Submanifolds with harmonic mean curvature function.} \vskip.1in

 In this section, we
would like to consider submanifolds  whose {\it mean
curvature functions are harmonic.\/} From Hopf's lemma we
know that if a compact submanifold has harmonic mean
curvature function, then the  mean curvature
function is constant. So in this section we are only
interested in non-compact submanifolds.  

First we give the following 

{\bf Lemma 2.1.} {\it The only planar curves with harmonic curvature
function are open parts of a line,  a circle, or a Cornu
spiral. (By  a Cornu spiral curve, we mean a planar curve whose
curvature function in terms of an arc--length parameter $s$ is
given by $\kappa(s)=as+b$ for some real numbers, $a\not=0$ and
$b$.) }

\demo The Laplace operator of a unit speed curve $\gamma(s)$ is
given by $\Delta=-{{d^2}\over{ds^2}}$. Thus, $\gamma$ has harmonic
curvature function if and only if $\kappa''(s)=0$, or
equivalently, $\kappa(s)= as+b$ for some constants $a,b$. If
$a=b=0$, the curve is an open part of a line. If $a=0$ and
$b\not=0$, $\gamma$ is an open part of a circle. If $a\not=0$,
$\gamma$ is an open part of a Cornu spiral curve. \sq

Let $M$ be a flat surface in $\E^3$. Then, locally, $M$ is a
cylinder, a cone, or a tangential developable surface.  We shall
consider these three cases separately.

{\bf Theorem 2.2.} {\it  The only hypercylinders in
$\E^{n+1}$ with harmonic mean curvature function are
open parts of the following hypersurfaces:
\begin{itemize}
\item[(1)] hyperplanes; 

\item[(2)] circular hypercylinders; and 

\item[(3)] hypercylinders $C\times \E^{n-1}$,
where $C$ is a Cornu spiral.
\end{itemize}
}

\demo Let $M$ be a  hypercylinder in $\E^{n+1}$
 parametrized by
$$x(s,t_2,\ldots,t_n)=(u(s),v(s),t_2,\ldots,t_n),\leqno (2.1)$$
where $\gamma(s)=(u(s),v(s))$ is a planar curve with  $s$ as
its arclength parameter.  Then the Laplace operator of
$M$ is given by $$\Delta=-{{\partial^2}\over{\partial s^2}}
-\sum_{i=2}^n {{\partial^2}\over{\partial t_i^2}}.\leqno (2.2)$$
It is easy to see that the mean curvature function
$\alpha$ of $M$ in $\E^{n+1}$ is related with the plane curvature
function $\kappa$ of $\gamma$ by $\kappa(s)=n\alpha$. Thus, by the
expression of the Laplace operators of $M$ and $\gamma$, we know
that the hypercylinder $M$ has harmonic mean curvature function
if and only if $\gamma$ has harmonic mean curvature function.
Thus, by applying Lemma 2.1, we obtain the Theorem. \sq

{\bf Lemma 2.3.} {\it Open parts of a plane are the only tangential
developable surfaces in $\E^3$ with harmonic mean curvature function.
}

\demo Let $M$ be a tangential developable surface parametrized
by $$x(s,t)=\beta(s)+t\beta'(s),\leqno (2.3)$$ where $\beta(s)$ is a unit
speed curve in $\E^3$. Then, by a direct computation, we have
$$\aligned\Delta= -{1\over{t\kappa_1}}\Big\{ {{\partial}\over{\partial
s}} \Big({1\over{t\kappa_1}}{{\partial}\over{\partial
s}}\Big)-{{\partial}\over{\partial s}} ({1\over{t\kappa_1}}
{{\partial}\over{\partial t}}\Big)\\ -
{{\partial}\over{\partial t}}\Big({1\over{t\kappa_1}} 
{{\partial}\over{\partial s}}\Big)+
{{\partial}\over{\partial t}}\Big({{1+t^2\kappa_1^2}\over{t\kappa_1}}
{{\partial}\over{\partial t}}\Big)\Big\}.\endaligned\leqno (2.4)$$

From (2.3) we obtain
$${{\partial x}\over{\partial s}}=\beta_1+t\kappa_1\beta_2,\quad
{{\partial x}\over{\partial t}}= \beta_1.\leqno (2.5)$$
Thus, $\beta_1,\beta_2$ form an orthonormal frame field on
$M$. Therefore $\xi=\beta_3$
is a unit normal vector field of the surface. By direct computation,
we have 
$$A_\xi \beta_1=0,\quad A_\xi
\beta_2={{\kappa_2}\over{t\kappa_1}}\beta_2. \leqno (2.6)$$
This implies that the mean curvature function of $M$ in $\E^3$ is
given by 
$$\alpha={{\kappa_2}\over{2t\kappa_1}}
.\leqno (2.7)$$
By using (2.4) and (2.7) we conclude that $\Delta\alpha=0$ if and
only if  the Frenet curvatures of $\beta$ satisfy
$$\aligned(\kappa_1^2\kappa''_2-\kappa_1 \kappa''_1\kappa_2-3\kappa_1
\kappa'_1\kappa'_2+3{\kappa'_1}^2  \kappa_2 +\kappa_1^4\kappa_2)t^2\\
+(2\kappa_1^2\kappa'_2-4\kappa_1 \kappa'_1\kappa_2+\kappa_1^2
\kappa_2)t+3\kappa_1^2\kappa_2=0.\endaligned\leqno (2.8)$$
Since (2.8) holds for any $t$, we obtain $\kappa_2=0$. Therefore,
$\beta$ is a planar curve; and hence the tangential developable
surface is an open part of a plane. \sq

{\bf Lemma 2.4.} {\it  Besides open parts of  planes, the only
other cones in $\E^3$ with harmonic mean curvature function are
those whose rays cut the unit sphere $S^{2}_{p}(1)$  centered at
$p$ in a curve $\beta(s)$ of which the curvature in
$S^{2}_{p}(1)$ is given by $c_1 \cos s + c_2 \sin s$, where $s$
is an arclength parameter of $\beta$ and $c_1$ and $c_2$ are
two nonzero constants. }

\demo It is clear that open parts of planes have harmonic
mean curvature functions and circular cones do not have harmonic
mean curvature functions.

Let $M$ be a cone parametrized by
$$x(s,t)=t\beta(s),\quad <\beta,\beta>=<\beta',\beta'>=1.\leqno (2.9)$$
We assume $M$ is neither an open part of a plane or of a circular
cone. We put $\rho=\kappa_1^{-1}$ and $\sigma=\kappa_2^{-1}$. It
is  well--known from classical differential geometry that the unit speed
spherical curve $\beta$ satisfies
$$\beta=-\rho\beta_2-(\sigma\rho')\beta_3.\leqno (2.10)$$
From (2.9) and a direct computation, we have
$${{\partial x}\over{\partial s}}=t\beta_1,\quad
{{\partial x}\over{\partial t}}=\beta.\leqno (2.11)$$
Put
$$\xi=(\sigma\rho')\beta_2-\rho\beta_3.\leqno (2.12)$$
Then $\xi$ is a normal vector field of $M$ in $\E^3$. Moreover, by
direct computation, we  obtain
$$\tilde\nabla_{{\partial}\over{\partial s}}\xi=-\kappa_1\sigma
\rho'\beta_1,\quad \tilde\nabla_{{\partial}\over{\partial t}}\xi=0.
\leqno (2.13)$$
Therefore, the mean curvature function of $M$ is given by
$$\alpha=-{{\kappa'_1}\over{2t\kappa_1\kappa_2}}.\leqno (2.14)$$
By direct computation, we also have
$$\Delta=-{1\over t}\{ {{\partial}\over{\partial s}}
({1\over t}{{\partial}\over{\partial s}})+
{{\partial}\over{\partial t}}(t{{\partial}\over{\partial t}})\}.
\leqno (2.15)$$
Therefore, by using (2.14) and (2.15), we conclude that $M$ has harmonic
mean curvature function if and only if the Frenet curvatures of
$\beta$ satisfy
$$({{\kappa'_1}\over{\kappa_1\kappa_2}})''+
{{\kappa'_1}\over{\kappa_1\kappa_2}}=0.\leqno (2.16)$$
From (2.9) and (2.10) we obtain
$$\kappa_1^2\kappa_2^2(\kappa_1^2-1)=\kappa'_1{}^2.\leqno (2.17)$$
Since the curvature function $\bar\kappa$ of $\beta(s)$ in
$S^2(1)$ is given by $\bar\kappa=\sqrt{\kappa_1^2-1}$, (2.17)
yields $$\bar\kappa=\pm
{{\kappa'_1}\over{\kappa_1\kappa_2}}.\leqno (2.18)$$ After solving the
second order differential equation  (2.16), we find from (2.18)
that $\bar\kappa(s)=c_1\cos s+c_2\sin s$ for some constants
$c_1,c_2$. \sq

{\bf Theorem 2.5.} {\it  Let $M$ be a flat surface in
$\E^3$. Then $M$ has harmonic mean curvature function if
and only if $M$ is given by  open parts of the following
surfaces:
\begin{itemize}
\item[(1)] the planes;

\item[(2)] the circular cylinders;

\item[(3)] the Cornu cylinders $C\times {\Bbb R}$,
i.e., the right cylinders on the Cornu spiral curves;
\item[(4)] the cones with vertex at a point $p$
whose  rays cut the unit sphere $S^{2}_{p}(1)$  centered
at $p$ in a curve $\beta(s)$ of which the curvature in
$S^{2}_{p}(1)$ is given by $c_1 \cos s + c_2 \sin s$,
where $s$ is an arclength parameter of $\beta$ and $c_1$
and $c_2$ are two nonzero constants.
\end{itemize}
}

{\bf Proof.} This Theorem follows easily from Theorem 2.2, Lemma 2.3 and
Lemma 2.4. \sq

{\bf Remarks 2.1.} 
\begin{itemize}

\item[(i)] Surfaces given in (1) and (2) of 
Theorem 2.5 are of course surfaces with {\it constant\/}
mean curvature; surfaces in (3) and (4) have
{\it non-constant} harmonic mean curvature. \sq

\item[(ii)] The Cornu spirals are planar curves with
harmonic curvature; the curves $\beta(s)$ mentioned in (4)
are spherical curves with harmonic curvature. \sq

\item[(iii)] In view of Thereom 2.5, it is interesting to
point out that {\it the only complete flat surfaces in $\E^3$
with harmonic mean curvature functions are planes, circular
cylinders and Cornu cylinders.}

\end{itemize}

We recall that Delaunay's surfaces in $\E^3$ are surfaces
of revolution with constant mean curvature function. For
surfaces of revolution we have the following result.

{\bf Theorem 2.6.} {\it We have the following:
\begin{itemize}

\item[(1)] If $M$ is a surface
of revolution in $\E^3$ defined by
$$x(t,\theta ) = (t, f(t)\cos\theta ,f(t)\sin\theta
),\leqno(2.19)$$
then $M^2$ has harmonic mean curvature function if and
only if the generating function $f$ satisfies the
 following
 ordinary differential equ\-ation of third order:
$$ (1+{f'}^2 )\{ f^2 f'''+ f'  + ff' f'' \}
-3f^2 f' {f''}^2 + cf(1+ {f'}^2 )^3 =0,\leqno(2.20)$$
 where c is an arbitrary constant; and

\item[(2)] besides  the surfaces of Delaunay,
there exist infintely many surfaces of revolution in $\E^3$ with
harmonic mean curvature function.
\end{itemize}
}

\demo Let $M$ be a surface of revolution in $\E^3$ defined by
(2.19). Then the Laplace operator and the Laplace map of $M$ are
respectively given by
$$\Delta = -{1\over {1+{f'}^2 }}{{\partial}^{2}\over
{\partial t^2}} - {{{f'}+{f'}^3 -ff'f''}\over {f{(1+{f'}^2
)}^2 }} {\partial \over {\partial t}} - {1\over {f^2}}
{{\partial}^2 \over {\partial {{\theta}^2 }}},\leqno(2.21)$$ 

 $$L(t,\theta)=({{1+{f'}^2
-ff''}\over {f{(1+{f'}^2 )}^2}})(-f',\cos\theta,
\sin\theta).\leqno(2.22)$$ 
Thus the mean curvature function of $M$ is given by
 $$\alpha={{1+{f'}^2
-ff''}\over {2f{(1+{f'}^2 )}^{3/2}}}.\leqno (2.23)$$
(2.21) and (2.23) imply that the mean curvature function is
harmonic if and only if 
$${f\over{\sqrt{1+{f'}^2}}}{d\over{dt}}\left( 
{{1+{f'}^2
-ff''}\over {f{(1+{f'}^2 )}^{3/2}}}\right) =c,\leqno (2.24)$$
where $c$ is a constant.

It is easy to verify that (2.24) is equivalent to (2.20). This
proves (1).

 (2) From the existence and uniqueness theorem of ordinary
differential equations, we know that for a given initial point
$t_0$, there exists a unique solution $f=\phi(t)$ defined on
some interval about $t_0$ of the differential equation
defined by (2.20) that satisfies the prescribed initial
conditions: $$f(t_0)=f_0,\quad f'(t_0)=f'_0,\quad
f''(t_0)=f''_0.\leqno (2.25)$$ Furthermore, by using (2.23), we see
that the class of surfaces of revolution with constant mean
curvature depends on three parameters:
$t_0,f(t_0),f'(t_0)$; and the class of surfaces of
revolutions with harmonic mean curvature function depends on
four parameters $t_0,f(t_0),f'(t_0),f''(t_0)$. From these we
conclude that, besides the surfaces of Delaunay, there exist
infinitely many surfaces of revolution in $\E^3$ with harmonic
mean curvature function. 
\sq

\vfill\eject

\noindent {\bf Chapter VIII: LAPLACE AND GAUSS IMAGES}
\vskip.2in

\noindent{\bf  \S1.  Laplace and Gauss images  of
submanifolds in $\E^m$. } \vskip.1in 

 Let $x: M^n \rightarrow \E^m$ be an
isometric immersion of an $n$-dimensional, oriented
Riemannian manifold $M^n$ into $\E^m$ and
$e_{n+1},\ldots,e_m$  a local oriented orthonormal frame
of the normal bundle $T^{\perp}(M )$. Then the {\it
Gauss map} $G :M^n \rightarrow G(n,m-n),$
defined by 
$$G (p) =(e_{n+1}\wedge\ldots\wedge e_m)(p),\quad
p\in  M,\leqno(1.1)$$
is a differentiable map from $M$ into the
real Grassmannian of oriented $(m-n)$-planes in $\E^m$. We
call the image $G (M^n )$ of the Gauss map $G$ the {\it Gauss
image\/} of the immersion $x$. Whenever the Laplace map of $x$
is an immersion, then, locally, we have a  map, denoted by $$LG
:L(M^n )  \rightarrow G (M^n),\leqno(1.2)$$ from the Laplace image
into the Gauss image defined by $LG(L (p))=G (p)$ for  $p\in
M$. We call this map the {\it Laplace-Gauss transformation\/}
 of $x$, or for short, {\it LG-transformation\/}.

 For simplicity,  we shall 
 assume throughout this section that the Laplace map
of  $x$ is an immersion. The purpose of this section is to
study the geometry of $LG$-transformations of such
immersions. 

First we consider curves in Euclidean spaces.

{\bf Proposition 1.1.} {\it Let $\beta$ be a curve in
$\E^m$ parametrized by arclength. Then
\begin{itemize}
\item[(1)]  if $\beta$ is a planar curve,
then the LG-transformation of $\beta$ is homothetic if and
only if its curvature function $\kappa$ is of the following
form:
$$\kappa(s)={1\over 2}\Big(\sqrt{a^2e^{-s^2}+4(c+c^2)}-ae^{-s}\Big),
\leqno(1.3)$$ for some positive constants $a,c$;

\item[(2)] the LG-transformation of
$\beta$ is homothetic if and only if the first and the
second Frenet curvatures of $\beta$ satisfy the
differential equation $$\kappa_{2}^2 = c\{\kappa_1^2
+{(\ln \kappa_1)'}^2\},$$ where $c$ is a positive number.
\end{itemize}
}

\demo  (1) If $\beta(s)$ is a planar curve  parametrized
by arclength, then the Laplace and Gauss maps are given
respectively by 
$$L(s)=-\kappa(s) \beta_2, G(s)=\beta_2(s).\leqno(1.4)$$
Since the Laplace map is assumed to be regular,
$\kappa(s)\not=0.$ Moreover, we have from (1.4) that
the induced metrics of the Laplace and Gauss images of $\beta$
are given respectively by
$$g_{_L}=\kappa^4+{\kappa'}^2,\quad g_{_G}=\kappa^2.\leqno(1.5)$$
Therefore, the $LG$--transformation of the plane curve is
homothetic if and only if ${\kappa'}^2=c^2\kappa^2-\kappa^4$ for
some positive number $c$. By solving this differential
equation, we obtain (1.3).

(2) If $\beta$ is a curve in $\E^m$ with $m\geq 3$, then the
induced metrics of the Lalace and Gauss images of $\beta$ are
given respectively by
$$g_{_L}=\kappa_1^4+{\kappa'}_1^2+\kappa_1^2\kappa_2^2,\quad
g_{_G}=\kappa_1^2.\leqno(1.6)$$
This implies that the $LG$--transformation is homothetic if and
only if the first and second Frenet curvatures of $\beta$ satisfy
$\kappa_{2}^2 =  c\{\kappa_1^2
+{(\ln \kappa_1)'}^2\}$. \sq

{\bf Remark 1.1.} Proposition 1.1 implies that there exist
infinitely many curves whose $LG$--transformation are
homothetic. In particular, statement (1) of Proposition 1.1
says that, up to Euclidean motions, the class of planar
curves with homothetic $LG$--transformation depends on two
parameters $a,c>0$. \sq
 
For hypersurfaces of dimension $>1$ we have the following result.

{\bf Theorem 1.2.} {\it  Let $x: M \rightarrow
\E^{n+1}$ be an isometric immersion. Then
 the LG-transformation $LG :L (M^n ) 
\rightarrow G (M^n)$ of $x$ is weakly conformal if and
only if, locally, either $M$
 has  constant mean curvature function in $\E^{n+1}$ or $M$ is
the product of a planar curve and an affine $(n-1)$--space
$\E^{n-1}$.
}

\demo Let $M$ be a hypersurface of $\E^{n+1}$ with $n>1$. 
We choose a local orthonormal frame field
$e_1, \ldots,e_n, e_{n+1}$ such that $e_1,\ldots,e_n$ are
tangent vectors given by eigenvectors of $A_{n+1}$ with
$A_{n+1}e_i=\kappa_ie_i,i=1,\ldots,n$. Then the Laplace and the
Gauss maps of $M$ satisfy
 $$L_*(e_i)=n\alpha\kappa_ie_i-(e_i\alpha)e_{n+1},\quad G_*(
e_i)=-\kappa_ie_i,\quad i=1,\ldots,n.\leqno(1.7)$$
Therefore, the induced metrics on the Laplace and the Gauss
images of $M$ are given by
$$g_{_L}(e_i,e_j)=n^2\alpha^2\kappa_i\kappa_j\delta_{ij}
+(e_i\alpha)
(e_j\alpha),\quad  g_{_G}(e_i,e_j)=\kappa_i\kappa_j\delta_{ij}
.\leqno(1.8)$$

If the $LG$--transformation is weakly conformal,
then (1.8) implies 
$$n^2\alpha^2\kappa_i\kappa_j\delta_{ij}+(e_i\alpha)(e_j
\alpha) =f^2\kappa_i\kappa_j\delta_{ij},\quad
i,j=1,\ldots,n,\leqno(1.9)$$ for some function $f$. In
particular, we have $(e_i\alpha)(e_j\alpha)=0$ when $i\not=
j$. Therefore, the gradient of $\alpha$, $\nabla\alpha$ is
parallel to an eigenvector of $A_{n+1}$. Without loss of
generality, we may assume $\nabla\alpha$ is parallel to
$e_1$. Then, by (1.9), we have
$$n^2\alpha^2\kappa^2_1+|\nabla\alpha|^2=f^2\kappa_1^2,
\quad (n^2\alpha^2-f^2)\kappa_i^2=0, \quad
i=2,\ldots,n.\leqno(1.10)$$

{\bf Case 1.} If at least one of $\kappa_2,\ldots,\kappa_n$ is
nonzero, then the second equation of (1.10) yields
$f^2=n^2\alpha^2$ and hence, by the first equation of
(1.10), $\alpha$ is constant.

{\bf Case 2.}
If $\kappa_2=\cdots=\kappa_n=0$, then 
$$A_{n+1}e_1=n\alpha e_1,\quad A_{n+1}e_2=\cdots= A_{n+1
}e_n=0.\leqno(1.11)$$
Put  $U=\{p\in M:H(p)\not=0\}$. If $U$ is an empty set, then
$M$ is an open part of a hyperplane. So, we assume
$U\not=\emptyset$ and let $W$ be a connected component of $U$.
We claim that $W$ is an open part of the product  of a planar
curve and an affine $(n-1)$--space $\E^{n-1}$. This can be
seen as follows.

Let ${\cal D}_1=\hbox{Span}\{e_1\}$ and ${\cal
D}_2=\hbox{Span} \{e_2,\ldots,e_n\}$ on $W$. If $X,Y$ are
vector fields in ${\cal D}_2$, then we have
$A_{n+1}X=A_{n+1}Y=0$. Thus, by the equation of Codazzi, we
have $$A_{n+1}([X,Y])=
(\nabla_X  A_{n+1})(Y)-(\nabla_Y  A_{n+1})(X)=0.\leqno(1.12)$$
Hence, ${\cal D}_2$ is integrable. Moreover, since ${\cal
D}_1$ is 1--dimensional, ${\cal D}_1$ is trivially integrable.

Furthermore, since $\nabla \alpha$ is parallel to $e_1$,
(1.12) and the equation of Codazzi imply that
$$0=(\nabla_{e_1}A_{n+1})(e_i)-(\nabla_{e_i}A_{n+1})(e_1) =
n\alpha\omega_1^i(e_1)e_1-n\alpha\nabla_{e_i}e_1,\leqno(1.13)$$
for $i=2,\ldots,n$.
Therefore
$$\nabla_{e_1}e_1=0,\quad \nabla_{e_i}e_1=0.\leqno(1.14)$$
(1.14) imply that integral submanifolds of ${\cal D}_1,
{\cal D}_2$ are totally geodesic submanifolds of $W$.
Hence, by the de Rham decomposition theorem, $W$ is locally
the Riemannian product of a geodesic of $M$ and an
$(n-1)$--dimensional totally geodesic submanifold of $M$.
Moreover, by (1.11) and a lemma of Moore, we may conclude
that $W$ is the product submanifold of a planar curve and an
affine $(n-1)$--space $\E^{n-1}$.

The converse is easy to verify by using (1.8). 
\sq
 
 Theorem 1.2 implies immediately the following.

{\bf Corollary 1.3.} {\it   Surfaces of Delaunay
(i.e., the surfaces of revolution in $\E^3$ with
constant mean curvature functions) have homothetic
LG-transforma\-tions.
\sq}

For hypersurfaces with homothetic $LG$--transformation, we
have the following.

{\bf Theorem 1.4.} {\it  Let $x: M \rightarrow
\E^{n+1}$ be an isometric immersion. Then
 the LG-transformation $LG :L (M^n ) 
\rightarrow G (M^n)$ of $x$ is homothetic if and
only if, locally, either $M$
 has  constant mean curvature function in $\E^{n+1}$ or $M$
is the product of an affine $(n-1)$--space
$\E^{n-1}$ with a planar curve whose curvature function is
given by
$$\kappa(s)={n\over 2}(\sqrt{a^2e^{-s^2}+4(c+c^2)}-ae^{-s}),
\leqno(1.15)$$ for some positive constants $a,c$.
}
 
 \demo If $M$ is a hypersurface of $\E^{n+1}$ with
homothetic $LG$--transforma\-tion, then Theorem 1.2 implies
locally either $M$ has constant mean curvature or $M$ is the
product of an affine $(n-1)$--space $\E^{n-1}$ and a planar
curve $\gamma(s)$. Moreover, from the proof of Theorem 1.2, we
know that if the second case occurs, the curvature function
$\kappa(s)$ of the planar curve satisfies
the differential equation:
$$\kappa^4+n^2{\kappa'}^2=n^2c^2\kappa^2,\leqno(1.16)$$
for some positive constant $c$. By solving differential
equation (1.16), we obtain (1.15).

The  converse can be easily verified.
\sq
 
 \vskip.2in
 \noindent{\bf \S2. $LG$-transformation of spherical 
 submanifolds.} \vskip.1in

In this section we study the $LG$--transformation of
spherical hypersurfaces.
 
{\bf Theorem 2.1.} {\it Let $x: M\rightarrow S^{n+1}\subset
\E^{n+2}$ be a hypersurface of a hypersphere $S^{n+1}$ of $
\E^{n+2}$. Then $M$ has constant mean curvature function and
homothetic LG-transformation if and only if  either 
\begin{itemize}
\item[(1)]  $M$ is  an open part of the
product of two spheres: $M^k (a)\times M^{n-k}(b)$ with some
suitable radii a and b; or
\item[(2)]  $M^n$ is an open part of a
hypersphere of $S^{n+1}.$ 
\end{itemize}
} 

\demo Assume $M$ is a hypersurface of $S^{n+1}$ of
$\E^{n+2}$ with constant mean curvature and homothetic
$LG$--transformation. Without loss of generality, we may
assume the hypersphere is centered at the origin and with
radius one. Denote by $\bar\alpha$ the mean curvature of $M$
in $S^{n+1}$ and by $\xi$  a unit normal vector of $M$ in
$S^{n+1}$. Then we have
$$H=\bar\alpha\xi-x.\leqno(2.1)$$
Let $e_1,\ldots,e_n$ be a local orthonormal frame field of $M$
such that
$$A_\xi e_i=\mu_i e_i,\quad i=1,\ldots,n.\leqno(2.2)$$
Then from (2.1) and (2.2) we have
$$dL(e_i)=n(\bar\alpha \mu_i+1)e_i,\quad
i=1,\ldots,n.\leqno(2.3)$$

From the definition of the Gauss map, we have $G(p)=\xi\wedge
x$. Thus, we have
$$dG(e_i)=-\mu_ie_i\wedge x+\xi\wedge e_i,\quad
i=1,\ldots, n.\leqno(2.4)$$
From (2.3) and (2.4) we see that the $LG$--transformation of
$M$ is homothetic if and only if 
$$(\bar\alpha\mu_i+1)^2=c^2(\mu_i^2+1), \quad
i=1,\ldots,n,\leqno(2.5)$$ 
for some positive constant $c$.
(2.5) implies that $M$ has at most two constant principal
curvatures in $S^{n+1}$. If $M$ has exactly one constant
principal curvature, $M$ is an open portion of a hypersphere
of $S^{n+1}$. If $M$ has exactly two constant principal
curvatures, then $M$ is an open portion of the product of two
spheres with some suitable radii which are completely
determined by the two roots of the quadratic equation given by
(2.5). 

The converse is easy to verify.
\sq

If $M$ has non--constant mean curvature function in $S^{n+1}$,
then we have the following result.

{\bf Theroem 2.2.} {\it  Let $ M$ be a
hypersurface of a hypersphere $S^{n+1}$ of $\E^{n+2}$. If
$M$ has non-constant mean curvature function and the
LG-transformation $LG  :L (M)  \rightarrow G
(M^n)$ is conformal, then
\begin{itemize}
\item[(1)] the gradient of the mean
curvature function is an eigenvector of the
shape operator $A_{\xi}$ of $M^n$ in $S^{n+1}$;
\item[(2)]
the shape operator $A_{\xi}$ has at most three distinct
eigenvalues; and
\item[(3)]  the eigenvalue corresponding to
the eigenvector given by the gradient of the mean
curvature function is of multiplicity one.
\end{itemize}
}

\demo Assume $M$ is a hypersurface of $S^{n+1}$ of
$\E^{n+2}$ with non--constant mean curvature and conformal
$LG$--transformation. Without loss of generality, we may
assume the hypersphere is centered at the origin and with
radius one. Denote by $\bar\alpha$ the mean curvature of $M$
in $S^{n+1}$ and by $\xi$ the a unit normal vector of $M$ in
$S^{n+1}$ as before.
And let $e_1,\ldots,e_n$ be a local orthonormal frame field of
$M$ such that
$$A_\xi e_i=\mu_ie_i,\quad i=1,\ldots,n.\leqno(2.6)$$
Then  we have
$$dL(e_i)=n(\bar\alpha
\mu_i+1)e_i-n(e_i\bar\alpha)\xi,\quad
i=1,\ldots,n.\leqno(2.7)$$

Also for the Gauss map, we have as before the following
$$dG(e_i)=-\mu_ie_i\wedge x+\xi\wedge e_i,\quad
i=1,\ldots, n.\leqno( 2.8)$$
From (2.7) and (2.8) we see that the induced metrics of the
Laplace and Gauss images satisfy
$$g_{_L}(e_i,e_j)=n^2(\bar\alpha\mu_i+1)(\bar\alpha\mu_j+1)
\delta_{ij}+n^2(e_i\bar\alpha)(e_j\bar\alpha),\leqno(2.9)$$
$$g_{_G}(e_i,e_j)=(\mu_i\mu_j+1)\delta_{ij}.\leqno(2.10)$$
From (2.9) and (2.10), we see that $\nabla\bar\alpha$ is an
eigenvector of $A_\xi$. Without loss of generality, we may
assume $\nabla\bar\alpha$ is parallel to $e_1$. Thus, from
(2.9) and (2.10), we may get
$$(\bar\alpha\mu_1+1)^2+|\nabla\bar\alpha |^2=f^2(\mu_1^2+1),
\leqno(2.11)$$
$$(\bar\alpha\mu_i+1)^2=f^2(\mu_i^2+1),\quad i=2,\ldots,n,\leqno(2.12)$$ for some positive function $f$. 
From (2.11) and (2.12) we conclude that the shape operator
$A_\xi$ has at most three distinct eigenvalues and the
eigenvalue $\mu_1$ is of multiplicity one.
\sq

For surfaces, we prove the following

{\bf Theorem 2.3.} {\it  Let $ M$ be a surface
of a hypersphere $S^{3}$ of $\E^{4}$. Then $M$ has constant
 mean curvature  and  conformal
LG-transformation if and only if  $M$ is either a totally
umbilical surface or a minimal surface in $S^3$;
}

\demo  Without loss of
generality, we may assume the radius of $S^3$ is one.

If $M$ is a totally umbilical surface in $S^3$, then
$\mu_1=\mu_2=\bar\alpha$ which is a constant. Thus, $M$ has
constant mean curvature; moreover,(2.9) and (2.10) imply that the
$LG$--transformation is homothetic. 

If $M$ is a minimal surface in $S^3$, then $\mu_1=-\mu_2$.
Thus, (2.9) implies $M$ has conformal $LG$--transformation.

Conversely, if $M$ has constant mean curvature 
and conformal $LG$--transfor\-mation, then (2.9) and (2.10) imply
$$(\bar\alpha \mu_1+1)^2=f^2(\mu_1^2+1),\leqno(2.13)$$
$$(\bar\alpha \mu_2+1)^2=f^2(\mu_2^2+1),\leqno(2.14)$$
for some positive function $f$. Combining (2.13) and (2.14) we
get
$$(\mu_1^2+1)(\bar\alpha\mu_2+1)^2=(\mu_2^2+1)(\bar\alpha
\mu_1+1)^2.\leqno(2.15)$$
Simplifying (2.15) we find
$$(\mu_1-\mu_2)\bar\alpha=0\leqno(2.16)$$
which implies  $M$ is neither a totally umbilical or a minimal
surface in $S^3$. 
\sq

Recall that a Clifford torus is the product of two plane
circles with the same radius.

{\bf Theorem 2.4.} {\it  Let $M$ be a surface of a hypersphere
$S^3$ of $\E^4$. Then $M$ has constant Gauss curvature and
homothetic $LG$--transformation if and only if $M$ is either
 a totally umbilical surface or  a Clifford torus  in $S^3$.
}

\demo Assume $M$ is a surface in a unit hypersphere $S^3$
with constant Gauss curvature and homothetic
$LG$--transformation. From Theorem 2.2 we may choose $e_1,e_2$
to be eigenvectors of $A_\xi$ such that $\nabla\bar\alpha$ is
parallel to $e_1$.  From the proof of  Theorem 2.2, we obtain
$$(\bar\alpha\mu_1+1)^2+|\nabla\bar\alpha |^2=c^2(\mu_1^2+1),
\leqno(2.17)$$
$$(\bar\alpha\mu_2+1)^2=c^2(\mu_2^2+1),\leqno(2.18)$$
for some positive constant $c$. Since the Gauss curvature $G$
of $M$ is constant, we get
$$\mu_1\mu_2+1=a,\leqno(2.19)$$
where $G=a$ is constant. 

{\bf Case 1.} If $\mu_2=0$, then  (2.18) yields $c=1$. Hence, 
by (2.17)  and $\mu_1=\bar\alpha$, we obtain $\mu_1=0$.
Thus, in this case, $M$ is a totally geodesic surface in $S^3$;
hence $M$ is an open part of a great sphere of $S^3$.

{\bf Case 2.} If $\mu_2\not=0$, then $\mu_1=\mu_2^2(a-1)$.
Hence, by (2.18), we get
$(\mu_2^2+a+1)^2=4c^2(\mu_2^2+1),$ from which we conclude that
$\mu_2$ is constant. Hence $\mu_1$ and $\bar\alpha$ are also
constant. It is known that the only surfaces in $S^3$ with
constant mean curvature and constant Gauss curvature are  open
parts of the product of two circles or open parts  of a
totally umbilical surface of $S^3$.
If $M$ is an open part of the product of two circles, then,
according to Theorem 2.3, $M$ is a minimal surface of $S^3$.
Since the only minimal surfaces of $S^3$ with constant Gauss
curvature are either open portions of a great sphere or open
portions of a Clifford torus, we conclude that $M$ is an
open part of a Clifford torus.

The converse is easy to verify.   \sq

\vfill\eject

\noindent{\bf Chapter IX: LAPLACE MAP AND 2-TYPE IMMERSIONS}
\vskip.2in

\noindent{\bf \S1.  Spherical Laplace map.} \vskip.1in 

If an isometric immersion $x: M\rightarrow \E^m$ has
nonzero constant mean curvature function  $\alpha$, then
the Laplace image $L(M )$ is contained in a
hypersphere $S^{m-1}(n\alpha)$ centered at the origin
and with radius $n\alpha$. In this case we have an
associated spherical map 
$$L_{S} :M \rightarrow  S^{m-1}(n\alpha ).\leqno(1.1)$$
 We call this spherical map the {\it spherical Laplace
map\/} of the immersion. 

In this section we study the following problem.

{\bf Problem 1.1} {\it  ``When is the spherical Laplace map
$L_{S} :M^n \rightarrow S^{m-1}$ of an isometric immersion
harmonic?''}

First we give the following result.

{\bf Theorem 1.1.} {\it  Let $x :M \rightarrow  \E^m$
be an isometric spherical immersion from an $n$--dimensional
Riemannian manifold $M$ into  $\E^m$. If $x$ has   constant mean
curvature function $\alpha$, then
\begin{itemize}

\item[(1)] the associated spherical Laplace map 
$L_{S} :M\rightarrow S^{m-1}(n\alpha )$ of $x$ has
positive constant energy density;  

\item[(2)] the spherical Laplace map $L_S$ of
$x$ is a harmonic map if and only if $M$ is
immersed as a minimal submanifold of a hypersphere of
$\E^m$ via $x$.
\end{itemize}
}

\demo Let $x :M \rightarrow \E^m$
be an isometric spherical immersion with nonzero constant mean
curvature function $\alpha$. Without loss of generality, we may
assumed $M$ is immersed in the hypersphere $S^{m-1}(r)$ centered at
the origin and with radius $r$. Since a spherical submanifold
in $\E^m$ has nonzero mean curvature function,  by the constancy of
the mean curvature,  we
may assume $n\alpha=1$ for simplicity. In this case, the Laplace map
is the composition $j\circ L_S$ of the spherical Laplace map
$L_S$ followed by the inclusion $j$ of the unit hypersphere
$S^{m-1}(1)$ in $\E^m$. The second fundamental forms $h_L,
h_{L_S}$ and $h_j$ of the maps $L, L_S$ and $j$,
respectively, satisfy
$$h_L(X,Y)=j_*h_{L_S}(X,Y)+h_j(dL_S(X),dL_S (Y))\leqno(1.2)$$
for $X,Y$ tangent to $M$. Thus we have
$$\Delta L=-\tau(L)=-j_*\tau(L_S)-\sum_{i=1}^n
h_j(dL_S(e_i), dL_S(e_i)),\leqno(1.3)$$
where $e_1,\ldots,e_n$ is an orthonormal local frame field of
$M$, and $\tau(L),\tau(L_S)$ are the tension fields of the
Laplace map and the spherical Laplace map, respectively.
Since $j$ is a totally umbilical isometric immersion, (1.3)
yields
$$\Delta L=-j_*\tau(L_S)+2e(L_S)L,\leqno(1.4)$$
where $e(L_S)$ is the energy density of the spherical
Laplace map. 

On the other hand, since $M$ is isometrically immersed into the
hypersphere $S^{m-1}(r)$ and the immersion $x$ has constant mean
curvature, the mean curvature of $x$ is a positive constant.
Moreover, from Theorem 1.3 of Chapter II, we have
$$\Delta H=\Delta^{\bar D}\bar H +(n+ ||A_\xi ||^2)\bar H +{\cal
A}(\bar H) +2\hbox{trace}\,A_{DH}-{{n\alpha^2}\over {r^2}}x,\leqno(1.5)$$
where $\bar H$ is the mean curvature vector of $M$ in $S^{m-1}(r)$, 
$\xi$  the unit vector in the direction of $\bar H$, and $\bar D$,
${\cal A}(\bar H)$
the normal connection and  the allied mean curvature vector of $M$ in
$S^{m-1}(r)$, respectively. 

Because $L=-nH=-n\bar H+{n\over{r}}x$, by comparing the
$x$-components of (1.4) and (1.5), we obtain
$e(L_S)={{2\alpha^2}\over 
{r^2}}$ which is a nonzero constant. This proves (1).

(2) If the spherical Laplace map $L_S$ is harmonic, then (1.4)
implies
$\Delta H=cH,$ where $c=2e(L_S)$ is a constant by statement (1). 
Thus, by Theorem 3.11 of Chapter II, 
$M$ is either immersed
as a minimal submanifold of a hyperphere of $E^m$ via $x$, or immersed
as a biharmonic submanifold of $E^m$, or immersed as a null 2--type
submanifold. Because $x$ is spherical, the second and the third cases
cannot occur (cf. [C16]).  
\sq

As  an immediate consequence of Theorem 1.1, we have the
following

{\bf Corollary 1.2.} {\it  Let $x: M^n \rightarrow
\E^{n+1}$ be an isometric immersion with nonzero constant
mean curvature. Then the spherical Laplace map $L_S$ is a
harmonic map if and only if $M^n$ is an open part of a
hypersphere of $\E^{n+1}$. \sq}

For curves we have the following.

{\bf Proposition 1.3.}  {\it Let $\beta (s)$ be a curve in
$\E^m$ whose  first Frenet curvature function is a nonzero
constant.
 Then  the spherical Laplace map of $\beta$ is a
harmonic map if and only either $\beta$ is an open part of a circle
in a plane or $\beta$ is an open part of a circular helix lying
in an affine 3--space of $\E^m$. }

\demo Without loss of generality, we may assume $\beta$ is a
unit speed curve in $\E^m$ with $m\geq 4$. Denote by $\kappa_i$ and
$\beta_i$ the $i$--the Frenent curvature and the $i$--the Frenet
vector of $\beta$. By the hypothesis and direct computation, we have
$$\tau(L_S)=\kappa_1(\kappa_2\beta_3+\kappa_2\kappa_3\beta_4).\leqno(
1.6)$$
From this we see that the second Frenent curvature $\kappa_2$ is a
constant. If $\kappa_2=0$, $\beta$ is an open part of a circle in a
plane. And if $\kappa_2\not=0$, then $\kappa_3=0$. In this case,
$\beta$ is an open part of a circular helix. \sq

\eject

\noindent{\bf \S2.  2--type immersions.} \vskip.1in

The main purpose of this section is to give some special
properties of 2--type immersions. 

Let $x:M\rightarrow \E^m$ be an isometric immersion of
$k$--type with spectral decomposition given by 
$$x = c + x_1 +\cdots+ x_k,\quad \Delta x_1 = \lambda_1
x_1,\ldots \Delta x_k = \lambda_k x_k,\quad  \lambda_1
<\cdots< \lambda_k,$$ where $c$ is a constant map and
$x_1,\ldots,x_k$ are non--constant maps. Put
$$E_i=\hbox{Span}\{x_i(p):p\in M\}.$$
Recall that the immersion $x$ is said to be {\it
orthogonal\/} if the subspaces $$E_1,\ldots,E_k$$ are
mutually orthogonal in $\E^m$.

For  orthogonal 2--type immersions, we have the following
result.

{\bf Propostion 2.1.} {\it Let $x:M\rightarrow \E^m$ be 
an isometric immersion of a compact manifold $M$ into
$\E^m$. Assume $x$ is an orthogonal 2--type immersion with
spectral decomposition given by
$$x=x_1+x_2,\quad \Delta x_1 = \lambda_1 x_1,\quad
\Delta x_2 = \lambda_2 x_2.\leqno(2.1)$$
Then
\begin{itemize}
\item[(1)] $M$ lies in a hypersphere $S^{m-1}(r)$
centered at the origin with radius, say $r$, as  a
mass--symmetric submanifold;
\item[(2)]  $M$ has constant mean curvature $\alpha$ in $\E^m$
given by
$$\alpha^2={1\over n}(\lambda_1+\lambda_2)-\Big({r\over
n}\Big)^2\lambda_1\lambda_2;$$
 \item[(3)] both $x_1,x_2$ are mass--symmetric spherical
maps; 
 \item[(4)] $x_1,x_2$ and the Laplace map of $x$ satisfy
$$x_1={{\lambda_2 x-L}\over{\lambda_2-\lambda_1}},\quad
x_2={{\lambda_1x-L}\over{\lambda_1-\lambda_2}};\leqno(2.2)$$
and
\item[(5)] the  spherical maps $\bar x_1,\bar x_2$,
induced from $x_1,x_2$, are harmonic maps which have
constant energy density. \end{itemize}
}

\demo Let $x:M\rightarrow \E^m$ be 
an orthogonal, 2--type, isometric immersion of a compact
manifold $M$ into $\E^m$ whose  spectral decomposition
is given by (2.1).
For any vector $X$ tangent to $M$, we put
$$X=X_1+X_2,\quad X_1\in E_1,\quad X_2\in E_2.\leqno(2.3)$$
Then, by the orthogonality of $E_1,E_2$, we have
$$\tilde\nabla_X x_i=X_i,\quad i=1,2.\leqno(2.4)$$

We put
$$f=\<\Delta x,x\>.\leqno(2.5)$$
Then, from (2.5), Beltrami's formula and  the hypothesis,   we
have 
$$Xf=X(\lambda_1\<x_1,x_1\>+\lambda_2\<x_2,x_2\>)$$
$$=2(\lambda_1\<x_1,X_1\>+\lambda_2\<x_2,X_2\>)$$
$$=2\<\Delta x,X\>=0.$$
Therefore, $f$ is a constant. 
Hence,
$$\Delta\<x,x\>=2\<\Delta x,x\>-2n=2f-2n\leqno(2.6)$$
is a constant.
Because $M$ is compact, this implies that  $\<x,x\>$ is a
constant. Thus, $M$ is contained in a hypersphere
$S^{m-1}(r)$ of $\E^m$ centered at the origin and with
radius, say $r$. This shows that $M$ is a mass--symmetric, 2--type
spherical submanifold which proves (1).

(2) Follows from (1) and a Theorem 4.1 of [C5, page 274].

(3) From (1) we have
$$\<x_1,x_1\>+\<x_2,x_2\>=r^2.\leqno(2.7)$$
On the other hand, by using Beltrami's formula, we also
have
$$\lambda_1\<x_1,x_1\>+\lambda_2\<x_2,x_2\>=\<-nH,x\>=n.\leqno(2.8)$$
Combining (2.7) and (2.8) we obtain 
$$\<x_1,x_1\>={{\lambda_2r^2-n}\over{\lambda_2-\lambda_1}},
\quad
\<x_2,x_2\>={{\lambda_1r^2-n}\over{\lambda_1-\lambda_2}},
\leqno(2.9)$$
which implies (3).

(4) Follows from (2.1) and the fact: $L=\lambda_1x_1+
\lambda_2x_2.$

(5) From (3) we have
$$x_1: M\rightarrow S^{m-1}(r_1)\subset \E^m,\leqno(2.10)$$
where 
$$r_1=\sqrt{(\lambda_2r^2-n)/(\lambda_2-\lambda_1)}.$$
Denote by $h_j$ the second fundamental form of
$S^{m-1}(r_1)$ in $\E^m$. Then (2.10) implies that the
tension fields of $x_1$ is given by
$$\tau(x_1)=j_*\tau(\bar x_1)-2e(\bar x_1)x_1,\leqno(2.11)$$
where $\bar x_1$ is the map $M\rightarrow S^{m-1}(r_1)$
induced from $x_1:M\rightarrow \E^m$, $j$ the inclusion of
$S^{m-1}(r_1)$ into $\E^m$, and $e(\bar x_1)$ the energy
density of $\bar x_1$.

On the other hand, we have 
$$\Delta x_1=-\tau(x_1)=\lambda_1 x_1.\leqno(2.12)$$
Comparing (2.11) and (2.12), we conclude that $\bar x_1$
is a harmonic map and it has constant energy density.
The same argument can be applied to $x_2$.
\sq

{\bf Remark 2.1.} The condition of compactness in
Proposition 2.1 is essential. For instance, the circular
cylinder defined by 
$$x(t,s)=(t,\cos s,\sin s)$$
is an orthogonal 2--type surface in $\E^3$ which is not
spherical. \sq

If  $x : M \rightarrow E^m$ is a map of 2-type,  we 
define the notion of {\it conjugation\/} of  $x$  as
follows. 

{\bf Definition 2.1.} Let $x$ be a map of 2-type whose 
 spectral decomposition takes  the following form: 
$$x = c + x_1 + x_2,\quad \Delta x_1 = \lambda_1 x_1,\quad
\Delta x_2 = \lambda_2 x_2,\quad  \lambda_1 < \lambda_2,
$$
where $c$ is a constant map and $x_1$ and $x_2$ are
non--constant maps. Then the map ${\bar
x}: M^n \rightarrow \E^m$ given by 
$${\bar x} =c +x_1-x_2\leqno(2.13)$$
 is called the {\it conjugate\/} of the 2--type map $x$. \sq

The conjugate of a 2--type isometric immersion is not
isometric, in general. For example, for the conjugate of
a 2--type curve in $\E^m$, we have the following result.

{\bf Proposition 2.2.} {\it The conjugate $\bar\beta$ of a
2-type unit speed curve $\beta(s)$ in $\E^m$ is   a unit speed
curve  if and only if $\beta$ is either an open part of a
circular helix lying in an affine 3--space or an open part of
the diagonal immersion of two circles. }

\demo First we recall that if $\beta(s)$ is curve of finite
type, then $\beta$ can be written as (cf. [CDVV1])
$$\beta(s)=a_0+b_0s+\sum_{t=1}^k(a_t\cos(p_ts)+b_t\sin(p_ts)).\leqno(2.14)$$
We may assume $a_0=0$ by choosing $a_0$  as the origin of $\E^m$. 
Thus, if  $\beta$ is of 2--type, then $\beta$ can be espressed as
one of the following forms:
$$\beta(s)=b_0s+(a_1\cos(ps)+b_1\sin(ps)),
\leqno(2.15)$$
or 
$$\beta(s)=\sum_{t=1}^2(a_t\cos(p_ts)+b_t\sin(p_ts)),
\leqno(2.16)$$
where $b_0,a_1,a_2,b_1,b_2$ are constant vectors in $\E^m$ and
$p,p_1,p_2$ are nonzero constants.

If $\beta(s)$ is   a unit speed curve whose conjugate is also a
unit speed curve, then we have
$$\<\beta'(s),\beta'(s)\>=\<\bar\beta'(s),\bar\beta'(s)\>=1.
\leqno(2.17)$$

{\bf Case 1.} If (2.15) holds, then, by (2.17), we may
prove that $b_0,a_1,b_1$ are mutually orthogonal and
$|a_1|=|b_1|$. Hence, in this case, $\beta$ is an open
part of a circular helix.
 
{\bf Case 2.} If (2.16) holds, then (2.17) implies that
$a_1,a_2,b_1,b_2$ are mutually orthogonal and,
moreover, $|a_1|=|b_1|, |a_2|=|b_2|$. Hence, $\beta$ is the
diagonal immersion of two circles.

The converse is easy to verify.
\sq

Similar to the notion of orthogonal immersions, we give the
following

{\bf Definition 2.2.} Let $x: M\rightarrow \E^m$ be an immersion of
$k$--type whose spectral decomposition is given by
$$x = c + x_1 +\cdots+ x_k,\quad \Delta x_1 = \lambda_1 x_1,\ldots
\Delta x_k = \lambda_k x_k,\quad  \lambda_1 <\cdots< \lambda_k,$$
where $c$ is a constant map and $x_1,\ldots,x_k$ are
non--constant maps. Then the immersion $x$ is said to be {\it
pointwise orthogonal\/} (respectively, {\it strongly pointwise
orthogonal\/}) if for each $p\in M$, $x_1(p),\ldots,x_k(p)$ are
mutually orghogonal (respectively,
the  image subspaces ${x_1}_*
(T_pM),\ldots,{x_k}_* (T_p M)$ are mutually orthogonal) in $\E^m$.
\sq

For  2-type isometric immersions, we have the following result.

{\bf Theorem 2.3.} {\it  Let $x : M \rightarrow \E^m$  be a 2-type
isometric immersion, then
\begin{itemize}
\item[(1)]  $x$ is  strongly pointwise  orthogonal if and
only if  the conjugate ${\bar x}$ of $x$ is an
isometric immersion;
\item[(2)]  $x$ is pointwise orthogonal if the isometric
immersion  $x$ is a  mass-symmetric spherical immersion;
and \item[(3)]   $x$ is orthogonal
 if and only if the conjugate ${\bar x}$ of $x$ is congruent to
 immersion $x$ up to an orthogonal transformation of $\E^m$.
 \end{itemize}
}  

\demo (1) Let $x$ be a 2--type isometric immersion whose
conjugate $\bar x$ is also an isometric immersion. Assume the
spectral decomposition of $x$ is given by (2.13). Then we
have
$$\<dx,dx\>=\<dx_1,dx_1\>+\<dx_2,dx_2\>+2\<dx_1,dx_2\>,\leqno(2.18)$$
$$\<d\bar x,d\bar
x\>=\<dx_1,dx_1\>+\<dx_2,dx_2\>-2\<dx_1,dx_2\>.\leqno(2.19)$$ Since
both $x$ and $\bar x$ are isometric immersions, (2.18) and
(2.19) yield $$\<dx_1,dx_2\>=0$$ identically. This implies
that the immersion $x$ is strongly pointwise orthogonal. 

Conversely, if $x$ is strongly pointwise orthogonal, then (2.18)
and (2.19) imply $\<dx,dx\>=\<d\bar x,d\bar x\>$. Therefore, the
conjugate of $x$ is isometric.

(2) Let $x:M\rightarrow S^{m-1}(r)\subset \E^m$ be a spherical,
mass--symmetric, 2--type isometric immersion. Assume the
spectral decomposition of $x$ is given by
$$x =  x_1 + x_2,\quad \Delta x_1 = \lambda_1 x_1,\quad
\Delta x_2 = \lambda_2 x_2,\quad  \lambda_1 < \lambda_2,$$
Then we have
$$\<x_1,x_1\>+\<x_2,x_2\>+2\<x_1,x_2\>=r^2.\leqno(2.20)$$
By  Beltrami's formula, we have 
$$\<\lambda_1x_1+\lambda_2x_2,x\>=\<\Delta x,x\>=\<-nH,x\>=n.\leqno(2.21)$$
Furthermore, by Lemma 4.2 of [C5, p.273], we have
$$\<\lambda_1^2x_1+\lambda_2^2x_2,x\>=-n\<\Delta
H,x\>=({{n\alpha}\over {r}})^2.\leqno(2.22)$$
From (2.20), (2.21) and (2.22) we obtain $\<x_1,x_2\>=0$. Thus,
$x$ is pointwise orthogonal.

(3) Let $x$ be a 2--type immersion. Without loss of generality,
we may assume the spectral decomposition of $x$ is given by
(2.1). Assume, up to orthogonal transformations of $E^m$,
 the conjugate $\bar x$ of $x$ is congruent to the immersion $x$. Then
there  is  $A\in O(m)$ such that $Ax=\bar x$. Thus, we have
$$(Ax_1-x_1)+(Ax_2+x_2)=0.\leqno(2.23)$$
Because $\Delta(Ax_i)=A(\Delta x_i)=\lambda_iAx_i$, (2.23)
yields $$\lambda_1(Ax_1-x_1)+\lambda_2(Ax_2+x_2)=0.\leqno(2.24)$$
From (2.23) and (2.24) we get
$$Ax_1=x_1,\quad Ax_2=-x_2.\leqno(2.25)$$
Let $E_i=\hbox{Span}\{x_i(p):p\in M\}, i=1,2,$ and let $V_1,V_2$
be the eigenspaces of $A$ with eigenvalues $1,-1$, respectively.
Then (2.25) implies
$$E_1\subset V_1,\quad E_2\subset V_2.\leqno(2.26)$$
Because $A$ is an orthogonal transformation, (2.26) implies that,
for any $u\in V_1, v\in V_2$, we have
$\<u,v\>=\<Au,Av\>=-\<u,v\>.$ Therefore, $V_1$ and $V_2$ are
orthogonal in $\E^m$. Consequently, by (2.26), we conclude that 
$x$ is an orthogonal immersion.

The converse of this is clear.
\sq

\eject
\noindent{\bf \S3. Laplace map of 2--type immersions.} \vskip.1in

The purpose of this section is to study relations between the Laplace
maps of a 2--type immersion and its conjugate.

{\bf Lemma 3.1.} {\it  Let $x: M\rightarrow
\E^m$ be a 2-type isometric immersion. Then, up
to $\pm 1$, the
Laplace map of $x$ and the Laplace map of the
conjugate ${\bar x}$ of $x$ in $\E^m$  are equal  if and only if
the immersion $x$  is of null 2-type.}

\demo Let $x$ be a 2--type immersion with spectral
decomposition given by 
$$x=c+x_1+x_2,\quad \Delta x_1 = \lambda_1 x_1,\quad
\Delta x_2 = \lambda_2 x_2.\leqno(3.1$$
Then the Laplace map  of $x$ and its
conjugate $\bar x$ of $x$ are given by
$$L=\lambda_1x_1+\lambda_2x_2,\quad \bar
L=\lambda_1x_1-\lambda_2x_2.\leqno(3.2$$
Therefore, $L=\bar L$ (respectively, $L=-\bar L$) if and only if
$\lambda_2=0$ (respectively, $\lambda_1=0$). Hence, $L=\bar L$
if and only if  immersion $x$ is of null 2--type. \sq

{\bf Lemma 3.2.} {\it  If  $x: M \rightarrow \E^m$ is a 2-type
isometric immersion whose conjugate ${\bar x}$ is an  isometric
immersion, then the Laplace map of the conjugate ${\bar x}$ is the
conjugate of the Laplace map.
}

\demo If  $x: M \rightarrow \E^m$ is a 2-type
isometric immersion whose conjugate ${\bar x}$ is an  isometric
immersion, then the Laplace operator $\bar\Delta$ on $M$
induced from the conjugate $\bar x$ equals  the Laplace
operator $\Delta$ of $M$. Hence, 
the Laplace map of $\bar x$ is given by $\lambda_1 x_1-\lambda_2
x_2$ which is nothing but the conjugate of the Laplace map of
$x$.
\sq

We now introduce the following
\vskip.1in
{\bf Definition 3.1.} A 2--type immersion (respectively, a
2--type map) $x:M\rightarrow \E^m$ is called a {\it
dual 2-type immersion\/} (respectively, a dual 2--type map) if
the spectral decomposition of immersion $x$ takes the form:
$x  =  c+x_1 + x_2$ such that $\Delta x_1 =-\lambda x_1$ and
$\Delta x_2 =\lambda x_2$, where $c$ is a constant map and
$x_1$ and $x_2$ are two nonconstant maps and 
 $\lambda$ is a positive  number. \sq
\vskip.1in
 Geometrically, up to a constant, a dual 2-type immersion is
an immersion whose the Laplace transformation is
involutive.

{\bf Lemma 3.3.} {\it  Let $x: M\rightarrow \E^m$ be a 2-type
isometric immersion. Then, up to nonzero constants, the Laplace
map  of $x$ is the conjugate $\bar x$ of $x$  if and only if
immersion $x$ is of dual 2-type. }

\demo Let $x$ be a 2--type immersion with spectral
decomposition given by (3.1). If there is a  constant
$\mu\not=0$ such that the Laplace map $L$ of $x$ and  the
conjugate $\bar x$ of $x$ satisfy $L=\mu\bar x$, then
 we have
$$(\lambda_1-\mu)x_1+(\lambda_2+\mu)x_2=\mu  c.\leqno(3.3)$$
From (3.3), we get
$$\lambda_1(\lambda_1-\mu)x_1+\lambda_2(\lambda_2+\mu)x_2=0,
\leqno(3.4)$$
$$\lambda_1^2(\lambda_1-\mu)x_1+\lambda_2^2(\lambda_2+\mu)x_2=0.
\leqno(3.5)$$
Since $\lambda_1\not=\lambda_2$, (3.3), (3.4) and (3.5) 
imply
  $$c=0,\quad \lambda_1=-\lambda_2.\leqno(3.6)$$
Thus, the immersion $x$ is of dual 2--type. 

Conversely, if $x$ is of dual 2--type, then
$$ x= c+x_1+x_2,\quad \Delta x_1=-\lambda x_1,\quad \Delta
x_2=\lambda x_2.\leqno( 3.7)$$
From (3.7), we get
$L=-\lambda(x_1-x_2)=\lambda\bar x.$
\sq

{\bf Remark 3.1.} For an isometric immersion $x : M \rightarrow
\E^m$  of a Riemannian manifold $(M^n, g)$ into $\E^m$ and for
 a positive real number $c$, we have  the immersion $x^c$
defined by $x^c =cx$. The induced
metric $g^c$ on $M$ with respect to $x^c$ is given by
$g^c = c^2 g$ and the Laplace operator $\Delta^c$ of
$(M ,g^c )$ is given by $\Delta^c = c^{-2}\Delta$. It
is easy to see that the Laplace image  of $cx$
is given by $c\, L (M)$. Therefore, when the
Laplace transformation ${\cal L} : M^n \rightarrow L(M )
$ is homothetic,  one may multiply a suitable constant to
the immersion $x$ to obtain an isometric immersion whose
Laplace transformation is also an isometry. \sq

In terms of conjugation, we also have the following  result.

{\bf Lemma 3.4.} {\it  Let  $x : M\rightarrow \E^m$  be a dual 2-type
isometric immersion. If the Laplace transformation $\cal L$
of $x$ is isometric, then 
\begin{itemize}
\item[(1)] the immersion $x$ is constructed from eigenspaces of
$\Delta$ belonging to eigenvalues $1$ and $-1$;
\item[(2)] the immersion $x$ is strongly pointwise orthogonal;
\item[(3)] up to sign, the Laplace map $L : M\rightarrow \E^m$ is
an immersion which is given by the conjugation of the immersion
$x : M \rightarrow \E^m$; 
\item[(4)] the immersion $x : M \rightarrow
\E^m$  is the conjugate of the Laplace immersion  $L : M
\rightarrow \E^m$ of $x$; and 
\item[(5)] the Laplace transformation ${\cal L}$ is
idempotent.
\end{itemize}
}

\demo  Let  $x= c+x_1+x_2$  be a dual 2-type
isometric immersion. Then $\Delta x_1=\lambda x_1,\,\Delta x_2=
-\lambda x_2$, for some  number $\lambda\not=0.$ Thus, the
Laplace map is given by 
$$L=\lambda (x_1-x_2).\leqno(3.8)$$
Assume the Laplace transformation ${\cal L}:M\rightarrow L(M)$
is isometric. Then we have
$$\aligned\lambda^2(\<dx_1,dx_1\>-2\<dx_1,dx_2\>+\<dx_2,dx_2\>\\
=\<dx_1,dx_1\>+2\<dx_1,dx_2\>+\<dx_2,dx_2\>.\endaligned\leqno(3.9)$$
Let $X$ be a tangent vector of $M$ such that $(x_1)_*X\not=0$ and 
$(x_2)_*X=0.$ Then from (3.9), we find
$\lambda^2=1$. This proves (1).

The remaining parts of this Lemma follow easily from (1).  
\sq

\vfill\eject

\end{document}